\documentclass[11pt,letterpaper,leqno]{amsart}
\usepackage{amsmath,amssymb,mathtools}
\usepackage[hidelinks]{hyperref}
\usepackage{comment}
\usepackage{mathrsfs}
\usepackage{xcolor}
\colorlet{sectionrevision}{red}
\hypersetup{pdftitle={Uniformly bounded real spherical harmonic bases in dimensions three and higher},
  pdfsubject={Full eigenspaces of the round sphere},
  pdfauthor={Miao He, Penghui Wang, Fan Zhang and Zeyou Zhu},
  pdfkeywords={Spherical harmonics, uniformly bounded orthonormal bases, Gaussian beams}}
\numberwithin{equation}{section}
\newtheorem{theorem}{Theorem}[section]
\newtheorem{lemma}[theorem]{Lemma}
\newtheorem{proposition}[theorem]{Proposition}

\theoremstyle{remark}

\title{Uniformly bounded real spherical harmonic bases in dimensions three and higher}
\author{Miao He}
\author{Penghui Wang}
\author{Fan Zhang}
\author{Zeyou Zhu}
\date{}

\begin{document}
\begin{abstract}
	Let $\mathcal H_n(S^d)$ be the space of real spherical harmonics of
	degree $n$ on $S^d$. Write $D_n=\dim\mathcal H_n(S^d)$. 
	Fix $d\geq 3$. 
	For each $n\geq 0$, we construct an orthonormal basis
	$\{Y_{n,j}\}_{j=1}^{D_n}$ of $\mathcal H_n(S^d)$ so that
	\[
	\sup_{n\geq 0}\max_{1\leq j\leq D_n}
	\|Y_{n,j}\|_{L^\infty(S^d)}\leq C_d,
	\]
	where $C_d$ depends only on $d$.
\end{abstract}
\maketitle


\section{Introduction}

Let $\sigma$ be the normalized surface measure on the unit sphere
$S^d\subset\mathbb R^{d+1}$. For $n\geq0$, let $\mathcal H_n(S^d)$
be the space of restrictions to $S^d$ of real homogeneous harmonic
polynomials of degree $n$ and put $D_n=\dim\mathcal H_n(S^d)$.
The space $\mathcal H_n(S^d)$ is the eigenspace of $-\Delta_{S^d}$
with eigenvalue $n(n+d-1)$ and
\[
L^2(S^d,\sigma;\mathbb R)
=\bigoplus_{n\geq0}\mathcal H_n(S^d).
\]

Toth and Zelditch~\cite{TZ} asked
whether $S^2$ admits a uniformly bounded orthonormal basis of Laplace
eigenfunctions. In 2001, Jakobson, Nadirashvili and
Toth~\cite[Section~6.1]{JNT} recorded the same question for the
standard sphere $S^d$. In the real setting, this amounts to choosing
an orthonormal basis of each $\mathcal H_n(S^d)$ with a common
supremum norm bound depending only on $d$.

VanderKam~\cite{VK} obtained an upper bound of order $(\log n)^2$
for the supremum norms of random $L^2$-normalized spherical harmonics
of degree $n$ on $S^2$.
For individual spherical harmonics on $S^2$, Han~\cite[Theorem~1.2]{Han2} constructed one uniformly bounded,
$L^2$-normalized spherical harmonic in each degree on $S^2$,
under separation and counting conditions on the poles of Gaussian beams.
For finite $p$, Burq and Lebeau~\cite{BL} proved that, for each
$d\geq2$ and each fixed $2\leq p<\infty$, almost every orthonormal
basis of spherical harmonics on $S^d$ satisfies
\[
\sup_{n\geq0}\max_{1\leq j\leq D_n}
\|Y_{n,j}\|_{L^p(S^d,\sigma)}<\infty.
\]
The bound may depend on $p$.
They also showed that, for almost every such basis, there is a
constant $C$ such that, for every $n\geq 0$,
\[
\max_{1\leq j\leq D_n}
\|Y_{n,j}\|_{L^\infty(S^d)}
\leq C\sqrt{\log(n+2)}.
\]

Ryll and Wojtaszczyk~\cite{RyllWojtaszczyk} constructed, in each
degree, a single $L^2$-normalized homogeneous holomorphic polynomial
with a uniform supremum norm bound on the complex unit sphere.
Shiffman~\cite[Theorem~1.1]{Shiffman} obtained uniformly bounded
orthonormal systems of such polynomials whose cardinalities are at
least a fixed positive fraction of the dimension in each degree.
Both constructions yield uniformly bounded real orthonormal systems
of spherical harmonics in each positive degree on odd-dimensional
spheres, by taking real and imaginary parts and normalizing.
Shiffman noted that the existence of a
uniformly bounded sequence of $L^2$-normalized real spherical
harmonics, one in each degree, was open at the time on
even-dimensional spheres.

For $k\geq1$, a related question of Rudin~\cite{Rudin} concerns
uniformly bounded orthonormal bases of homogeneous holomorphic
polynomials on the unit sphere in $\mathbb C^{k+1}$.
Bourgain answered the question affirmatively for $k=1$ and $k=2$
in \cite{Bourgain85} and \cite{Bourgain}, respectively.
The general case was settled in \cite{GFR,WZZ}.
The real and imaginary parts of the homogeneous holomorphic
polynomials of degree $n$ span a subspace $V_n$ of
$\mathcal H_n(S^{2k+1})$. The solution to Rudin's problem gives
uniformly bounded real orthonormal bases for $V_n$.

Our main theorem gives the uniform bound for the full spaces
$\mathcal H_n(S^d)$ when $d\geq3$.

\begin{theorem}\label{thm:main}
	For each $d\geq3$, there exists a constant $C_d>0$ such that,
	for every $n\geq0$, the space $\mathcal H_n(S^d)$ has a real
	$L^2(\sigma)$-orthonormal basis $\{Y_{n,j}\}_{j=1}^{D_n}$
	satisfying
	\[
	\max_{1\leq j\leq D_n}\|Y_{n,j}\|_{L^\infty(S^d)}\leq C_d.
	\]
\end{theorem}

The paper is organized as follows.
Section~2 gives a local approximation for the inner products of
Gaussian beams. In Section~3, we use this approximation to construct
local orthonormal systems. We combine these systems in Section~4 to
obtain a global orthonormal system with at least
$(1-C_dn^{-\frac{9}{100}})D_n$ functions. In Section~5, we rotate groups of functions and pair them with their
conjugates. After orthonormalization, the sum of their squared
absolute values is at least $(1-C_dn^{-\frac{1}{20}})D_n$ at every point.
Section~6 gives a real orthonormal basis of their span with a bound
on the sum of its squared supremum norms. We complete this basis
in Section~7 and prove Theorem~\ref{thm:main}.

\section{Great circles and Gaussian beams}
Let $d\ge 2$.
For \(n \geq 0\), let $\mathcal P_n$ denote the space of homogeneous polynomials of degree $n$ on $\mathbb R^{d+1}$ with complex coefficients and let $\mathcal H_n=\mathcal H_n(S^d)\otimes_{\mathbb R}\mathbb C$.
On $H_n$ we use the inner product
\[
\langle u,v\rangle
=\int_{S^d}u(x)\overline{v(x)}\,d\sigma(x),
\]
which is linear in the first variable. Let $\Delta$ denote the
Euclidean Laplacian on $\mathbb R^{d+1}$.
For $n\ge 2$,
\[
\mathcal H_n \cong \ker(\Delta:\mathcal P_n\to\mathcal P_{n-2}).
\]
According to \cite[Theorem~A.1]{HelzelTzavaras2017}, we have $\Delta:\mathcal P_n\to\mathcal P_{n-2}$ is surjective. 
Therefore, by  the
rank--nullity theorem \cite[Theorem~3.21]{Axler2024} and Pascal's identity,

\begin{equation}\label{Dn}
	\begin{aligned}
		D_n &=\dim \mathcal H_n
		=\dim\ker\Delta\\
		&=\dim\mathcal P_n-\dim\mathcal P_{n-2}\\
		&=\binom{n+d}{d}-\binom{n+d-2}{d}\\
		&=\left(\binom{n+d}{d}-\binom{n+d-1}{d}\right)
		+\left(\binom{n+d-1}{d}-\binom{n+d-2}{d}\right)\\
		&=\binom{n+d-1}{d-1}+\binom{n+d-2}{d-1}\\
		&=\frac{2n+d-1}{(d-1)!}\prod_{j=1}^{d-2}(n+j)
		\asymp_d n^{d-1},
	\end{aligned}
\end{equation}
where $D_n\asymp_d n^{d-1}$ means that there exist positive constants
$c_d$ and $C_d$, depending only on $d$, such that \(c_d n^{d-1}\le D_n\le C_d n^{d-1}\) for all $n\geq 1$.
For $n=0,1$, we have $\mathcal H_n=\mathcal P_n|_{S^d}$.
Thus $D_0=1$ and $D_1=d+1$, in agreement with \eqref{Dn}.
For $z,w\in\mathbb C^{d+1}$, write \(z\cdot w=\sum\limits_{j=1}^{d+1}z_jw_j\).
For linearly independent vectors $e,f\in\mathbb R^{d+1}$, write
$\operatorname{span}^{+}(e,f)$ for the plane spanned by $e$ and $f$,
with the orientation determined by the ordered basis $(e,f)$.

\begin{lemma}\label{lem:grassmannian-quadric}
	Let $Q=G_2^+(\mathbb R^{d+1})$ be the Grassmannian of oriented two-dimensional subspaces of $\mathbb R^{d+1}$. The map
	\begin{align*}
	Q &\longrightarrow \{[z]\in\mathbb{CP}^{d}:z\cdot z=0\}\\
	E=\operatorname{span}^{+}(e,f)
	&\longmapsto
	\left[\frac{e+if}{\sqrt2}\right],
	\end{align*}
	where $(e,f)$ is an oriented orthonormal basis of $E$, is a bijection.
\end{lemma}

\begin{proof}
	For $E\in Q$, choose an oriented orthonormal basis $(e,f)$ and set
	$z=\frac{e+if}{\sqrt2}$. 
	It is easy to verify that  \(z\cdot z = 0\) and \(z\cdot\bar z=1\).
	If $(e',f')$ is another oriented orthonormal basis of $E$, then there exists $\theta\in\mathbb R$ such that $e'+if'=e^{i\theta}(e+if)$, and hence
	$\left[\frac{e'+if'}{\sqrt2}\right]=[z]$. This proves that the map is well defined.
	
	Suppose that
	\[
	\left[\frac{e+if}{\sqrt2}\right]
	=
	\left[\frac{e'+if'}{\sqrt2}\right].
	\]
	Since both representatives have norm one, there exists
	$\theta\in\mathbb R$ such that \(e'=\cos\theta\,e-\sin\theta\,f\) and \(f'=\sin\theta\,e+\cos\theta\,f\).
	Hence $\operatorname{span}^{+}(e',f')=\operatorname{span}^{+}(e,f)$,
	which proves injectivity.
	Let $[z]\in\mathbb{CP}^{d}$ satisfy $z\cdot z=0$ and choose a unit representative $z$. 
	From $z\cdot z=0$ and $z\cdot\bar z=1$, we obtain \(\operatorname{span}^{+}(\sqrt2\,\operatorname{Re}z,\sqrt2\,\operatorname{Im}z) \in Q\).
	Obviously, 
	\[
	\frac{\sqrt2\,\operatorname{Re}z+i\sqrt2\,\operatorname{Im}z}{\sqrt2}=z.
	\]
	Hence, the map is surjective.
\end{proof}

For $E=\operatorname{span}^{+}(e,f)$, set \(\omega_E=e\wedge f\) and \(s_E(x)=|(I-P_E)x|,\ x\in S^d\),
where $P_E$ is the orthogonal projection onto $E$. 
Denote $\overline E$ = $\operatorname{span}^{+}(e,-f)$. 
Thus \(\omega_{\bar E}=-\omega_E\) and \(P_{\bar E}=P_E\).
For unit representatives $z$ and $w$ of $E$ and $F$, respectively,
define
\[
d_Q(E,F)=\arccos|z\cdot\bar w|.
\]
If $z$ and $w$ are replaced by
$e^{i\alpha}z$ and $e^{i\beta}w$, respectively, then
\[
(e^{i\alpha}z)\cdot\overline{e^{i\beta}w}
=e^{i(\alpha-\beta)}(z\cdot\bar w).
\]
Taking absolute values shows that $d_Q(E,F)$ is independent of
the choice of unit representatives.
Moreover, the function $d_Q$ is the restriction of the Fubini--Study distance on
$\mathbb{CP}^d$ to $Q$. Hence $d_Q$ is a metric on $Q$.

\begin{lemma}\label{lem:grassmannian-comparisons}
	For all $E,F\in Q$ and $x\in S^d$, we have
	\begin{align}\label{we-wf}
	\frac{2\sqrt{2}}{\pi}\,d_Q(E,F)
	\leq |\omega_E-\omega_F|
	\leq 2\,d_Q(E,F)
	\end{align}
	and
	\[
	|s_E(x)-s_F(x)|
	\leq \|P_E-P_F\|
	\leq 2\sqrt{2}\,d_Q(E,F).
	\]
\end{lemma}

\begin{proof}
	Let \(E=\operatorname{span}^+ (e,f)\) and \(F=\operatorname{span}^+(g,k)\).
	Write $z=\frac{e+if}{\sqrt{2}}$ and $w=\frac{g+ik}{\sqrt{2}}$.
	Then
	\[
	z\wedge\bar z
	=\frac{1}{2}(e+if)\wedge(e-if)
	=-i\,e\wedge f
	=-i\,\omega_E.
	\]
	Similarly, $w\wedge\bar w=-i\,\omega_F$. 
	Therefore,
	\begin{equation}\label{w}
	\begin{aligned}
		|\omega_E-\omega_F|^2
		&=2-2\langle\omega_E,\omega_F\rangle\\
		&=2-2\langle z\wedge\bar z,w\wedge\bar w\rangle\\
		&=2-2\det
		\begin{pmatrix}
			z\cdot\bar w & z\cdot w\\
			\bar z\cdot\bar w & \bar z\cdot w
		\end{pmatrix}\\
		&=2-2|z\cdot\bar w|^2+2|z\cdot w|^2\\
		&=2\sin^2 d_Q(E,F)+2|z\cdot w|^2.
	\end{aligned}
	\end{equation}
	Moreover, since $z\cdot z=0$ and $z\cdot\bar z=1$, we have
	\[
	|z\cdot\bar w|^2+|z\cdot w|^2\leq |w|^2 = 1.
	\]
	It follows that $|z\cdot w|^2\leq \sin^2 d_Q(E,F)$.
	Hence by \eqref{w},
	\begin{align}\label{we-wf1}
	\sqrt{2}\sin d_Q(E,F)
	\leq |\omega_E-\omega_F|
	\leq 2\sin d_Q(E,F).
	\end{align}
	By the Cauchy--Schwarz inequality, $|z\cdot\bar w|\le1$, so
	$d_Q(E,F)\in[0,\frac{\pi}{2}]$.
	Since \(\frac 2 \pi t\leq \sin t \leq t\) for \(t \in [0,\frac{\pi}{2}]\), we obtain \eqref{we-wf} by \eqref{we-wf1}.

	Since $(e,f)$ is an orthonormal basis of $E$, for every
	$x\in\mathbb R^{d+1}$ we have
	\[
	P_E x
	=(x\cdot e)e+(x\cdot f)f
	=e(e^{\mathsf T}x)+f(f^{\mathsf T}x)
	=(ee^{\mathsf T}+ff^{\mathsf T})x=2\operatorname{Re}(zz^*)x.
	\]
	Hence \(P_E
	=2\operatorname{Re}(zz^*)\).
	Similarly, \(P_F=2\operatorname{Re}(ww^*)\).
	Consequently,
	\begin{equation}\label{pe-pf}
	\begin{aligned}
		\|P_E-P_F\|
		&\leq \|P_E-P_F\|_{\mathrm{HS}}\\
		&= 2\|\operatorname{Re}(zz^*-ww^*)\|_{\mathrm{HS}}\\
		&\leq 2\|zz^*-ww^*\|_{\mathrm{HS}}\\
		&=2\Bigl[\operatorname{tr}\bigl((zz^*-ww^*)^2\bigr)\Bigr]^{\frac 12}\\
		&=2\bigl(2-2\operatorname{tr}(zz^*ww^*)\bigr)^{\frac 12}\\
		&=2\bigl(2-2|z\cdot\bar w|^2\bigr)^{\frac 12}\\
		&=2\sqrt{2}\sin d_Q(E,F),
	\end{aligned}
	\end{equation}
	where $\|\cdot\|_{\mathrm{HS}}$ denotes the Hilbert--Schmidt norm.
	For $x\in S^d$, by \eqref{pe-pf}, we have
	\[
	\begin{aligned}
		|s_E(x)-s_F(x)|
		&=\bigl|| (I-P_E)x|-|(I-P_F)x|\bigr|\\
		&\leq |(P_E-P_F)x|\\
		&\leq \|P_E-P_F\|\\
		&\leq 2\sqrt{2}\,\sin d_Q(E,F)\\
		&\leq 2\sqrt{2}\, d_Q(E,F).
	\end{aligned}
	\]
\end{proof}

Let $\nu$ be the rotation-invariant probability measure on $Q$.
For a unit representative $z$ of $E$, define
\begin{align}\label{def_bz}
b_z(x)=c_{n,d-1}\left(\sqrt2\,x\cdot\bar z\right)^n,\ x\in S^d,
\end{align}
where
\[
c_{n,d-1}=
\left(\frac{\Gamma(n+1+\frac {d-1}{2})}
{\Gamma(1+\frac {d-1}{2})\Gamma(n+1)}\right)^{\frac 12}.
\]
\begin{lemma}\label{lemma2.3}
	The following statements hold.
	\begin{enumerate}
		\item[\textup{(i)}]
		For every unit representative $z$ of $E\in Q$, we have
		\begin{enumerate}
			\item[\textup{(a)}]
			$b_z\in\mathcal H_n$.
			\item[\textup{(b)}]$\|b_z\|_2=1$.
			\item[\textup{(c)}]For $n\geq1$ and $x\in S^d$,
			\[|b_z(x)|
			\leq \sqrt{\frac{d+1}{2}}\,n^{\frac{d-1}{4}}
			\exp\left(-\frac{n\,s_E(x)^2}{2}\right).\]
			\item[\textup{(d)}]
			Define the operator
			$b_z\otimes b_z$ on $\mathcal H_n$ by
			\((b_z\otimes b_z)u=\langle u,b_z\rangle b_z\).
			Then
			\[
			D_n\int_Q b_z\otimes b_z\,d\nu([z])=I,
			\]
			where $\nu$ is the rotation-invariant probability measure on $Q$.
		\end{enumerate}
		
		\item[\textup{(ii)}]
		For unit representatives $z$ and $w$ of $E,F\in Q$,
		respectively, we have 
		\begin{enumerate}
			\item[\textup{(a)}]
			$\langle b_z,b_w\rangle=(w\cdot\bar z)^n$.
			\item[\textup{(b)}]
			\(|\langle b_z,b_w\rangle|
			=\bigl(\cos d_Q(E,F)\bigr)^n
			\leq \exp\left(-\frac{n\,d_Q(E,F)^2}{2}\right)\).
		\end{enumerate}
		
		\item[\textup{(iii)}]
		For every orthonormal basis $\{Y_j\}_{j=1}^{D_n}$ of $\mathcal H_n$ and \(x\in S^d\),
		we have \(\sum\limits_{j=1}^{D_n} |Y_j(x)|^2=D_n\).
	\end{enumerate}
\end{lemma}

\begin{proof}
	\textup{(i)} Let \(E=\operatorname{span}^{+}(e,f)\) and write $z=\frac{e+if}{\sqrt{2}}$.\\
	\textup{(a)}  
	For $n\geq 2$, we have
	\[
	\Delta(\sqrt{2}\,x\cdot\bar z)^n
	=2n(n-1)(\bar z\cdot\bar z)(\sqrt{2}\,x\cdot\bar z)^{n-2}
	=0.
	\]
	For $n=1$, $\Delta(\sqrt{2}\,x\cdot\bar z)=0$.
	Thus $b_z\in\mathcal H_n$.\\
	\textup{(b)} 
	Write $x=ae+bf+y$, where $a,b\in\mathbb R$ and $y\in E^\perp$.
	Since the decomposition is orthogonal, we have
	\(|x|^2=a^2+b^2+|y|^2\), \(\sqrt{2}\,x\cdot\bar z=a-ib\)
	and $dx=da\,db\,dy$.
	Therefore, by Tonelli's theorem,
	\begin{equation}\label{int1}
	\begin{aligned}
		\int_{\mathbb R^{d+1}}
		e^{-|x|^2}|\sqrt{2}\,x\cdot\bar z|^{2n}\,dx
		&=\left(\int_{E^\perp}e^{-|y|^2}\,dy\right)
		\left(2\pi\int_0^\infty e^{-r^2}r^{2n+1}\,dr\right)\\
		&\overset{(*)}{=}\left(\prod_{j=1}^{d-1}\int_{\mathbb R}e^{-y_j^2}\,dy_j\right)
		\left(\pi\int_0^\infty e^{-t}t^n\,dt\right)\\
		&=(\sqrt{\pi})^{d-1}\,\pi\,\Gamma(n+1)\\
		&=\pi^{1+\frac{d-1}{2}}\Gamma(n+1),
	\end{aligned}
	\end{equation}
	where \((*)\) follows from $\dim E^\perp=d-1$ and the substitution $t=r^2$.
	On the other hand, we have
	\begin{equation}\label{int2}
	\begin{aligned}
		&\int_{\mathbb R^{d+1}}
		e^{-|x|^2}|\sqrt{2}\,x\cdot\bar z|^{2n}\,dx\\
		=&|S^d|
		\left(\int_0^\infty e^{-r^2}r^{2n+d}\,dr\right)
		\left(\int_{S^d}|\sqrt{2}\,\omega\cdot\bar z|^{2n}
		\,d\sigma(\omega)\right)\\
		=&\frac{\pi^{\frac{d+1}{2}}}{\Gamma\left(\frac{d+1}{2}\right)}
		\left(\int_0^\infty e^{-t}t^{n+\frac{d-1}{2}}\,dt\right)
		\left(\int_{S^d}|\sqrt{2}\,\omega\cdot\bar z|^{2n}
		\,d\sigma(\omega)\right)\\
		=&\frac{\pi^{\frac{d+1}{2}}
			\Gamma\left(n+\frac{d+1}{2}\right)}
		{\Gamma\left(\frac{d+1}{2}\right)}
		\int_{S^d}|\sqrt{2}\,\omega\cdot\bar z|^{2n}\,d\sigma(\omega).
	\end{aligned}
	\end{equation}
	By \eqref{int1} and \eqref{int2}, we obtain
	\[
	\int_{S^d}|\sqrt{2}\,x\cdot\bar z|^{2n}\,d\sigma(x)
	=\frac{\Gamma(1+\frac{d-1}{2})\Gamma(n+1)}
	{\Gamma(n+1+\frac{d-1}{2})}
	=c_{n,d-1}^{-2}.
	\]
	Hence $\|b_z\|_2=1$.\\
	\textup{(c)} For $x\in S^d$, we have
	\begin{align}\label{xz}
	|\sqrt{2}\,x\cdot\bar z|^2
	=(x\cdot e)^2+(x\cdot f)^2
	=|P_Ex|^2
	=1-s_E(x)^2.
	\end{align}
	By the relation $\Gamma(t+1)=t\Gamma(t)$ and the
	inequality $1+t\leq \exp(t)$ for $t\in\mathbb R$, we have, for $n\geq 1$,
	\[
	\begin{aligned}
		\frac{c_{n+1,d-1}^{\,2}(n+1)^{-\frac{d-1}{2}}}
		{c_{n,d-1}^{\,2}n^{-\frac{d-1}{2}}}
		&=
		\frac{\Gamma\left(n+1+\frac{d+1}{2}\right)}
		{\Gamma\left(n+\frac{d+1}{2}\right)}
		\frac{\Gamma(n+1)}{\Gamma(n+2)}
		\left(\frac{n}{n+1}\right)^{\frac{d-1}{2}}\\
		&=
		\left(1+\frac{d-1}{2(n+1)}\right)
		\left(1-\frac{1}{n+1}\right)^{\frac{d-1}{2}}\\
		&\leq
		\exp\left(\frac{d-1}{2(n+1)}\right)
		\exp\left(-\frac{d-1}{2(n+1)}\right)
		=1.
	\end{aligned}
	\]
	Thus
	\[
	c_{n,d-1}^{\,2}n^{-\frac{d-1}{2}}
	\leq c_{1,d-1}^{\,2}
	=\frac{\Gamma\left(\frac{d+3}{2}\right)}
	{\Gamma\left(\frac{d+1}{2}\right)\Gamma(2)}
	=\frac{d+1}{2},
	\]
	which implies that
	\begin{align}\label{bz}
	c_{n,d-1}\leq \sqrt{\frac{d+1}{2}}\,n^{\frac{d-1}{4}}.
	\end{align}
	Therefore, by \eqref{xz} and \eqref{bz},
	\[
	|b_z(x)|
	\leq\sqrt{\frac{d+1}{2}}\,n^{\frac{d-1}{4}}\bigl(1-s_E(x)^2\bigr)^{\frac{n}{2}}
	\leq \sqrt{\frac{d+1}{2}}\,n^{\frac{d-1}{4}}
	\exp\left(-\frac{n\,s_E(x)^2}{2}\right).
	\]
	\textup{(d)} 
	For every $u\in\mathcal H_n$ and \(\theta \in \mathbb R\), we have
	\[
	\begin{aligned}
		(b_{e^{i\theta}z}\otimes b_{e^{i\theta}z})u
		=\langle u,b_{e^{i\theta}z}\rangle b_{e^{i\theta}z}
		=e^{in\theta}\langle u,b_z\rangle e^{-in\theta}b_z
		=(b_z\otimes b_z)u.
	\end{aligned}
	\]
	Thus $b_z\otimes b_z$ depends only on $[z]$. 
	Set
	\[
	T=\int_Q b_z\otimes b_z\,d\nu([z]).
	\]
	Let $R$ be a rotation of $\mathbb R^{d+1}$ and define $U_R: \mathcal H_n \to \mathcal H_n$ by
	$(U_Ru)(x)=u(R^{-1}x)$. 
	The rotation invariance of \(d\sigma\) implies
	that $U_R$ is unitary. 
	Moreover,
	\begin{align*}
	(U_Rb_z)(x)
	&=c_{n,d-1}\bigl(\sqrt{2}\,R^{-1}x\cdot\bar z\bigr)^n\\
	&=c_{n,d-1}\bigl(\sqrt{2}\,x^{\mathsf T}(R^{-1})^{\mathsf T}\bar z\bigr)^n \\
	&=c_{n,d-1}\bigl(\sqrt{2}\,x\cdot\overline{Rz}\bigr)^n\\
	&=b_{Rz}(x).
	\end{align*}
	Thus, for every $u\in\mathcal H_n$,
	\begin{align}\label{URbz}
	U_R(b_z\otimes b_z)U_R^{-1}u
	=\langle U_R^{-1}u,b_z\rangle U_Rb_z
	=\langle u,b_{Rz}\rangle b_{Rz}
	=(b_{Rz}\otimes b_{Rz})u.
	\end{align}
	By \eqref{URbz} and the rotation invariance of $\nu$, it follows that
	\begin{align}\label{commute}
	U_RTU_R^{-1}
	=\int_Q b_{Rz}\otimes b_{Rz}\,d\nu([z])
	=\int_Q b_z\otimes b_z\,d\nu([z])
	=T.
	\end{align}
	By \cite[Theorem~1.7.2]{DaiXu2013}, the representation
	$R\mapsto U_R$ of the rotation group on the complex space
	$\mathcal H_n$ is irreducible.
	By \eqref{commute} and Schur's lemma \cite[Theorem~16.40(2)]{Hall2013}, there exists $\lambda\in\mathbb C$ such that 
	\begin{align}\label{T=lambdaI}
		T=\lambda I.
	\end{align}
	Let $\{e_j\}_{j=1}^{D_n}$ be an orthonormal basis of $\mathcal H_n$.
	By \textup{(i)(b)},
	\begin{align}\label{trbz}
	\operatorname{tr}(b_z\otimes b_z)
	=\sum_{j=1}^{D_n} \bigl\langle (b_z\otimes b_z)e_j,e_j\bigr\rangle
	=\sum_{j=1}^{D_n} |\langle e_j,b_z\rangle|^2
	=\|b_z\|_2^2
	=1.
	\end{align}
	Therefore, by \eqref{T=lambdaI} and \eqref{trbz},
	\[
	\lambda {D_n}
	=\operatorname{tr}T
	=\int_Q\operatorname{tr}(b_z\otimes b_z)\,d\nu([z])
	=\nu(Q)
	=1.
	\]
	Consequently, \({D_n}T=I\).\\
	\textup{(ii)} 
	 Let \(E=\operatorname{span}^{+}(e,f)\) and \(F=\operatorname{span}^{+}(g,k)\).
	 Write $z=\frac{e+if}{\sqrt{2}}$ and $w=\frac{g+ik}{\sqrt{2}}$.\\
	\textup{(a)} 
	For $x\in\mathbb R^{d+1}$, write $x=ae+bf+y$, where
	$a,b\in\mathbb R$ and $y\in E^\perp$.
	Define
	\[
	R_\theta x
	=(a\cos\theta-b\sin\theta)e
	+(a\sin\theta+b\cos\theta)f+y.
	\]
	For $x\in S^d$, set $x_\theta=R_\theta x$. 
	Then
	\begin{equation}\label{eq:beam-rotation}
			\begin{aligned}
			&(x_\theta\cdot\bar z)^n(x_\theta\cdot w)^n\\
			=&\bigl(e^{-i\theta}(x\cdot\bar z)\bigr)^n
			(x_\theta\cdot w)^n\\
			=&(x\cdot\bar z)^n
			\bigl(e^{-i\theta}(x_\theta\cdot w)\bigr)^n\\
			=&(x\cdot\bar z)^n
			\biggl[e^{-i\theta}\Bigl(
			e^{i\theta}(x\cdot z)(w\cdot\bar z)
			+e^{-i\theta}(x\cdot\bar z)(w\cdot z)
			+((I-P_E)x)\cdot w
			\Bigr)\biggr]^n\\
			=&(x\cdot\bar z)^n
			\Bigl(
			(x\cdot z)(w\cdot\bar z)
			+e^{-2i\theta}(x\cdot\bar z)(w\cdot z)
			+e^{-i\theta}((I-P_E)x)\cdot w
			\Bigr)^n.
		\end{aligned}
	\end{equation}
	By \eqref{eq:beam-rotation}, we obtain
	\begin{equation}\label{eq:beam-angular-average}
		\begin{aligned}
			\frac{1}{2\pi}\int_0^{2\pi}
			(x_\theta\cdot\bar z)^n(x_\theta\cdot w)^n\,d\theta
			=(x\cdot\bar z)^n\bigl((x\cdot z)(w\cdot\bar z)\bigr)^n
			=|x\cdot\bar z|^{2n}(w\cdot\bar z)^n.
		\end{aligned}
	\end{equation}
	Thus, by \eqref{eq:beam-angular-average} and
	$\|b_z\|_2=1$,
	\[
	\begin{aligned}
		\langle b_z,b_w\rangle
		&=2^n c_{n,d-1}^{\,2}
		\int_{S^d} (x\cdot\bar z)^n(x\cdot w)^n\, d\sigma(x)\\
		&=2^n c_{n,d-1}^{\,2}
		\int_{S^d}\frac{1}{2\pi}\int_0^{2\pi}
		(x_\theta\cdot\bar z)^n(x_\theta\cdot w)^n\,
		d\theta\,d\sigma(x)\\
		&=2^n c_{n,d-1}^{\,2}(w\cdot\bar z)^n
		\int_{S^d}|x\cdot\bar z|^{2n}\,d\sigma(x)\\
		&=(w\cdot\bar z)^n\|b_z\|_2^2
		=(w\cdot\bar z)^n.
	\end{aligned}
	\]
	\textup{(b)} By \textup{(ii)(a)}, we have
	\begin{equation}\label{|bz|}
	|\langle b_z,b_w\rangle|
	=|w\cdot\bar z|^n
	=\bigl(\cos d_Q(E,F)\bigr)^n.
	\end{equation}
	For $0\leq t\leq \frac \pi 2$, it is easy to verify that
	$\cos t\leq \exp(-\frac{t^2}{2})$.
	Therefore, by \eqref{|bz|},
	\[
	|\langle b_z,b_w\rangle|
	\leq \exp\left(-\frac{n\,d_Q(E,F)^2}{2}\right).
	\]
	\textup{(iii)}  For any $x,y\in S^d$, choose a rotation $R$ such that $Ry=x$.
	Then $\{U_RY_j\}_{j=1}^{D_n}$ is also an orthonormal basis. 
	Thus there is a unitary matrix $(a_{jk})_{j,k=1}^{D_n}$ such that
	$U_RY_j=\sum\limits_{k=1}^{D_n} a_{jk}Y_k$. 
	Therefore, 
	\begin{equation*}
	\begin{aligned}
		\sum_{j=1}^{D_n} |Y_j(y)|^2
		&=\sum_{j=1}^{D_n} |Y_j(R^{-1}x)|^2\\
		&=\sum_{j=1}^{D_n} |(U_RY_j)(x)|^2\\
		&=\sum_{j=1}^{D_n} \left|\sum_{k=1}^D a_{jk}Y_k(x)\right|^2\\
		&=\sum_{k,\ell=1}^{D_n}
		Y_k(x)\overline{Y_\ell(x)}
		\left(\sum_{j=1}^{D_n} a_{jk}\overline{a_{j\ell}}\right)\\
		&=\sum_{k=1}^{D_n} |Y_k(x)|^2.
	\end{aligned}
	\end{equation*}
	Consequently, $\sum\limits_{j=1}^{D_n} |Y_j(x)|^2$ is
	constant on $S^d$.
	Moreover,
	\[
	\sum\limits_{j=1}^{D_n} |Y_j(x)|^2=\int_{S^d}\sum_{j=1}^{D_n} |Y_j(x)|^2\,d\sigma(x)
	=\sum_{j=1}^{D_n}\|Y_j\|_2^2
	={D_n}.
	\]
\end{proof}

Let $d\geq 2$ and fix an orthogonal decomposition \(\mathbb R^{d+1}=\mathbb R^2\oplus\mathbb R^{d-1}\).
Let $(\varepsilon_1,\varepsilon_2)$ be the standard ordered basis of
$\mathbb R^2$, and let $M_{d-1,2}(\mathbb R)$ denote the space of
$(d-1)\times 2$ real matrices. For
$A=[A_1,A_2]\in M_{d-1,2}(\mathbb R)$, set
\[
E_A=\{(u,Au):u\in\mathbb R^2\}.
\]
We orient $E_A$ by declaring the ordered basis
\[
(v_1,v_2)=\bigl((\varepsilon_1,A_1),(\varepsilon_2,A_2)\bigr)
\]
to be positive. 
Thus an ordered basis $(w_1,w_2)$ of $E_A$
is positively oriented if and only if $\det C>0$, where
$[w_1,w_2]=[v_1,v_2]C$.
Define
\[
[e(A),f(A)]
=\binom{I_2}{A}(I_2+A^{\mathsf T}A)^{-\frac 1 2}.
\]
Then $(e(A),f(A))$ is a positively oriented orthonormal basis of $E_A$.
Set
\[
z(A)=\frac{e(A)+if(A)}{\sqrt{2}}.
\]
By Lemma~\ref{lem:grassmannian-quadric}, $[z(A)]$ represents $E_A$. 
We simply write $s_A=s_{E_A}$ and $b_A=b_{z(A)}$.

For $A\in M_{d-1,2}(\mathbb R)$, set
$|A|=(\operatorname{tr}(A^{\mathsf T}A))^{\frac 12}$.
Let $0<\ell\leq R\ll 1$ and fix a matrix
$A_0\in M_{d-1,2}(\mathbb R)$ with $|A_0|\leq R$.
Set
\begin{align*}
\Omega
=\left\{
A\in M_{d-1,2}(\mathbb R):
|A_{ij}-(A_0)_{ij}|\leq \ell
\text{ for all }1\leq i\leq d-1,\,1\leq j\leq 2
\right\}.
\end{align*}
For $A\in\Omega$, set
\[
z'(A)=e^{-i\arg(z(A)\cdot\overline{z(A_0)})}z(A).
\]
We simply write \(b'_A=b_{z'(A)}\).
For $A,B\in\Omega$, set $X=A-A_0$ and $Y=B-A_0$.
Writing $X=[X_1,X_2]$ and $Y=[Y_1,Y_2]$, define
\begin{align*}
\omega(X,Y)=X_1\cdot Y_2-X_2\cdot Y_1.
\end{align*}
We also set
\[
F_{A_0}(X,Y)
=\log\!\left(z'(A_0+Y)\cdot\overline{z'(A_0+X)}\right),
\]
where $\log$ denotes the branch of the complex logarithm near $1$ with $\log 1=0$.
\begin{lemma}\label{lem:local-log-expansion}
	With the notation above, assume that $R$ is sufficiently small
	depending only on $d$. Then
	\begin{equation*}
		F_{A_0}(X,Y)
		=-\frac14|X-Y|^2+\frac{i}{2}\omega(X,Y)
		+O_d(R\ell^2+\ell^3),
	\end{equation*}
	where, for $g\geq 0$, the notation $f=O_d(g)$ means that
	$|f|\leq C_d g$ for some constant $C_d>0$ depending only on $d$.
\end{lemma}

\begin{proof}
	By the definition of $z(A)$,
	\[
	z(A)=\frac{1}{\sqrt{2}}\binom{I_2}{A}
	(I_2+A^{\mathsf T}A)^{-\frac 12}\binom{1}{i}.
	\]
	Since $(I_2+A^{\mathsf T}A)^{-\frac 12}$ is real symmetric and positive definite, we have
	\[
	z(A)\cdot\overline{z(0)}
	=\frac12(1,-i)(I_2+A^{\mathsf T}A)^{-\frac 1 2}\binom{1}{i}=\frac12\operatorname{tr}(I_2+A^{\mathsf T}A)^{-\frac1 2}>0.
	\]
	When $A_0=0$, we have $z'(A)=z(A)$. 
	The same argument gives	$z'(B)=z(B)$.
	Write
	\[
	A^{\mathsf T}A
	=U\operatorname{diag}(\lambda_1,\lambda_2)U^{\mathsf T},
	\]
	where $U$ is orthogonal and $\lambda_1,\lambda_2\geq 0$.
	Since \(0\leq \lambda_j\leq \lambda_1+\lambda_2
	=\operatorname{tr}(A^{\mathsf T}A)=|A|^2\), 
	we have
	\begin{equation}\label{I+AA}
	\begin{aligned}
		(I_2+A^{\mathsf T}A)^{-\frac 12}
		&=U\operatorname{diag}\bigl((1+\lambda_1)^{-\frac 12},
		(1+\lambda_2)^{-\frac 12}\bigr)U^{\mathsf T}\\
		&=U\operatorname{diag}\!\left(
		1-\tfrac12\lambda_1+O(\lambda_1^2),
		1-\tfrac12\lambda_2+O(\lambda_2^2)
		\right)U^{\mathsf T}\\
		&=I_2-\tfrac12A^{\mathsf T}A+O_d(|A|^4).
	\end{aligned}
	\end{equation}
	Therefore, by \eqref{I+AA},
	\begin{equation}\label{zYzX}
	\begin{aligned}
		&z(Y)\cdot\overline{z(X)}\\
		=&\frac12(1,-i)(I_2+X^{\mathsf T}X)^{-\frac 12} (I_2+X^{\mathsf T}Y)
		(I_2+Y^{\mathsf T}Y)^{-\frac 12}\binom{1}{i}\\
		=&\frac12(1,-i)
		\left(I_2+X^{\mathsf T}Y-\frac12X^{\mathsf T}X
		-\frac12Y^{\mathsf T}Y\right)\binom{1}{i}
		+O_d(\ell^4)\\
		=&1+\frac12\operatorname{tr}(X^{\mathsf T}Y)
		-\frac14\operatorname{tr}(X^{\mathsf T}X)-\frac14\operatorname{tr}(Y^{\mathsf T}Y)
		+\frac{i}{2}(X_1\cdot Y_2-X_2\cdot Y_1)
		+O_d(\ell^4)\\
		=&1-\frac14|X-Y|^2+\frac{i}{2}\omega(X,Y)+O_d(\ell^4).
	\end{aligned}
	\end{equation}
	Consequently, by \eqref{zYzX} and the fact that
	$\log(1+w)=w+O(|w|^2)$ as $w\to 0$, we obtain
	\begin{align}\label{F_0}
	F_0(X,Y)
	=-\frac14|X-Y|^2+\frac{i}{2}\omega(X,Y)+O_d(\ell^4).
	\end{align}

	For general $A_0$ and $t$ near $0$, by the definition of $z'$,
	\begin{align*}
		z'(A_0+tX)\cdot\overline{z'(A_0)}
		&=e^{-i\arg\left(z(A_0+tX)\cdot\overline{z(A_0)}\right)}
		\bigl(z(A_0+tX)\cdot\overline{z(A_0)}\bigr)\\
		&=\left|z(A_0+tX)\cdot\overline{z(A_0)}\right|.
	\end{align*}
	Thus \(0 < z'(A_0+tX)\cdot\overline{z'(A_0)} \leq 1\),
	which implies that
	\[
	F_{A_0}(tX,0)
	=\log\!\left|z(A_0+tX)\cdot\overline{z(A_0)}\right|\leq 0.
	\]
	Since $F_{A_0}(0,0)=0$, \(F_{A_0}(tX,0)\) attains a local
	maximum at $t=0$. Hence
	\[
	\left.\frac{d}{dt}F_{A_0}(tX,0)\right|_{t=0}
	=0.
	\]
	Similarly, \(	\left.\frac{d}{dt}F_{A_0}(0,tY)\right|_{t=0}=0\).
    Therefore,
	\begin{align*}
		\left.\frac{d}{dt}F_{A_0}(tX,tY)\right|_{t=0}
		=\left.\frac{d}{dt}F_{A_0}(tX,0)\right|_{t=0}+\left.\frac{d}{dt}F_{A_0}(0,tY)\right|_{t=0}=0.
	\end{align*}
	 Since \(|A_0|+|X|+|Y|\le (1+2\sqrt{2(d-1)})R\) and $(A_0,X,Y)\mapsto F_{A_0}(X,Y)$ is smooth near the origin, for $R$ sufficiently small, we have
	\begin{align}\label{san}
	\left|\frac{d^3}{dt^3}F_{A_0}(tX,tY)\right|
	\le C_d(|X|+|Y|)^3,\ 0\le t\le1,
	\end{align}
	and
	\begin{align}\label{er}
	\left|
	\frac{d}{ds}\left(
	\left.\frac{d^2}{dt^2}F_{sA_0}(tX,tY)\right|_{t=0}
	\right)
	\right|
	\le C_dR(|X|+|Y|)^2,\ 0\le s\le1.
	\end{align}
	By \eqref{er} and \eqref{F_0}, 
	\begin{equation}\label{er0}
	\begin{aligned}
	&\left.\frac{d^2}{dt^2}F_{A_0}(tX,tY)\right|_{t=0}\\
	=&
	\left.\frac{d^2}{dt^2}F_0(tX,tY)\right|_{t=0}
	+O_d\bigl(R(|X|+|Y|)^2\bigr)\\
	=&-\frac12|X-Y|^2+i\omega(X,Y)+O_d\bigl(R(|X|+|Y|)^2\bigr).
	\end{aligned}
	\end{equation}
	Therefore, by Taylor's formula at $t=0$, \eqref{san} and \eqref{er0},
	\begin{align*}
	&F_{A_0}(X,Y)\\
	=&\frac12\left.\frac{d^2}{dt^2}
	F_{A_0}(tX,tY)\right|_{t=0}
	+O_d\bigl((|X|+|Y|)^3\bigr)\\
	=&-\frac14|X-Y|^2+\frac{i}{2}\omega(X,Y)
	+O_d\bigl(R(|X|+|Y|)^2\bigr)+O_d\bigl((|X|+|Y|)^3\\
	=&-\frac14|X-Y|^2+\frac{i}{2}\omega(X,Y)
	+O_d(R\ell^2+\ell^3),
	\end{align*}
	where the last equation follows from \(|X|+|Y|\le 2\sqrt{2(d-1)}\ell\).
\end{proof}
\begin{lemma}\label{lem:graph-coordinate-distance}
	For all $A,B\in M_{d-1,2}(\mathbb R)$ with
	$|A|,|B|\le 1$,
	\begin{align}\label{dQ}
	\frac{1}{12}|A-B|\le d_Q(E_A,E_B)\le \frac{3\sqrt2 \pi}{2}|A-B|.
	\end{align}
	Moreover, for every $n\ge 1$, 
	\[
	|\langle b_A',b_B'\rangle|=|\langle b_A,b_B\rangle|
	\le \exp\!\left(-\frac{n}{288}|A-B|^2\right).
	\]
\end{lemma}
\begin{proof}
	For $A=[A_1,A_2]$, set \(W(A)=(\varepsilon_1,A_1)\wedge(\varepsilon_2,A_2)\).
	Since $[e(A),f(A)]=[v_1,v_2](I_2+A^{\mathsf T}A)^{-\frac12}$,
	we have
	\[
	\omega_{E_A}
	=e(A)\wedge f(A)
	=\det\!\bigl((I_2+A^{\mathsf T}A)^{-\frac12}\bigr)W(A).
	\]
	Moreover, notice that $|\omega_{E_A}|=1$. 
	Hence
	\begin{align}\label{we}
		\omega_{E_A}=\frac{W(A)}{|W(A)|}.
	\end{align}
	It is easy to see that
	\begin{align}\label{waw0=1}
	\langle W(A),W(0)\rangle
	=
	\det\begin{pmatrix}
		(\varepsilon_1,A_1)\cdot(\varepsilon_1,0)
		&(\varepsilon_1,A_1)\cdot(\varepsilon_2,0)\\
		(\varepsilon_2,A_2)\cdot(\varepsilon_1,0)
		&(\varepsilon_2,A_2)\cdot(\varepsilon_2,0)
	\end{pmatrix}
	=\det I_2=1.
	\end{align}
	Thus
	\begin{align}\label{wage1}
	1=|\langle W(A),W(0)\rangle|
	\le |W(A)|\,|W(0)|=|W(A)|.
	\end{align}
	Therefore, by \eqref{we} and \eqref{wage1},
	\begin{equation}\label{leq6}
	\begin{aligned}
		&|\omega_{E_A}-\omega_{E_B}|\\
		\le& \frac{|W(A)-W(B)|}{|W(A)|}
		+ |W(B)|\left|\frac1{|W(A)|}-\frac1{|W(B)|}\right|\\
		\le& \frac{2|W(A)-W(B)|}{|W(A)|}\\
		\le& 2|W(A)-W(B)|\\
		\leq& 2\bigl(|A-B| + |(0,A_1-B_1)\wedge(0,A_2)
		+(0,B_1)\wedge(0,A_2-B_2)|\bigr)\\
		\le& 2\bigl(|A-B|+|A_1-B_1|\,|A_2|
		+|B_1|\,|A_2-B_2|\bigr)\\
		\le& 2(1+|A|+|B|)|A-B|\\
		\le& 6|A-B|.
	\end{aligned}
	\end{equation}
	On the other hand, notice that
	\[
	|W(A)|
	\le |(\varepsilon_1,A_1)|\,|(\varepsilon_2,A_2)|
	=\sqrt{1+|A_1|^2}\sqrt{1+|A_2|^2}
	\le 2.
	\]
	Then by \eqref{we} and \eqref{waw0=1}, we have
	\[
	\langle \omega_{E_A},W(0)\rangle=\frac{ \langle W(A), W(0)\rangle }{|W(A)|}=\frac{1}{|W(A)|}\ge\frac12,
	\]
	which implies that
	\[
	W(A)=\omega_{E_A} |W(A)| =\frac{\omega_{E_A}}{\langle\omega_{E_A},W(0)\rangle}.
	\]
	The same identities hold for $B$. 
	Therefore,
	\begin{equation}\label{geq6}
		\begin{aligned}
			|W(A)-W(B)|
			&\le
			\frac{|\omega_{E_A}-\omega_{E_B}|}
			{\langle\omega_{E_A},W(0)\rangle}
			+
			\frac{|\langle\omega_{E_A}-\omega_{E_B},W(0)\rangle|}
			{\langle\omega_{E_A},W(0)\rangle
				\langle\omega_{E_B},W(0)\rangle}\\
			&\le 2|\omega_{E_A}-\omega_{E_B}| +4 |\omega_{E_A}-\omega_{E_B}| |W(0)|\\
			&\le 6|\omega_{E_A}-\omega_{E_B}|,
		\end{aligned}
	\end{equation}
	 Consequently, by \eqref{leq6} and \eqref{geq6},
	\begin{align}\label{wa-wbfanwei}
	\frac16|A-B|\le |\omega_{E_A}-\omega_{E_B}|
	\le 6|A-B|.
	\end{align}
	By \eqref{wa-wbfanwei} and \eqref{we-wf}, \eqref{dQ} holds.
	
	By Lemma~\ref{lemma2.3} (ii)(a), (ii)(b) and \eqref{dQ}, we have
	\[
	|\langle b_A',b_B'\rangle|
	=|\langle b_A,b_B\rangle|
	\le \exp\!\left(-\frac n2 d_Q(E_A,E_B)^2\right)
	\le \exp\!\left(-\frac{n}{288}|A-B|^2\right).
	\]
\end{proof}

For $u=(q,p)\in\mathbb R^{d-1}\times\mathbb R^{d-1}$, set
\[
\psi_u(t)
=2^{\frac{d-1}{4}}
\exp\!\left(2\pi i\,p\cdot\left(t-\frac{q}{2}\right)
-\pi|t-q|^2\right),
\ t\in\mathbb R^{d-1}.
\]
For $F\in L^1(\mathbb R^{d-1})$, we define its Fourier transform by
\[
\widehat F(p)
=\int_{\mathbb R^{d-1}}F(t)e^{-2\pi i\,p\cdot t}\,dt,
\ p\in\mathbb R^{d-1}.
\]

\begin{lemma}\label{lem:gaussian-packet-inner-product}
	For all $u=(q,p),v=(q',p')\in\mathbb R^{d-1}\times\mathbb R^{d-1}$, 
	\begin{equation*}
		\langle\psi_u,\psi_v\rangle
		=
		\exp\!\left(
		-\frac{\pi}{2}|u-v|^2
		+\pi i(p\cdot q'-q\cdot p')
		\right).
	\end{equation*}
\end{lemma}
\begin{proof}
	Let $u=(q,p)$ and $v=(q',p')$. By the definition of $\psi_u$,
	\begin{align*}
		\langle\psi_u,\psi_v\rangle
		&=2^{\frac{d-1}{2}}\int_{\mathbb R^{d-1}}
		e^{-\pi(|t-q|^2+|t-q'|^2)}
		e^{2\pi i\left(p\cdot(t-\frac{q}{2})
			-p'\cdot(t-\frac{q'}{2})\right)}\,dt\\
		&\overset{(*)}{=}\exp\!\left(-\frac{\pi}{2}|q-q'|^2
		+\pi i(p\cdot q'-q\cdot p')\right)
		\int_{\mathbb R^{d-1}}
		e^{-\pi|y|^2}e^{2\pi i\frac{p-p'}{\sqrt{2}}\cdot y}\,dy\\
		&=\exp\!\left(-\frac{\pi}{2}|q-q'|^2
		+\pi i(p\cdot q'-q\cdot p')\right)\widehat{e^{-\pi |\,\cdot\,|^2}}\left(-\frac{p-p'}{\sqrt 2}\right)\\
		&=\exp\!\left(-\frac{\pi}{2}|q-q'|^2
		+\pi i(p\cdot q'-q\cdot p')\right)
		\exp\!\left(-\frac{\pi}{2}|p-p'|^2\right)\\
		&=\exp\!\left(-\frac{\pi}{2}|u-v|^2
		+\pi i(p\cdot q'-q\cdot p')\right),
	\end{align*}
	where $(*)$ follows from the change of variables
	$y=\sqrt{2}\left(t-\frac{q+q'}{2}\right)$.
\end{proof}

For $n\ge 1$ and $u=(q,p)\in\mathbb R^{d-1} \times \mathbb R^{d-1}$, set
\[A(u)=A_0+\sqrt{\frac{2\pi}{n}}\,[q,-p],\]
and
\[
\Omega_*=
\left[-\sqrt{\frac{n}{2\pi}}\,\ell,\,
\sqrt{\frac{n}{2\pi}}\,\ell\right]^{2(d-1)}.
\]
Then $A(\Omega_*)=\Omega$.

\begin{lemma}\label{lem:wq-gram-model}
	If $e_n=n(R\ell^2+\ell^3)\le 1$,
	then there exists a constant $C_d>0$, depending only on $d$,
	such that
	\begin{equation*}
		\bigl|\langle b'_{A(u)},b'_{A(v)}\rangle
		-\langle\psi_u,\psi_v\rangle\bigr|
		\le C_d e_n\exp\!\left(-\frac{\pi}{2}|u-v|^2\right)
	\end{equation*}
	for all $u,v\in\Omega_*$.
\end{lemma}

\begin{proof}
	Let $u=(q,p)$ and $v=(q',p')$ belong to $\Omega_*$. 
	By Lemma~\ref{lem:local-log-expansion} and the definitions of $A(u)$
	and $\omega$, we have
	\begin{equation}\label{nFa}
		\begin{aligned}
			&nF_{A_0}\bigl(A(u)-A_0,A(v)-A_0\bigr)\\
			=&-\frac{n}{4}|A(u)-A(v)|^2
			+\frac{in}{2}\omega\bigl(A(u)-A_0,A(v)-A_0\bigr)
			+O_d(e_n)\\
			=&-\frac{\pi}{2}\bigl(|q-q'|^2+|p-p'|^2\bigr)
			+\pi i\bigl(q\cdot(-p')-(-p)\cdot q'\bigr)
			+O_d(e_n)\\
			=&-\frac{\pi}{2}|u-v|^2
			+\pi i(p\cdot q'-q\cdot p')
			+O_d(e_n).
		\end{aligned}
	\end{equation}
	Hence, by Lemma~\ref{lemma2.3} (ii)(a), \eqref{nFa} and Lemma~\ref{lem:gaussian-packet-inner-product},
	\begin{equation}\label{2.8.1}
			\begin{aligned}
			\langle b'_{A(u)},b'_{A(v)}\rangle
			&=\left(z'(A(v))\cdot\overline{z'(A(u))}\right)^n\\
			&=\exp\!\left(nF_{A_0}\bigl(A(u)-A_0,A(v)-A_0\bigr)\right)\\
			&=\exp\!\left(-\frac{\pi}{2}|u-v|^2
			+\pi i(p\cdot q'-q\cdot p')\right)\exp\!\bigl(O_d(e_n)\bigr)\\
			&=\langle\psi_u,\psi_v\rangle\exp\!\bigl(O_d(e_n)\bigr).
		\end{aligned}
	\end{equation}
	Choose $C_d'>0$ depending only on $d$ such that
	$|O_d(e_n)|\le C_d'e_n$. Since $e_n\le 1$, we have
	\begin{equation}\label{2.8.2}
		\begin{aligned}
			\left|\exp\!\bigl(O_d(e_n)\bigr)-1\right|
			&\le \left|O_d(e_n)\right|
			\exp\!\left(\left|O_d(e_n)\right|\right)\\
			&\le C_d'e_n \exp(C_d' e_n)\\
			&\le C_d' \exp(C_d')e_n.
		\end{aligned}
	\end{equation}
	Therefore, by \eqref{2.8.1}, \eqref{2.8.2} and Lemma~\ref{lem:gaussian-packet-inner-product},
	\begin{align*}
		\bigl|\langle b'_{A(u)},b'_{A(v)}\rangle
		-\langle\psi_u,\psi_v\rangle\bigr|
		&=\left|\exp\bigl(O_d(e_n)\bigr)-1\right| |\langle\psi_u,\psi_v\rangle|\\
		&\le C_d' \exp(C_d') e_n|\langle\psi_u,\psi_v\rangle|\\
		&=C_d' \exp(C_d') e_n\exp\!\left(-\frac{\pi}{2}|u-v|^2\right).
	\end{align*}
	The conclusion follows by taking $C_d=C_d' \exp(C_d')$.
\end{proof}

\section{The local Wilson construction}

We use the Wilson construction of Daubechies, Jaffard, and
Journ\'e~\cite{Wilson}. 
Choose a real-valued even function  $g_0\in C_c^\infty(\mathbb R)$ such that
\(\operatorname{supp}g_0
=\left[-\frac{1}{2},\frac{1}{2}\right]\) and \(g_0(t)>0\) for \(|t|<\frac{1}{2}\).
Define
\[
g(t)={\sqrt{2}\,g_0(t)}
{\left(\displaystyle\sum_{k\in\mathbb Z}
	g_0\!\left(t-\frac{k}{2}\right)^2\right)^{-\frac12}}.
\]
It is easy to verify that $g$ is real-valued
and even, with \(g\in C_c^\infty(\mathbb R)\) and \(\operatorname{supp}g=\left[-\frac12,\frac12\right]\).

For $j\in\mathbb Z$ and $\kappa\in\mathbb N=\{1,2,\ldots\}$,
define $w_{j,\kappa}\colon\mathbb R\to\mathbb R$ by
\[
w_{j,\kappa}(t)
=
\begin{cases}
	\sqrt{2}\,g\!\left(t-\frac{j}{2}\right)\cos\!\left(2\pi\kappa\left(t-\frac{j}{2}\right)\right),
	& j+\kappa \text{ is even},\\[2pt]
	\sqrt{2}\,g\!\left(t-\frac{j}{2}\right)\sin\!\left(2\pi\kappa\left(t-\frac{j}{2}\right)\right),
	& j+\kappa \text{ is odd}.
\end{cases}
\]
\begin{lemma}\label{lem:one-dimensional-wilson}
	The family $\{w_{j,\kappa}\}_{j\in\mathbb Z,\,\kappa\geq1}$
	is orthonormal in $L^2(\mathbb R)$.
\end{lemma}

\begin{proof}
	For every $t\in\mathbb R$,
	\begin{equation}\label{=2}
		\sum_{j\in\mathbb Z}g\!\left(t-\frac{j}{2}\right)^2
		=
		\frac{2\displaystyle\sum_{j\in\mathbb Z}
			g_0\!\left(t-\frac{j}{2}\right)^2}
		{\displaystyle\sum_{k\in\mathbb Z}
			g_0\!\left(t-\frac{j}{2}-\frac{k}{2}\right)^2}
		=
		\frac{2\displaystyle\sum_{j\in\mathbb Z}
			g_0\!\left(t-\frac{j}{2}\right)^2}
		{\displaystyle\sum_{\ell\in\mathbb Z}
			g_0\!\left(t-\frac{\ell}{2}\right)^2}
		=2.
	\end{equation}
	Hence 
	\[
	\begin{aligned}
		\|g\|_2^2
		=\sum_{j\in\mathbb Z}
		\int_{-\frac j2}^{\frac{1-j}{2}}g(u)^2\,du
		=\int_0^{\frac 12}\sum_{j\in\mathbb Z}
		g\!\left(t-\frac j2\right)^2\,dt
		=\int_0^{\frac 12}2\,dt=1.
	\end{aligned}
	\]
	For $j,k\in\mathbb Z$, set
	$G_{j,k}(t)=e^{2\pi i kt}g(t-\frac j2)$.
	For \(f\in L^2(\mathbb R)\),
	by Parseval's identity, \(\operatorname{supp}g=\left[-\frac12,\frac12\right]\) and \eqref{=2},
	\begin{equation*}
		\begin{aligned}
			\sum_{j,k\in\mathbb Z}|\langle f,G_{j,k}\rangle|^2
			&=\sum_{j\in\mathbb Z}\sum_{k\in\mathbb Z}
			\left|\int_{\left[\frac{j-1}{2},\frac{j+1}{2}\right]} f(t)g\!\left(t-\frac j2\right)
			e^{-2\pi i kt}\,dt\right|^2\\
			&=\sum_{j\in\mathbb Z}\int_{\left[\frac{j-1}{2},\frac{j+1}{2}\right]}|f(t)|^2g\!\left(t-\frac j2\right)^2\,dt\\
			&=\sum_{j\in\mathbb Z}\int_{\mathbb R}
			|f(t)|^2g\!\left(t-\frac j2\right)^2\,dt\\
			&=\int_{\mathbb R}|f(t)|^2
			\sum_{j\in\mathbb Z}g\!\left(t-\frac j2\right)^2\,dt\\
			&=2\|f\|_2^2.
		\end{aligned}
	\end{equation*}
	Therefore, by the Wilson basis theorem
	\cite[Theorem~1.1]{GrochenigWalnut}, $\{w_{j,\kappa}\}_{j\in\mathbb Z,\,\kappa\geq1}$
	is orthonormal in $L^2(\mathbb R)$.
\end{proof}

For $\lambda=(q,\kappa)\in
(\frac{1}{2}\mathbb Z)^{d-1}\times\mathbb N^{d-1}$, define 
\[
O_\lambda=
\bigl\{(q,\sigma\kappa):\sigma\in\{-1,1\}^{d-1}\bigr\},
\] 
where $\sigma\kappa=(\sigma_1\kappa_1, \ldots, \sigma_{d-1}\kappa_{d-1})$.
Set
\[
W_\lambda(t)=\prod_{i=1}^{d-1}w_{2q_i,\kappa_i}(t_i).
\]
	Let $\lambda=(q,\kappa)$ and $\mu=(q',\kappa')$. By Lemma~\ref{lem:one-dimensional-wilson},
\begin{equation}\label{Wlambda_orth}
\begin{aligned}
	\langle W_\lambda,W_\mu\rangle
	&=\int_{\mathbb R^{d-1}}
	\prod_{i=1}^{d-1}
	\left(
	w_{2q_i,\kappa_i}(t_i)
	\overline{w_{2q_i',\kappa_i'}(t_i)}
	\right)\,dt_1\cdots dt_{d-1}\\
	&=\prod_{i=1}^{d-1}
	\int_{\mathbb R}
	w_{2q_i,\kappa_i}(t_i)
	\overline{w_{2q_i',\kappa_i'}(t_i)}\,dt_i\\
	&=\prod_{i=1}^{d-1}
	\langle w_{2q_i,\kappa_i},w_{2q_i',\kappa_i'}\rangle	=\delta_{\lambda\mu}.
\end{aligned}
\end{equation}
Therefore, $\{W_\lambda\}_{\lambda \in (\frac{1}{2}\mathbb Z)^{d-1}\times\mathbb N^{d-1}}$ is orthonormal in
$L^2(\mathbb R^{d-1})$.
\begin{lemma}\label{lem:wilson-packet-decay}
	 For every $M>0$, there exists
	$C_{d,M}>0$, independent of $\lambda$ and $u$, such that
	\begin{equation*}
		|\langle W_\lambda,\psi_u\rangle|
		\leq C_{d,M}\sum_{\xi\in O_\lambda}(1+|u-\xi|)^{-M},
		\ u\in\mathbb R^{d-1} \times \mathbb R^{d-1} .
	\end{equation*}
\end{lemma}
\begin{proof}
	By Euler's formulas, for each $i=1,\ldots,d-1$,
	\[
	w_{2q_i,\kappa_i}(t_i)
	=
	\begin{cases}
		\frac{g(t_i-q_i)}{\sqrt{2}}
		\left(
		e^{2\pi\mathrm{i}\kappa_i(t_i-q_i)}
		+e^{-2\pi\mathrm{i}\kappa_i(t_i-q_i)}
		\right),
		& 2q_i+\kappa_i \text{ is even},\\[6pt]
		\frac{g(t_i-q_i)}{\mathrm{i}\sqrt{2}}
		\left(
		e^{2\pi\mathrm{i}\kappa_i(t_i-q_i)}
		-e^{-2\pi\mathrm{i}\kappa_i(t_i-q_i)}
		\right),
		& 2q_i+\kappa_i \text{ is odd}.
	\end{cases}
	\]
	Hence,
	\begin{equation}\label{Wlambda}
	\begin{aligned}
		W_\lambda(t)
		&=\prod_{i=1}^{d-1}w_{2q_i,\kappa_i}(t_i)\\
		&=2^{-\frac{d-1}{2}}
		\prod_{i=1}^{d-1}g(t_i-q_i)
		\prod_{\substack{1\leq i\leq d-1\\
				2q_i+\kappa_i\text{ is even}}}
		\left(
		\sum_{\sigma_i\in\{-1,1\}}
		e^{2\pi\mathrm{i}\sigma_i\kappa_i(t_i-q_i)}
		\right)\\
		&\qquad\times
		\prod_{\substack{1\leq i\leq d-1\\
				2q_i+\kappa_i\text{ is odd}}}
		\left(
		\sum_{\sigma_i\in\{-1,1\}}
		\frac{\sigma_i}{\mathrm{i}}\,
		e^{2\pi\mathrm{i}\sigma_i\kappa_i(t_i-q_i)}
		\right)\\
		&=2^{-\frac{d-1}{2}}
		\prod_{i=1}^{d-1}g(t_i-q_i)
		\sum_{\sigma\in\{-1,1\}^{d-1}}
		\left(
		\prod_{\substack{1\leq i\leq d-1\\
				2q_i+\kappa_i\text{ is odd}}}
		\frac{\sigma_i}{\mathrm{i}}
		\right)
		e^{2\pi\mathrm{i}(\sigma\kappa)\cdot(t-q)}\\
		&=\left(\prod_{i=1}^{d-1}g(t_i-q_i)\right)
		\sum_{\sigma\in\{-1,1\}^{d-1}}
		c_{(q,\sigma\kappa)}e^{2\pi\mathrm{i}(\sigma\kappa)\cdot t},
	\end{aligned}
	\end{equation}
	where
	\[
	c_{(q,\sigma\kappa)}
	=
	2^{-\frac{d-1}{2}}
	e^{-2\pi\mathrm{i}(\sigma\kappa)\cdot q}
	\prod_{\substack{1\le i\le d-1\\
			2q_i+\kappa_i\text{ is odd}}}
	\frac{\sigma_i}{\mathrm{i}}.
	\]
	By the definition of $\psi_u$, for $u=(x,y)$ we have
	\begin{align}\label{psiu=psi0}
	\psi_u(t)
	=e^{-\pi\mathrm{i}y\cdot x}
	e^{2\pi\mathrm{i}y\cdot t}\psi_0(t-x).
	\end{align}
	Therefore, by \eqref{Wlambda} and \eqref{psiu=psi0},
	\begin{equation}\label{Wpsi}
	\begin{aligned}
		&|\langle W_\lambda,\psi_u\rangle|\\
		=&
		\left|
		\sum_{\sigma\in\{-1,1\}^{d-1}} c_{(q,\sigma\kappa)}
		\int_{\mathbb R^{d-1}}
		\left(\prod_{i=1}^{d-1}g(t_i-q_i)\right)
		\overline{\psi_0(t-x)}
		e^{-2\pi\mathrm{i}(y-\sigma\kappa)\cdot t}\,dt
		\right|\\
		\leq&
		\sum_{\sigma\in\{-1,1\}^{d-1}}
		\left|
		\int_{\mathbb R^{d-1}}
		\left(\prod_{i=1}^{d-1}g(t_i)\right)
		\overline{\psi_0(t+q-x)}
		e^{-2\pi\mathrm{i}(y-\sigma\kappa)\cdot t}\,dt
		\right|.
	\end{aligned}
	\end{equation}
	For every $M>0$, by \cite[Theorem~2.3]{GrochenigZimmermann}, there exists \(C_{d,M}>0\), independent of $\lambda$ and $u$, such that
	\begin{equation}\label{psiin}
	\begin{aligned}
	&\left|
	\int_{\mathbb R^{d-1}}
	\left(\prod_{i=1}^{d-1}g(t_i)\right)
	\overline{\psi_0(t+q-x)}
	e^{-2\pi\mathrm{i}(y-\sigma\kappa)\cdot t}\,dt
	\right|\\
	\leq &C_{d,M}
	\bigl(1+|x-q|+|y-\sigma\kappa|\bigr)^{-M}.
	\end{aligned}
	\end{equation}
	Consequently, by \eqref{Wpsi} and \eqref{psiin}, we have 
	\begin{align*}
		|\langle W_\lambda,\psi_u\rangle| &\leq C_{d,M} \sum_{\sigma\in\{-1,1\}^{d-1}}\bigl(1+|x-q|+|y-\sigma\kappa|\bigr)^{-M}\\
		&\leq C_{d,M}\sum_{\xi\in O_\lambda}(1+|u-\xi|)^{-M}.
	\end{align*}
\end{proof}

Choose $B$ with $1\ll B<\frac12\sqrt{\frac{n}{2\pi}}\,\ell$.
For $u\in\Omega_*$, write
\[
\operatorname{dist}(u,\partial\Omega_*)
=\inf_{z\in\partial\Omega_*}|u-z|,
\]
where $|\cdot|$ denotes the Euclidean norm on $\mathbb R^{d-1} \times \mathbb R^{d-1}$.
Set
\[
\Lambda_B=
\left\{
\lambda=(q,\kappa)\in
\Omega_*\cap
\left((\tfrac12\mathbb Z)^{d-1}\times\mathbb N^{d-1}\right)
:\operatorname{dist}(\lambda,\partial\Omega_*)\geq B
\right\}.
\]
For each $\lambda=(q,\kappa)\in\Lambda_B$ and \(\sigma\in\{-1,1\}^{d-1}\), define
\begin{align}\label{def_vlambda}
v_\lambda
=\int_{\Omega_*}\langle W_\lambda,\psi_u\rangle\,b'_{A(u)}\,du,
\end{align}
and
\[
E_{\lambda,\sigma}=E_{A\bigl((q,\sigma\kappa)\bigr)}
=\left\{
\left(w,\left(A_0+\sqrt{\frac{2\pi}{n}}\,[q,-\sigma\kappa]\right)w\right)
:w\in\mathbb R^2
\right\}.
\] 
For $\lambda=(q,\kappa),\,\mu=(q',\kappa')\in\Lambda_B$, set
\[
\rho_*(\lambda,\mu)=|q-q'|+|\kappa-\kappa'|.
\]
\begin{lemma}\label{lem:local-packet-estimates}
	Assume that
	$e_n=n(R\ell^2+\ell^3)\leq 1$.
	For every $M,s>0$, there exists a constant $C_{d,M,s}>0$,
	independent of $n,R,\ell,B$ and the choice of $\Omega$, such that for all $\lambda,\mu\in\Lambda_B$,
	\begin{equation}\label{eq:wq-local-gram}
		\bigl|\langle v_\lambda,v_\mu\rangle-\delta_{\lambda\mu}\bigr|
		\leq C_{d,M,s}(e_n+B^{-s})
		\bigl(1+\rho_*(\lambda,\mu)\bigr)^{-M}.
	\end{equation}
	For every $M>0$, there exists a constant $C_{d,M}>0$,
	independent of $n,R,\ell,B$ and the choice of $\Omega$, such that for all $\lambda\in\Lambda_B$ and $x\in S^d$,
	\begin{equation}\label{eq:wq-local-field}
		|v_\lambda(x)|
		\leq C_{d,M}n^{\frac{d-1}{4}}
		\sum_{\sigma\in\{-1,1\}^{d-1}}
		\bigl(1+\sqrt n\,s_{E_{\lambda,\sigma}}(x)\bigr)^{-M}.
	\end{equation}
\end{lemma}
\begin{proof}
	By \eqref{psiu=psi0}, we have
	\[
	\begin{aligned}
		&\langle W_\lambda,\psi_u\rangle\psi_u(t)\\
		=&
		\left(
		e^{\pi i y\cdot x}
		\int_{\mathbb R^{d-1}}
		W_\lambda(s)\overline{\psi_0(s-x)}e^{-2\pi i y\cdot s}\,ds
		\right)
		e^{-\pi i y\cdot x}e^{2\pi i y\cdot t}\psi_0(t-x)\\
		=&
		\left(
		\int_{\mathbb R^{d-1}}
		W_\lambda(s)\overline{\psi_0(s-x)}e^{-2\pi i y\cdot s}\,ds
		\right)
		e^{2\pi i y\cdot t}\psi_0(t-x).
	\end{aligned}
	\]
	Hence, for \(t\in\mathbb R^{d-1}\), by \cite[Lemma~1.4(ii)]{GrochenigZimmermann}, 
	\begin{equation}\label{Wint}
	\begin{aligned}
		&W_\lambda(t)\\
		=&\langle\psi_0,\psi_0\rangle W_\lambda(t)\\
		=&\iint_{\mathbb R^{d-1}\times\mathbb R^{d-1}}
		\left(
		\int_{\mathbb R^{d-1}}
		W_\lambda(s)\,\overline{\psi_0(s-x)}\,
		e^{-2\pi i y\cdot s}\,ds
		\right)
		e^{2\pi i y\cdot t}\psi_0(t-x)\,dx\,dy\\
		=&\iint_{\mathbb R^{d-1}\times\mathbb R^{d-1}} \langle W_\lambda,\psi_u\rangle\psi_u\,du.
	\end{aligned}
	\end{equation}
	Therefore, by \eqref{Wlambda_orth} and \eqref{Wint},
	\begin{equation}\label{v-delta}
		\begin{aligned}
			&\langle v_\lambda,v_\mu\rangle-\delta_{\lambda\mu}\\
			=&\langle v_\lambda,v_\mu\rangle-\langle W_\lambda,W_\mu\rangle\\
			=&\iint_{\mathbb R^{4(d-1)}}
			\langle W_\lambda,\psi_u\rangle
			\overline{\langle W_\mu,\psi_v\rangle}
			\Bigl(
			\mathbf 1_{\Omega_*\times\Omega_*}(u,v)
			\langle b'_{A(u)},b'_{A(v)}\rangle
			-\langle\psi_u,\psi_v\rangle
			\Bigr)\,du\,dv,
		\end{aligned}
	\end{equation}
	where $\mathbf{1}_{\Omega_*\times\Omega_*}$ denotes the indicator
	function of $\Omega_*\times\Omega_*$.
	Set $k=M+2(d-1)+1$.
	Then, by \eqref{v-delta}, Lemma~\ref{lem:gaussian-packet-inner-product}, Lemma~\ref{lem:wq-gram-model} and Lemma~\ref{lem:wilson-packet-decay}, there exists \(C_{d,k,s}>0\) such that
	\begin{equation}\label{|v-delta|}
	\begin{aligned}
		\bigl|\langle v_\lambda,v_\mu\rangle-\delta_{\lambda\mu}\bigr|
		&\leq C_{d,k,s}
		\sum_{\xi\in O_\lambda}\sum_{\eta\in O_\mu}
		\iint_{\mathbb R^{4(d-1)}}
		\frac{e^{-\frac{\pi|u-v|^2}{2}}}
		{(1+|u-\xi|)^{k+s}(1+|v-\eta|)^{k+s}}\\
		&\quad\times\Bigl(
		e_n\mathbf 1_{\Omega_*\times\Omega_*}(u,v)
		+\mathbf 1_{(\Omega_*\times\Omega_*)^c}(u,v)
		\Bigr)\,du\,dv .
	\end{aligned}
	\end{equation}
	If \((u,v) \notin \Omega^* \times \Omega^*\), then $|u-\xi|\geq B$ or $|v-\eta|\geq B$. 
	Thus
	\begin{align}\label{bu}
	\frac{1}{(1+|u-\xi|)^{k+s}(1+|v-\eta|)^{k+s}}
	\leq
	\frac{B^{-s}}{(1+|u-\xi|)^k(1+|v-\eta|)^k}.
	\end{align}	
	Since \(1+|\xi-\eta|\leq(1+|u-\xi|)(1+|u-v|)(1+|v-\eta|)\), we have
	\begin{equation}\label{2eq}
	  \begin{aligned}
			&\iint_{\mathbb R^{4(d-1)}}
			\frac{e^{-\frac{\pi|u-v|^2}{2}}}
			{(1+|u-\xi|)^k(1+|v-\eta|)^k}\,du\,dv\\
			\leq&
			\frac{1}{(1+|\xi-\eta|)^M}
			\iint_{\mathbb R^{4(d-1)}}
			\frac{(1+|u-v|)^M e^{-\frac{\pi|u-v|^2}{2}}}
			{(1+|u-\xi|)^{k-M}(1+|v-\eta|)^{k-M}}\,du\,dv\\
			\leq&
			C_{M,k}{(1+|\xi-\eta|)^{-M}},
		\end{aligned}
	\end{equation}
	where 
	\[
	C_{M,k}= \sup\limits_{r \geq0} (1+r)^M e^{-\frac{\pi r^2}{2}}
	\left(
	\int_{\mathbb R^{2(d-1)}}
	(1+|w|)^{-(k-M)}\,dw
	\right)^2<\infty.
	\]
	Consequently, by \eqref{|v-delta|}, \eqref{bu} and \eqref{2eq}, 
	\begin{align*}
		&\bigl|\langle v_\lambda,v_\mu\rangle-\delta_{\lambda\mu}\bigr|\\
		\leq& C_{d,k,s}(e_n+B^{-s})
		\sum_{\xi\in O_\lambda}\sum_{\eta\in O_\mu}
		\iint_{\mathbb R^{4(d-1)}}
		\frac{e^{-\frac{\pi|u-v|^2}{2}}}
		{(1+|u-\xi|)^k(1+|v-\eta|)^k}\,du\,dv\\
		\leq& C_{d,k,s}C_{M,k}(e_n+B^{-s})
		\sum_{\xi\in O_\lambda}\sum_{\eta\in O_\mu}
		(1+|\xi-\eta|)^{-M}\\
		\leq& C_{d,k,s}C_{M,k} 2^{2(d-1)+\frac M2}(e_n+B^{-s})
		(1+\rho_*(\lambda,\mu))^{-M},
	\end{align*}
	where \((*)\) is given by
	\[
	|\xi-\eta|
	\geq \sqrt{|q-q'|^2+|\kappa-\kappa'|^2}
	\geq \frac{|q-q'|+|\kappa-\kappa'|}{\sqrt2}
	=\frac{\rho_*(\lambda,\mu)}{\sqrt2}.
	\]
	Taking $C_{d,M,s}=C_{d,k,s}C_{M,k} 2^{2(d-1)+\frac M2}$ gives \eqref{eq:wq-local-gram}.

	For $u\in\Omega_*$, $\xi=(q,\sigma\kappa)\in O_\lambda$ and
	$x\in S^d$, by Lemma~\ref{lem:grassmannian-comparisons} and Lemma~\ref{lem:graph-coordinate-distance},
	\[
	\begin{aligned}
		\bigl|s_{E_{\lambda,\sigma}}(x)-s_{A(u)}(x)\bigr|
		&\le 2\sqrt2\,d_Q\bigl(E_{A(\xi)},E_{A(u)}\bigr)\\
		&\le 2\sqrt2\cdot\frac{3\sqrt2\,\pi}{2}\,|A(\xi)-A(u)|\\
		&=6\pi\sqrt{\frac{2\pi}{n}}\,|u-\xi|.
	\end{aligned}
	\]
	Therefore, for  $k=M+2(d-1)+1$, we have
	\begin{equation}\label{M}
	\begin{aligned}
		&\bigl(1+\sqrt n\,s_{E_{\lambda,\sigma}}(x)\bigr)^M
		\int_{\Omega_*}
		\frac{\exp\!\left(-\frac{n\,s_{A(u)}(x)^2}{2}\right)}
		{(1+|u-\xi|)^k}\,du\\
		\le&
		\int_{\Omega_*}
		\frac{\bigl(1+\sqrt n\,s_{A(u)}(x)+6\pi\sqrt{2\pi}\,|u-\xi|\bigr)^M
			\exp\!\left(-\frac{n\,s_{A(u)}(x)^2}{2}\right)}
		{(1+|u-\xi|)^k}\,du\\
		\le&(6\pi\sqrt{2\pi})^M
		\int_{\Omega_*}
		\frac{\bigl(1+\sqrt n\,s_{A(u)}(x)\bigr)^M
			\exp\!\left(-\frac{n\,s_{A(u)}(x)^2}{2}\right)}
		{(1+|u-\xi|)^{k-M}}\,du\\
		\le& C_{d,k,M}
		<\infty,
	\end{aligned}
	\end{equation}
	where
	\[
	C_{d,k,M}=(6\pi\sqrt{2\pi})^M
	\sup_{r\ge 0}\left((1+r)^M\exp\!\left(-\frac{r^2}{2}\right)\right)
	\int_{\mathbb R^{2(d-1)}}
	\frac{1}{(1+|w|)^{k-M}}\,dw.
	\]
	Then, by Lemma~\ref{lemma2.3}(i)(c) and Lemma~\ref{lem:wilson-packet-decay}, there exists \(C_{d,k}>0\) such that
	\begin{align*}
		|v_\lambda(x)|
		&\le \int_{\Omega_*}
		|\langle W_\lambda,\psi_u\rangle|\,
		|b'_{A(u)}(x)|\,du\\
		&\leq C_{d,k}\sqrt{\frac{d+1}{2}}n^{\frac{d-1}{4}}
		\sum_{\sigma\in\{-1,1\}^{d-1}}
		\int_{\Omega_*}
		\frac{\exp\left(-\frac{n s_{A(u)}(x)^2}{2}\right)}
		{(1+|u-(q,\sigma\kappa)|)^k}\,du \\
		&\leq  C_{d,k}\sqrt{\frac{d+1}{2}}
		C_{d,k,M}\ n^{\frac{d-1}{4}}
		\sum_{\sigma\in\{-1,1\}^{d-1}}
		\bigl(1+\sqrt n\,s_{E_{\lambda,\sigma}}(x)\bigr)^{-M}.
	\end{align*}
	Taking $C_{d,M}=C_{d,k}\sqrt{\frac{d+1}{2}}
	C_{d,k,M}$ gives \eqref{eq:wq-local-field}.
\end{proof}

Let $\Omega$ and $\widetilde\Omega$ be graph cubes with centers
$A_0$ and $\widetilde A_0$, respectively. Write $\Omega_*$ and
$\widetilde\Omega_*$ for their parameter cubes. 
For $\lambda=(q,\kappa)\in\Lambda_B(\Omega)$ and
$\mu=(\widetilde q, \widetilde{\kappa})\in\Lambda_B(\widetilde\Omega)$, let \(v_\lambda\) and \(\widetilde v_\mu\) be defined by \eqref{def_vlambda}.
For each $\sigma,\tau\in\{-1,1\}^{d-1}$, let $E_{\lambda,\sigma}=E_{A(q,\sigma\kappa)}$ and $\widetilde E_{\mu,\tau}
=E_{\widetilde A(\widetilde q,\tau\widetilde\kappa)}$.
Define $\mathcal C$ by $(\mathcal C f)(x)=\overline{f(x)}$. 
For rotations $R$ and $S$ and $\varepsilon,\delta\in\{0,1\}$, set
\begin{align}\label{def_UV}
U=T_R\mathcal C^\varepsilon,
\qquad
V=T_S\mathcal C^\delta,
\end{align}
where \((T_Rf)(x)=f(R^{-1}x)\) and \((T_Sf)(x)=f(S^{-1}x)\).
For $E=\operatorname{span}^{+}(e,f)\in Q$, set $T_RE=RE$, $T_SE=SE$ and
$\mathcal CE=\overline E=\operatorname{span}^{+}(e,-f)$.

\begin{lemma}\label{lem:cross-packet-estimate}
	For every $M>0$, there exists a constant $C_{d,M}>0$,
	independent of $n$, $\Omega$, $\widetilde\Omega$, $\lambda$, $\mu$, $U$ and $V$, such that
	\begin{equation*}
		\bigl|\langle Uv_\lambda,V\widetilde v_\mu\rangle\bigr|
		\leq C_{d,M}
		\sum_{\sigma,\tau\in\{-1,1\}^{d-1}}
		\frac{1}{
			\bigl(1+\sqrt n\,d_Q(UE_{\lambda,\sigma},
			V\widetilde E_{\mu,\tau})\bigr)^M}.
	\end{equation*}
\end{lemma}

\begin{proof}
	Fix $k=M+2(d-1)+1$. For $u\in\Omega_*$ and
	$v\in\widetilde\Omega_*$, set
	\[
	D_{uv}=\sqrt n\,d_Q\bigl(UE_{A(u)},VE_{\widetilde A(v)}\bigr).
	\] 
	Let $z$ and $w$ be unit representatives of $UE_{A(u)}$ and $VE_{\widetilde A(v)}$, respectively. 
	We claim that there exist \(\alpha, \beta \in \mathbb C\) with \(|\alpha|=|\beta|=1\) such that $Ub'_{A(u)}=\alpha b_z$ and
	$Vb'_{\widetilde A(v)}=\beta b_w$.
	In fact,
	since \(z\) and $R\mathcal C^\varepsilon z'(A(u))$ are unit
	representatives of $UE_{A(u)}$, we have
	$R\mathcal C^\varepsilon z'(A(u))=e^{\mathrm i\theta}z$ for some $\theta\in\mathbb R$.
	Therefore, by the definition of $b'_{A(u)}$,
	\[
	Ub'_{A(u)}
	=T_R\mathcal C^\varepsilon b_{z'(A(u))}
	=b_{R\mathcal C^\varepsilon z'(A(u))}
	=b_{e^{\mathrm i\theta}z}
	=e^{-\mathrm i n\theta}b_z.
	\]
	Setting $\alpha=e^{-\mathrm i n\theta}$, we obtain
	$Ub'_{A(u)}=\alpha b_z$ with $|\alpha|=1$.
	Similarly, $Vb'_{\widetilde A(v)}=\beta b_w$ with $|\beta|=1$.
	The claim is proved.
	Consequently, by Lemma~\ref{lemma2.3}(ii)(b),
	\begin{align}\label{UU}
	\bigl|\langle Ub'_{A(u)},Vb'_{\widetilde A(v)}\rangle\bigr|
	=\bigl|\langle b_z,b_w\rangle\bigr|
	=\left(\cos\!\left(\frac{D_{uv}}{\sqrt n}\right)\right)^n
	\leq \exp\!\left(-\frac{D_{uv}^2}{2}\right).
	\end{align}
	
	By the definitions of $v_\lambda$ and $\widetilde v_\mu$,
	\eqref{UU} and Lemma~\ref{lem:wilson-packet-decay},
	\begin{equation}\label{<u,v>}
	\begin{aligned}
		&\bigl|\langle Uv_\lambda,V\widetilde v_\mu\rangle\bigr|\\
		=&
		\left|
		\iint_{\Omega_*\times\widetilde\Omega_*}
		\Big\langle
		U\bigl(\langle W_\lambda,\psi_u\rangle b'_{A(u)}\bigr),
		V\bigl(\langle W_\mu,\psi_v\rangle b'_{\widetilde A(v)}\bigr)
		\Big\rangle\,du\,dv
		\right|\\
		\leq&
		\iint_{\Omega_*\times\widetilde\Omega_*}
		\bigl|\langle W_\lambda,\psi_u\rangle\bigr|
		\bigl|\langle W_\mu,\psi_v\rangle\bigr|
		\left|\langle Ub'_{A(u)}, Vb'_{\widetilde A(v)} \rangle\right|\,du\,dv\\
		\leq&
		\iint_{\Omega_*\times\widetilde\Omega_*}
		\bigl|\langle W_\lambda,\psi_u\rangle\bigr|
		\bigl|\langle W_\mu,\psi_v\rangle\bigr|
		\exp\left(-\frac{D_{uv}^2}{2}\right)\,du\,dv\\
		\leq& C_{d,k}
		\sum_{\xi\in O_\lambda}\sum_{\eta\in O_\mu}
		\iint_{\Omega_*\times\widetilde\Omega_*}
		\frac{\exp\left(-\frac{D_{uv}^2}{2}\right)}
		{(1+|u-\xi|)^k(1+|v-\eta|)^k}\,du\,dv.
	\end{aligned}
	\end{equation}
	For $\xi\in O_\lambda\subset\Omega_*$ and
	$\eta\in O_\mu\subset\widetilde\Omega_*$, by Lemma~\ref{lem:graph-coordinate-distance},
	\begin{align*}
		D_{\xi\eta}
		&=\sqrt n\,d_Q\bigl(UE_{A(\xi)},VE_{\widetilde A(\eta)}\bigr)\\
		&\leq \sqrt n\,d_Q\bigl(E_{A(\xi)},E_{A(u)}\bigr)
		+D_{uv}
		+\sqrt n\,d_Q\bigl(E_{\widetilde A(v)},E_{\widetilde A(\eta)}\bigr)\\
		&\leq D_{uv}
		+\frac{3\sqrt2\,\pi}{2}\sqrt n
		\bigl(|A(\xi)-A(u)|
		+|\widetilde A(v)-\widetilde A(\eta)|\bigr)\\
		&=D_{uv}+3\pi\sqrt{\pi}\bigl(|u-\xi|+|v-\eta|\bigr).
	\end{align*}
	Therefore,
	\begin{align*}
	1+D_{\xi\eta}
	\leq 3\pi\sqrt{\pi}\,
	(1+|u-\xi|)(1+D_{uv})(1+|v-\eta|).
	\end{align*}
	It follows that
	\begin{equation}\label{iint}
	\begin{aligned}
		&\iint_{\Omega_*\times\widetilde\Omega_*}
		\frac{\exp\left(-\frac{D_{uv}^2}{2}\right)}
		{(1+|u-\xi|)^k(1+|v-\eta|)^k}\,du\,dv\\
		\leq&
		\frac{(3\pi\sqrt{\pi})^M}{(1+D_{\xi\eta})^M}
		\iint_{\Omega_*\times\widetilde\Omega_*}
		\frac{(1+D_{uv})^M\exp\left(-\frac{D_{uv}^2}{2}\right)}
		{(1+|u-\xi|)^{k-M}(1+|v-\eta|)^{k-M}}\,du\,dv\\
		\leq&
		\frac{(3\pi\sqrt{\pi})^M}{(1+D_{\xi\eta})^M}
		\sup_{r\geq0}\left((1+r)^M\exp\left(-\frac{r^2}{2}\right)\right)
		\left(
		\int_{\mathbb R^{2(d-1)}}\frac{1}{(1+|w|)^{k-M}}\,dw
		\right)^2\\
		\leq& \frac{C_{d,M}'}{(1+D_{\xi\eta})^M},
	\end{aligned}
	\end{equation}
	where 
	\[
	C_{d,M}' = (3\pi\sqrt{\pi})^M 	\sup_{r\geq0}(1+r)^M\exp\left(-\frac{r^2}{2}\right)
	\left(
	\int_{\mathbb R^{2(d-1)}}\frac{1}{(1+|w|)^{k-M}}\,dw
	\right)^2.
	\]
	Consequently, by \eqref{<u,v>} and \eqref{iint}, we have
	\[
	\begin{aligned}
		\bigl|\langle Uv_\lambda,V\widetilde v_\mu\rangle\bigr|
		&\leq C_{d,k}C_{d,M}'
		\sum_{\xi\in O_\lambda}\sum_{\eta\in O_\mu}
		\frac{1}{(1+D_{\xi\eta})^M}\\
		&= C_{d,k}C_{d,M}'
		\sum_{\sigma,\tau\in\{-1,1\}^{d-1}}
		\frac{1}{
			\bigl(1+\sqrt n\,d_Q(UE_{\lambda,\sigma},
			V\widetilde E_{\mu,\tau})\bigr)^M}.
	\end{aligned}
	\]
	Taking $C_{d,M}=C_{d,k}
	C_{d,M}'$ gives Lemma~\ref{lem:cross-packet-estimate}.
\end{proof}
The following lemma follows from \cite[Theorem~1.1]{Kirwin2010}.
\begin{lemma}\label{lem:multidimensional-laplace}
	Let $\mathcal R\subset\mathbb R^m$ be a bounded open set containing
	$0$. Let $f$ and $g$ be real-valued functions of class $C^6$ and
	$C^4$, respectively, on a neighborhood of $\overline{\mathcal R}$.
	Suppose that $f(0)=0$, that $f(x)>0$ for all
	$x\in\overline{\mathcal R}\setminus\{0\}$, and that
	$\nabla^2f(0)$ is positive definite. Then, as $k\to\infty$,
	\[
	\int_{\mathcal R}e^{-kf(x)}g(x)\,dx
	=
	\left(\frac{2\pi}{k}\right)^{\frac{m}{2}}
	\frac{g(0)}{\sqrt{\det\nabla^2f(0)}}
	+O\!\left(\frac{1}{k^{\frac m2+1}}\right).
	\]
\end{lemma}
\begin{lemma}\label{lem:local-packet-count}
	Let $N_\Omega=\#\Lambda_B(\Omega)$. Then
	\begin{equation*}
		N_\Omega
		=D_n\nu(\Omega)
		\left(1+O_d\!\left(
		R+\frac 1n+\frac{B+1}{\sqrt n\,\ell}
		\right)\right).
	\end{equation*}
\end{lemma}
\begin{proof}
	Let $\rho(A)$ be the smooth positive density of $\nu$ in the
	graph coordinates $A\mapsto E_A$, so that $d\nu=\rho(A)\,dA$.
	Fix a sufficiently small $r>0$ and let $\mathcal V=\{E_A:|A|<r\}$.
	Since $\mathcal V$ is a neighborhood of $E_0$, \(\delta=\inf\limits_{E\in Q\setminus\mathcal V}d_Q(E_0,E)>0\).
	By Lemma~\ref{lemma2.3}(ii)(b),
	\begin{equation}\label{nu}
		\begin{aligned}
			\int_{Q\setminus\mathcal V}
			|\langle b_0,b_z\rangle|^2\,d\nu([z])
			\leq \int_{Q\setminus\mathcal V}e^{-n\delta^2}\,d\nu([z])
			=e^{-n\delta^2}\nu(Q\setminus\mathcal V)
			\leq e^{-n\delta^2}.
		\end{aligned}
	\end{equation}
	Set $\Phi(A)=-2\operatorname{Re}F_0(0,A)$.
	By Lemma~\ref{lemma2.3}(ii)(a) and the definition of $F_0$,
	\begin{equation}\label{b0bA}
		\begin{aligned}
			|\langle b_0,b_A\rangle|^2
			&=|\langle b'_0,b'_A\rangle|^2\\
			&=\bigl|z'(A)\cdot\overline{z'(0)}\bigr|^{2n}\\
			&=\bigl|\exp(F_0(0,A))\bigr|^{2n}\\
			&=\exp\!\bigl(2n\operatorname{Re}F_0(0,A)\bigr)\\
			&=\exp\!\left(-n\Phi(A)\right).
		\end{aligned}
	\end{equation}
	For $0<|A|$ sufficiently small, Lemma~\ref{lem:local-log-expansion} with
	$A_0=0$, $R=\ell=|A|$, $X=0$, and $Y=A$ gives
	\begin{align}\label{Fzhankai}
		F_0(0,A)
		=-\frac{|A|^2}{4}+O_d(R\ell^2+\ell^3)
		=-\frac{|A|^2}{4}+O_d(|A|^3).
	\end{align}
	Hence, by \eqref{Fzhankai}, we have
	\[
	\Phi(A)=\frac{|A|^2}{2}+O_d(|A|^3).
	\]
	In particular, $\Phi(0)=0$ and $\nabla^2\Phi(0)=I$.
	For $r$ sufficiently small, $\Phi$ and $\rho$ are smooth near
	$\{|A|\leq r\}$ and $\Phi(A)>0$ for $0<|A|\leq r$.
	Thus, by Lemma~\ref{lem:multidimensional-laplace} with $m=2(d-1)$, we have
	\begin{align}\label{lemma2.12}
	\int_{|A|<r}e^{-n\Phi(A)}\rho(A)\,dA
	=\rho(0)\left(\frac{2\pi}{n}\right)^{d-1}
	+O_d\!\left(\frac{1}{n^d}\right).
	\end{align}
	Therefore, by Lemma~\ref{lemma2.3}(i)(b) and (i)(d), \eqref{nu}, \eqref{b0bA} and \eqref{lemma2.12},
	\begin{align*}
		\frac{1}{D_n}
		&=\frac{\langle b_0,b_0\rangle}{D_n}\\
		&=\left\langle
		\left(\int_Q b_z\otimes b_z\,d\nu([z])\right)b_0,
		b_0
		\right\rangle\\
		&=\int_{\mathcal V}|\langle b_0,b_z\rangle|^2\,d\nu([z])
		+\int_{Q\setminus\mathcal V}
		|\langle b_0,b_z\rangle|^2\,d\nu([z])\\
		&=\int_{|A|<r}|\langle b_0,b_A\rangle|^2\rho(A)\,dA
		+O_d\left(e^{-n\delta^2}\right)\\
		&=\int_{|A|<r}e^{-n\Phi(A)}\rho(A)\,dA
		+O_d\left(e^{-n\delta^2}\right)\\
		&=\rho(0)\left(\frac{2\pi}{n}\right)^{d-1}
		+O_d\!\left(\frac{1}{n^d}\right)
		+O_d\left(e^{-n\delta^2}\right)\\
		&=\rho(0)\left(\frac{2\pi}{n}\right)^{d-1}
		\left(1+O_d\!\left(\frac{1}{n}\right)\right).
	\end{align*}
	Hence,
	\begin{equation}\label{Dnrho}
	\begin{aligned}
		D_n\rho(0)
		&=\left(\frac{n}{2\pi}\right)^{d-1}
		\frac{1}{1+O_d(\frac1n)}\\
		&=\left(\frac{n}{2\pi}\right)^{d-1}
		\left(1-\frac{O_d(\frac1n)}{1+O_d(\frac1n)}\right)\\
		&=\left(\frac{n}{2\pi}\right)^{d-1}
		\left(1+O_d\!\left(\frac{1}{n}\right)\right).
	\end{aligned}
	\end{equation}
For $A\in\Omega$, we have $|A_0|\leq R$ and
$|A_{ij}-(A_0)_{ij}|\leq\ell$ for all $i,j$. Since $\ell\leq R$,
\[
|A|\leq |A_0|+|A-A_0|
\leq R+\sqrt{2(d-1)}\,\ell
=O_d(R).
\]
Thus by the smoothness of $\rho$, we have
\begin{align}\label{rhoA}
\frac{\rho(A)}{\rho(0)}
=1+O_d(|A|)
=1+O_d(R).
\end{align}
Therefore, by \eqref{rhoA} and \eqref{Dnrho},
\begin{equation}\label{Dnnu}
	\begin{aligned}
		D_n\nu(\Omega)
		&=D_n\int_\Omega \rho(A)\,dA\\
		&=D_n\rho(0)(2\ell)^{2(d-1)}\bigl(1+O_d(R)\bigr)\\
		&=\left(2\sqrt{\frac{n}{2\pi}}\,\ell\right)^{2(d-1)}
		\left(1+O_d\!\left(R+\frac{1}{n}\right)\right).
	\end{aligned}
\end{equation}
	For $u=(u_1,\ldots,u_{2(d-1)})\in\Omega_*$, it is easy to see that
	\[
	\operatorname{dist}(u,\partial\Omega_*)
	=\sqrt{\frac{n}{2\pi}}\,\ell-\max_{1\le j\le 2(d-1)}|u_j|.
	\]
	Hence $\lambda=(q,\kappa)\in\Lambda_B(\Omega)$ if and only if,
	for each $1\leq i\leq d-1$,
	\[
	q_i\in\frac12\mathbb Z,\quad
	|q_i|\leq\sqrt{\frac{n}{2\pi}}\,\ell-B,
	\quad
	\kappa_i\in\mathbb Z,\quad
	1\leq\kappa_i\leq\sqrt{\frac{n}{2\pi}}\,\ell-B.
	\]
	Counting the admissible values of $2q_i$ and $\kappa_i$ gives
	\begin{equation}\label{Nomega}
		\begin{aligned}
			N_\Omega
			&=\left[
			\left(2\left\lfloor 2\sqrt{\frac{n}{2\pi}}\,\ell-2B
			\right\rfloor+1\right)
			\left\lfloor\sqrt{\frac{n}{2\pi}}\,\ell-B\right\rfloor
			\right]^{d-1}\\
			&=\left[
			\left(4\sqrt{\frac{n}{2\pi}}\,\ell+O(B+1)\right)
			\left(\sqrt{\frac{n}{2\pi}}\,\ell+O(B+1)\right)
			\right]^{d-1}\\
			&=
			\left(2\sqrt{\frac{n}{2\pi}}\,\ell\right)^{2(d-1)}
			\left[
			\left(1+O\!\left(\frac{B+1}{\sqrt n\,\ell}\right)\right)
			\left(1+O\!\left(\frac{B+1}{\sqrt n\,\ell}\right)\right)
			\right]^{d-1}\\
			&=
			\left(2\sqrt{\frac{n}{2\pi}}\,\ell\right)^{2(d-1)}
			\left[
			1+O\!\left(\frac{B+1}{\sqrt n\,\ell}
			+\frac{(B+1)^2}{n\ell^2}\right)
			\right]^{d-1}\\
			&=
			\left(2\sqrt{\frac{n}{2\pi}}\,\ell\right)^{2(d-1)}
			\left(1+O_d\!\left(\frac{B+1}{\sqrt n\,\ell}\right)\right),
		\end{aligned}
	\end{equation}
	where the last equality uses the fact that
	$\frac{B+1}{\sqrt n\,\ell}=O(1)$, which follows from
	$1\ll B<\frac12\sqrt{\frac{n}{2\pi}}\,\ell$.
	Consequently, by \eqref{Dnnu} and \eqref{Nomega},
	\begin{align*}
		\frac{N_\Omega}{D_n\nu(\Omega)}
		=\frac{1+O_d\!\left(\frac{B+1}{\sqrt n\,\ell}\right)}
		{1+O_d\!\left(R+\frac{1}{ n}\right)}
		=1+O_d\!\left(R+\frac{1}{ n}
		+\frac{B+1}{\sqrt n\,\ell}\right).
	\end{align*}
\end{proof}

For each $\Omega$ and $a\geq0$, we define the
weighted Schur norm of a matrix
$H=(H_{\lambda\mu})_{\lambda,\mu\in\Lambda_B(\Omega)}$ by
\[
\|H\|_a
=\max\left\{
\sup_\lambda\sum_\mu
(1+\rho_*(\lambda,\mu))^a|H_{\lambda\mu}|,
\ \sup_\mu\sum_\lambda
(1+\rho_*(\lambda,\mu))^a|H_{\lambda\mu}|
\right\}.
\]
This is the finite-matrix version of the norm in \cite[Definition~1]{GrochenigLeinert}.
By the triangle inequality for $\rho_*$,
\[
1+\rho_*(\lambda,\mu)
\leq (1+\rho_*(\lambda,\alpha))
(1+\rho_*(\alpha,\mu)).
\]
Thus for matrices $H$ and $K$, by \cite[Lemma~1(a)]{GrochenigLeinert}, we have
\begin{align}\label{aleqa}
\|HK\|_a\leq\|H\|_a\|K\|_a.
\end{align}
For $c\in\ell^2(\Lambda_B(\Omega))$, by the Cauchy--Schwarz
inequality, we have
\begin{align*}
	\|Hc\|_{\ell^2}^2
	&=\sum_\lambda
	\left|\sum_\mu H_{\lambda\mu}c_\mu\right|^2\\
	&\leq \sum_\lambda
	\left(\sum_\mu |H_{\lambda\mu}|\right)
	\left(\sum_\mu |H_{\lambda\mu}|\,|c_\mu|^2\right)\\
	&\leq \|H\|_0\sum_\mu |c_\mu|^2
	\sum_\lambda |H_{\lambda\mu}|\\
	&\leq \|H\|_0^2\|c\|_{\ell^2}^2.
\end{align*}
Thus 
\begin{align}\label{2to2}
	\|H\|_{2\to2}=\sup_{\|c\|_{\ell^2}=1}\|Hc\|_{\ell^2}\leq\|H\|_0.
\end{align}

\begin{lemma}\label{lem:local-packet-orthonormalization}
	Suppose that $e_n=n(R\ell^2+\ell^3)\to0$ and $B\to\infty$
	as $n\to\infty$. Then the following hold.
	\begin{enumerate}
		\item[\textup{(i)}] For all sufficiently large $n$ and each $\Omega$,
		there exists an orthonormal basis
		$\{w_\lambda\}_{\lambda\in\Lambda_B(\Omega)}$ of
		$\operatorname{span}\{v_\lambda:\lambda\in\Lambda_B(\Omega)\}$.
		
		\item[\textup{(ii)}] For every fixed $M>0$ and all sufficiently
		large $n$, the basis in \textup{(i)} satisfies
		\begin{equation}\label{eq:wq-orthonormal-field}
			|w_\lambda(x)|
			\leq C'_{d,M}n^{\frac{d-1}{4}}
			\sum_{\sigma\in\{-1,1\}^{d-1}}
			\bigl(1+\sqrt n\,s_{E_{\lambda,\sigma}}(x)\bigr)^{-M}
		\end{equation}
		for every $\Omega$, $\lambda\in\Lambda_B(\Omega)$ and $x\in S^d$,
		where $C'_{d,M}>0$ is independent of $n$, $\Omega$, $\lambda$ and $x$.
	\end{enumerate}
\end{lemma}
\begin{proof}
	Let $\{e_\lambda\}_{\lambda\in\Lambda_B(\Omega)}$ be the standard
	basis of $\ell^2(\Lambda_B(\Omega))$. Define the linear map
	$T_\Omega:\ell^2(\Lambda_B(\Omega))
	\to\mathcal H_n(S^d)\otimes_{\mathbb R}\mathbb C$
	by $T_\Omega e_\lambda=v_\lambda$.
	Set $G_\Omega=T_\Omega^*T_\Omega$. Then
	\begin{align}\label{Glambdanu}
		(G_\Omega)_{\lambda\mu}
		=\langle T_\Omega^*T_\Omega e_\mu,e_\lambda\rangle
		=\langle v_\mu,v_\lambda\rangle.
	\end{align}
	Since $\Lambda_B(\Omega)\subset
	\left(\tfrac12\mathbb Z\right)^{d-1}
	\times\mathbb Z_{\geq1}^{d-1}$, for every $\lambda\in\Lambda_B(\Omega)$
	and $t\geq0$, we have
	\begin{equation}\label{number_nu}
		\begin{aligned}
			\#\{\mu\in\Lambda_B(\Omega):\rho_*(\lambda,\mu)\leq t\}
			&\leq (2\lfloor 2t\rfloor+1)^{d-1}
			(2\lfloor t\rfloor+1)^{d-1}\\
			&\leq (4t+1)^{d-1}(2t+1)^{d-1}\\
			&\leq C_d(1+t)^{2(d-1)}.
		\end{aligned}
	\end{equation}
	Fix $a\geq0$ and $\lambda\in\Lambda_B(\Omega)$.
	By \eqref{Glambdanu}, \eqref{eq:wq-local-gram} with
	$M=a+2(d-1)+2$ and $s=1$, and \eqref{number_nu} with $t=j+1$,
	we have
	\begin{align*}
		&\sum_{\mu\in\Lambda_B(\Omega)}
		(1+\rho_*(\lambda,\mu))^a
		|(G_\Omega-I)_{\lambda\mu}|\\
		\leq& C_{d,a}\left(e_n+\frac{1}{B}\right)
		\sum_{j=0}^{\infty}
		\sum_{\substack{\mu\in\Lambda_B(\Omega)\\
				j\leq\rho_*(\lambda,\mu)<j+1}}
		(1+\rho_*(\lambda,\mu))^{-2(d-1)-2}\\
		\leq& C_{d,a}C_d\left(e_n+\frac{1}{B}\right)
		\sum_{j=0}^{\infty}(2+j)^{2(d-1)}(1+j)^{-2(d-1)-2}\\
		\leq& \frac{\pi^2}{6}\,2^{2(d-1)}C_{d,a}C_d
		\left(e_n+\frac{1}{B}\right).
	\end{align*}
	The same argument, with $\mu$ fixed and the sum taken over
	$\lambda$, gives the corresponding bound for the weighted
	column sums. Hence
	\begin{equation}\label{Ga}
		\|G_\Omega-I\|_a
		\leq C_{d,a}'\left(e_n+\frac{1}{B}\right),
	\end{equation}
	where $C_{d,a}'=\frac{\pi^2}{6}\,2^{2(d-1)}C_{d,a}C_d$.
	
	(i)
	Since $e_n+\frac{1}{B}\to0$, applying \eqref{Ga} with $a=0$
	gives, for all sufficiently large $n$, \(\|G_\Omega-I\|_0\leq\frac12\).
	Thus by \eqref{2to2}, we have $G_\Omega\geq\frac12 I$.
	For each $\lambda\in\Lambda_B(\Omega)$, define
	\begin{align}\label{def_wlambda}
	w_\lambda=T_\Omega G_\Omega^{-\frac12}e_\lambda
	=\sum_{\alpha\in\Lambda_B(\Omega)}
	\bigl(G_\Omega^{-\frac12}\bigr)_{\alpha\lambda}v_\alpha.
	\end{align}
	Since $G_\Omega^{-\frac12}$ is self-adjoint and
	$T_\Omega^*T_\Omega=G_\Omega$, we have
	\begin{align*}
		\langle w_\mu,w_\lambda\rangle
		=\left\langle
		G_\Omega^{-\frac12}G_\Omega G_\Omega^{-\frac12}e_\mu,
		e_\lambda
		\right\rangle
		=\langle e_\mu,e_\lambda\rangle
		=\delta_{\lambda\mu}.
	\end{align*}
	Moreover, the invertibility of $G_\Omega^{-\frac12}$ gives
	$\operatorname{ran}(T_\Omega G_\Omega^{-\frac12})
	=\operatorname{ran}T_\Omega$. Hence
	$\{w_\lambda\}_{\lambda\in\Lambda_B(\Omega)}$ is an orthonormal
	basis of $\operatorname{span}\{v_\lambda:
	\lambda\in\Lambda_B(\Omega)\}$.\\
	
	(ii)
	Fix $M>0$. Since $e_n+B^{-1}\to0$, \eqref{Ga} with $a=M$ gives
	\begin{equation}\label{a}
		\|G_\Omega-I\|_M\leq\frac12
	\end{equation}
	for all sufficiently large $n$. By the binomial series,
	\eqref{aleqa}, and \eqref{a},
	\begin{equation}\label{eq:wq-orthonormal-coeff}
		\begin{aligned}
			\|G_\Omega^{-\frac12}\|_M
			\leq 1+\sum_{r=1}^{\infty}
			\left|\binom{-\frac12}{r}\right|
			\|G_\Omega-I\|_M^r
			=\bigl(1-\|G_\Omega-I\|_M\bigr)^{-\frac12}
			\leq\sqrt{2}.
		\end{aligned}
	\end{equation}
	Fix $x\in S^d$. Let $\lambda=(q,\kappa)$ and
	$\alpha=(q',\kappa')$ belong to $\Lambda_B(\Omega)$.
	For each $\sigma\in\{-1,1\}^{d-1}$, by
	Lemmas~\ref{lem:grassmannian-comparisons}
	and~\ref{lem:graph-coordinate-distance},
	\[
	\begin{aligned}
		\sqrt n\,
		\bigl|s_{E_{\lambda,\sigma}}(x)-s_{E_{\alpha,\sigma}}(x)\bigr|
		&\leq 2\sqrt{2}\,\sqrt n\,
		d_Q(E_{\lambda,\sigma},E_{\alpha,\sigma})\\
		&\leq 6\pi\sqrt n\,
		\bigl|A(q,\sigma\kappa)-A(q',\sigma\kappa')\bigr|\\
		&=6\pi\sqrt{2\pi}\,
		\bigl(|q-q'|^2+|\kappa-\kappa'|^2\bigr)^{\frac12}\\
		&\leq 6\pi\sqrt{2\pi}\,\rho_*(\alpha,\lambda).
	\end{aligned}
	\]
	Hence
	\begin{equation}\label{1+sqrtn}
		1+\sqrt n\,s_{E_{\lambda,\sigma}}(x)
		\leq 6\pi\sqrt{2\pi}\,
		\bigl(1+\rho_*(\alpha,\lambda)\bigr)
		\bigl(1+\sqrt n\,s_{E_{\alpha,\sigma}}(x)\bigr).
	\end{equation}
	By the definition of $w_\lambda$,
	\eqref{eq:wq-local-field}, and \eqref{1+sqrtn},
	\[
	\begin{aligned}
		|w_\lambda(x)|
		&\leq \sum_{\alpha\in\Lambda_B(\Omega)}
		\bigl|\bigl(G_\Omega^{-\frac12}\bigr)_{\alpha\lambda}\bigr|
		\,|v_\alpha(x)|\\
		&\leq C_{d,M}n^{\frac{d-1}{4}}
		\sum_{\alpha\in\Lambda_B(\Omega)}
		\sum_{\sigma\in\{-1,1\}^{d-1}}
		\frac{\bigl|\bigl(G_\Omega^{-\frac12}\bigr)_{\alpha\lambda}\bigr|}
		{\bigl(1+\sqrt n\,s_{E_{\alpha,\sigma}}(x)\bigr)^M}\\
		&\leq C_{d,M}n^{\frac{d-1}{4}}
		\left(
		\sum_{\alpha\in\Lambda_B(\Omega)}
		\bigl(1+\rho_*(\alpha,\lambda)\bigr)^M
		\bigl|\bigl(G_\Omega^{-\frac12}\bigr)_{\alpha\lambda}\bigr|
		\right)\\
		&\quad \times \left(\sum_{\sigma\in\{-1,1\}^{d-1}}
		\frac{1}{\bigl(1+\sqrt n\,s_{E_{\lambda,\sigma}}(x)\bigr)^M}\right).
	\end{aligned}
	\]
	The sum over $\alpha$ is at most $\|G_\Omega^{-\frac12}\|_M$.
	Thus \eqref{eq:wq-orthonormal-coeff} gives
	\[
	|w_\lambda(x)|
	\leq C_{d,M}n^{\frac{d-1}{4}}
	\sum_{\sigma\in\{-1,1\}^{d-1}}
	\frac{1}{\bigl(1+\sqrt n\,s_{E_{\lambda,\sigma}}(x)\bigr)^M}.
	\]
\end{proof}
For arbitrary graph cubes $\Omega$ and $\widetilde\Omega$, let
$w_\lambda$ and $\widetilde w_\mu$ be the corresponding vectors
constructed in the proof of Lemma~\ref{lem:local-packet-orthonormalization}. 
Let $U$ and $V$ be
as in \eqref{def_UV}.

\begin{lemma}\label{lem:orthonormal-packet-cross}
	Suppose that $e_n\to0$ and $B\to\infty$ as $n\to\infty$.
	For every fixed $M>0$ and all sufficiently
	large $n$, there exists a constant $C_{d,M}''>0$, independent of
	$n$, $\Omega$, $\widetilde\Omega$, $\lambda$, $\mu$, $U$, and $V$, such that
	\begin{equation}\label{eq:wq-orthonormal-cross}
		|\langle Uw_\lambda,V\widetilde w_\mu\rangle|
		\leq C_{d,M}''
		\sum_{\sigma,\tau\in\{-1,1\}^{d-1}}
		\bigl(1+\sqrt n\,d_Q(UE_{\lambda,\sigma},
		V\widetilde E_{\mu,\tau})\bigr)^{-M}
	\end{equation}
	for all $\lambda\in\Lambda_B(\Omega)$ and
	$\mu\in\Lambda_B(\widetilde\Omega)$.
\end{lemma}

\begin{proof}
	For $\alpha\in\Lambda_B(\Omega)$ and
	$\gamma\in\Lambda_B(\widetilde\Omega)$,
	by Lemma~\ref{lem:graph-coordinate-distance}, we have
	\begin{equation}\label{1+n}
	\begin{aligned}
		&1+\sqrt n\,d_Q(UE_{\lambda,\sigma},V\widetilde E_{\mu,\tau})\\
		\leq& 1+\sqrt n\,d_Q(UE_{\alpha,\sigma},V\widetilde E_{\gamma,\tau})
		+\sqrt n\,d_Q(E_{\lambda,\sigma},E_{\alpha,\sigma})
		+\sqrt n\,d_Q(\widetilde E_{\gamma,\tau},\widetilde E_{\mu,\tau})\\
		\leq& 1+\sqrt n\,d_Q(UE_{\alpha,\sigma},V\widetilde E_{\gamma,\tau})
		+3\pi\sqrt{\pi}\,
		\bigl(\rho_*(\alpha,\lambda)+\widetilde\rho_*(\gamma,\mu)\bigr)\\
		\leq&
		3\pi\sqrt{\pi}\,
		(1+\rho_*(\alpha,\lambda))
		\bigl(1+\sqrt n\,d_Q(UE_{\alpha,\sigma},
		V\widetilde E_{\gamma,\tau})\bigr)
		(1+\widetilde\rho_*(\gamma,\mu)).
	\end{aligned}
	\end{equation}
	Then by \eqref{def_wlambda}, Lemma~\ref{lem:cross-packet-estimate} and \eqref{1+n},
	\begin{align*}
		&|\langle Uw_\lambda,V\widetilde w_\mu\rangle|\\
		\leq& \sum_{\alpha,\gamma}
		\bigl|\bigl(G_\Omega^{-\frac12}\bigr)_{\alpha\lambda}\bigr|
		\bigl|\bigl(G_{\widetilde\Omega}^{-\frac12}\bigr)_{\gamma\mu}\bigr|
		|\langle Uv_\alpha,V\widetilde v_\gamma\rangle|\\
		\leq& C_{d,M}
		\sum_{\alpha,\gamma}\sum_{\sigma,\tau}
		\frac{
			\bigl|\bigl(G_\Omega^{-\frac12}\bigr)_{\alpha\lambda}\bigr|
			\bigl|\bigl(G_{\widetilde\Omega}^{-\frac12}\bigr)_{\gamma\mu}\bigr|
		}{
			\bigl(1+\sqrt n\,d_Q(UE_{\alpha,\sigma},
			V\widetilde E_{\gamma,\tau})\bigr)^M
		}\\
		\leq& C_{d,M}(3\pi\sqrt{\pi})^M\\
		&\quad\times
		\sum_{\sigma,\tau}
		\frac{
			\displaystyle
			\left(\sum_\alpha
			(1+\rho_*(\alpha,\lambda))^M
			\bigl|\bigl(G_\Omega^{-\frac12}\bigr)_{\alpha\lambda}\bigr|
			\right)
			\left(\sum_\gamma
			(1+\widetilde\rho_*(\gamma,\mu))^M
			\bigl|\bigl(G_{\widetilde\Omega}^{-\frac12}\bigr)_{\gamma\mu}\bigr|
			\right)
		}{
			\bigl(1+\sqrt n\,d_Q(UE_{\lambda,\sigma},
			V\widetilde E_{\mu,\tau})\bigr)^M
		}\\
		\leq& C_{d,M}(3\pi\sqrt{\pi})^M
		\|G_\Omega^{-\frac12}\|_M
		\|G_{\widetilde\Omega}^{-\frac12}\|_M
		\sum_{\sigma,\tau}
		\bigl(1+\sqrt n\,d_Q(UE_{\lambda,\sigma},
		V\widetilde E_{\mu,\tau})\bigr)^{-M}.
	\end{align*}
	Therefore, for all sufficiently large \(n\), by \eqref{eq:wq-orthonormal-coeff}, we have
	\[
	|\langle Uw_\lambda,V\widetilde w_\mu\rangle| \leq 4C_{d,M}(3\pi\sqrt{\pi})^M
	\sum_{\sigma,\tau}
	\bigl(1+\sqrt n\,d_Q(UE_{\lambda,\sigma},
	V\widetilde E_{\mu,\tau})\bigr)^{-M}.
	\]
	Taking $C_{d,M}''=4C_{d,M}(3\pi\sqrt{\pi})^M
	$, we obtain
	\eqref{eq:wq-orthonormal-cross}.
\end{proof}

\section{Construction of a global orthonormal system}\label{sec:global-system}

Set \(R=n^{-\frac{1}{4}}\), \(\ell=n^{-\frac{2}{5}}\), \(h=n^{-\frac{1}{2}}\) and \(B=n^{\frac{1}{100}}\).
We identify $M_{d-1,2}(\mathbb R)$ with $\mathbb R^{2(d-1)}$.
Write $B_Q(E,r)=\{F\in Q:d_Q(E,F)<r\}$ and let $\partial_Q$
denote the boundary relative to $Q$. For vectors and matrices,
$\|\cdot\|_\infty$ denotes the maximum absolute value of the entries.
Let $E_0$ denote $E_A$ at $A=0$.
Fix a sufficiently small $\delta>0$. For each $F\in Q$, choose a
rotation $g_F$ of $\mathbb R^{d+1}$ such that $g_FE_0=F$ and set
\[
U_F=\{g_FE_A:A\in(-\delta,\delta)^{2(d-1)}\}.
\]
Then $\{U_F\}_{F \in Q}$ forms an open cover of $Q$.
By the compactness of \(Q\), choose $F_1,\ldots,F_k\in Q$ such that
\[
Q=\bigcup_{i=1}^k U_{F_i}.
\]
Set  \(V_i=U_{F_i}\setminus\left(\bigcup_{j<i}\overline{U_{F_j}}\right)\) and \(g_i=g_{F_i}\).
Write $\Phi_i(X)=g_iE_X$.
For each $i$ and $a\in R\mathbb Z^{2(d-1)}$, define
\[
C_{i,a}=\bigl\{g_iE_A \mid A\in a+[0,R]^{2(d-1)}\bigr\}\subset Q,
\]
and 
\[
\mathcal C=
\bigl\{C_{i,a}\mid C_{i,a}\subset V_i,\
1\leq i\leq k,\ 
a\in R\mathbb Z^{2(d-1)}
\bigr\}.
\]
\begin{lemma}\label{lem:coarse-cube-disjointness}
	For any two distinct coarse cubes $C,C'\in\mathcal C$, we have
	\[
	\operatorname{int}_Q C\cap\operatorname{int}_Q C'=\varnothing,
	\]
	where $\operatorname{int}_Q$ denotes the interior relative to $Q$.
\end{lemma}

\begin{proof}
	Suppose that $C_{i,a},C_{j,b}\in\mathcal C$ and
	\begin{align}\label{Fint}
		F\in\operatorname{int}_Q C_{i,a}
		\cap\operatorname{int}_Q C_{j,b}.
	\end{align} 
	Since $C_{i,a}\subset V_i$ and $C_{j,b}\subset V_j$, we have
	$F\in V_i\cap V_j$. The sets $V_i$ are pairwise disjoint, so $i=j$.
	Set  $A=\Phi_i^{-1}(F)$.
	Then by \eqref{Fint}, 
	\[
	A\in\bigl(a+(0,R)^{2(d-1)}\bigr)
	\cap\bigl(b+(0,R)^{2(d-1)}\bigr).
	\]
	Hence for each $r=1,\ldots,2(d-1)$, we have
	\[
	a_r-b_r=(A_r-b_r)-(A_r-a_r)\in(-R,R) \cap R\mathbb Z=\{0\}.
	\]
	Thus $a=b$. 
	Together with $i=j$, this gives $C_{i,a}=C_{j,b}$.
\end{proof}

For each $C=C_{i,a}\in\mathcal C$, choose a rotation $g_C$ such that
\[
g_CE_0=g_iE_{a+\frac{R}{2}\mathbf 1},
\]
where \(\mathbf 1=\underbrace{(1,\ldots,1)}_{2(d-1)}\). 
Write $\Phi_C(X)=g_CE_X$.
Define \(C^*=\bigl\{A\in\mathbb R^{2(d-1)}\mid g_CE_A\in C\bigr\}\).

\begin{lemma}\label{lem:recentered-graph-chart}
	There exists $R_0>0$ such that, for every $0<R<R_0$ and every
	$C\in\mathcal C$, the map $A\mapsto g_CE_A$ is a bijection from
	$C^*$ onto $C$.
\end{lemma}
\begin{proof}
	Fix $C=C_{i,a}\in\mathcal C$ and set
	$a_C=a+\frac{R}{2}\mathbf 1$. For any $F\in C$, there exists
	$X\in a+[0,R]^{2(d-1)}$ such that $F=g_iE_X$. 
	Set $H=g_C^{-1}F$. 
	By Lemma~\ref{lem:graph-coordinate-distance}, we have
	\begin{equation}\label{dR}
	\begin{aligned}
		d_Q(H,E_0)
		&=d_Q(F,g_CE_0)
		=d_Q(E_X,E_{a_C})\\
		&\leq \frac{3\sqrt{2}\pi}{2}|X-a_C|\\
		&\leq \frac{3\sqrt{2}\pi}{2}\cdot\frac{\sqrt{2(d-1)}}{2}R\\
		&=\frac{3\pi\sqrt{d-1}}{2}R.
	\end{aligned}
	\end{equation}
	Let $(e_1,e_2)$ and $(f_1,f_2)$ be positively oriented
	orthonormal bases of $E_0$ and $H$, respectively.
	In these bases,
	the matrix of $P_{E_0}|_H\colon H\to E_0$ is
	$[\langle e_i,f_j\rangle]_{i,j=1}^2$. 
	Let \(R_0=\frac{\sqrt{2}}{3\pi\sqrt{d-1}}\).
	For every $0<R<R_0$, by Lemma~\ref{lem:grassmannian-comparisons} and \eqref{dR}, we have
	\[
	\begin{aligned}
		\det[\langle e_i,f_j\rangle]_{i,j=1}^2
		&=\langle\omega_{E_0},\omega_H\rangle\\
		&=1-\frac{1}{2}|\omega_H-\omega_{E_0}|^2\\
		&\geq 1-2d_Q(H, E_0)^2\\
		&\geq 1-\frac{9\pi^2(d-1)}{2}R^2 >0.
	\end{aligned}
	\]
	Thus $P_{E_0}|_H$ is an orientation-preserving isomorphism.
	Define
	\[
	A=(P_{E_0^\perp}|_H)\circ(P_{E_0}|_H)^{-1}
	\colon E_0\to E_0^\perp.
	\]
	For every $u\in E_0$, it is easy to see that \((P_{E_0}|_H)^{-1}u=u+Au\).
	Then
	\[
	H=(P_{E_0}|_H)^{-1}(E_0)=\{u+Au:u\in E_0\}=E_A.
	\]
	This proves surjectivity.
	To prove injectivity, suppose that $g_CE_A=g_CE_B$ for
	$A,B\in C^*$. 
	For each $u\in E_0$, since $u+Au\in E_A=E_B$, there exists
	$v\in E_0$ such that \(u+Au=v+Bv\).
	Then
	\[
	u=P_{E_0}(u+Au)=P_{E_0}(v+Bv)=v.
	\]
	Hence $Au=Bu$ for every $u\in E_0$, which implies that \(A=B\).
\end{proof}

	Let $R_0$ be as in Lemma~\ref{lem:recentered-graph-chart} and assume
	that $0<R<R_0$. For each $C\in\mathcal C$, set
	\[
	\mathcal A_C
	=
	\left\{
	A_0\in(2\ell)\mathbb Z^{2(d-1)}
	\;\middle|\;
	A_0+[-\ell,\ell]^{2(d-1)}
	\subset
	\{A\in C^*:\operatorname{dist}(A,\partial C^*)\geq\ell\}
	\right\}.
	\]
	For each $A_0\in\mathcal A_C$, define
	\[
	\Omega_{A_0}=A_0+[-\ell,\ell]^{2(d-1)}.
	\]
	Set
	\begin{align}\label{def_T}
	\mathcal T
	=
	\bigcup_{C\in\mathcal C}
	\Bigl\{
	\{g_CE_A\mid A\in\Omega_{A_0}\}
	\,\Bigm|\,
	A_0\in\mathcal A_C
	\Bigr\}.
	\end{align}
\begin{lemma}\label{lem:graph-cube-disjointness}
	For any two distinct $T,T'\in\mathcal T$, we have
	\[
	\operatorname{int}_Q T\cap\operatorname{int}_Q T'=\varnothing.
	\]
\end{lemma}

\begin{proof}
	Suppose that $T,T'\in\mathcal T$ and 
	\[
	F\in\operatorname{int}_Q T\cap\operatorname{int}_Q T'.
	\]
	By \eqref{def_T}, there exist $C,C'\in\mathcal C$,
	$A_0\in\mathcal A_C$ and $B_0\in\mathcal A_{C'}$ such that \(T=\{g_CE_A\mid A\in\Omega_{A_0}\}\) and \(T'=\{g_{C'}E_A\mid A\in\Omega_{B_0}\}\).
	Then
	\[
	T=\Phi_C(\Omega_{A_0})
	\subset \Phi_C(\operatorname{int}(C^*))
	\subset \operatorname{int}_Q C.
	\]
	Similarly,
	$T'\subset\operatorname{int}_Q C'$. 
	Hence \(F\in\operatorname{int}_Q C\cap\operatorname{int}_Q C'\).
	By Lemma~\ref{lem:coarse-cube-disjointness}, we have $C=C'$.
	Let $A=\Phi_C^{-1}(F)$. Then
	\[
	A\in
	\bigl(A_0+(-\ell,\ell)^{2(d-1)}\bigr)
	\cap
	\bigl(B_0+(-\ell,\ell)^{2(d-1)}\bigr).
	\]
	Thus, for each $r=1,\ldots,2(d-1)$,
	\[
	(A_0)_r-(B_0)_r
	=\bigl(A_r-(B_0)_r\bigr)-\bigl(A_r-(A_0)_r\bigr)
	\in(2\ell)\mathbb Z\cap(-2\ell,2\ell)=\{0\}.
	\]
	Therefore $A_0=B_0$, and hence $T=T'$.
\end{proof}

\begin{lemma}\label{lem:wq-tiling}
	There are constants $C_1,C_2>0$, independent of $n$, such that
	for all sufficiently large $n$,
	\[
	C^*\subset\{A\in\mathbb R^{2(d-1)}:|A|\leq C_1R\},\ C\in\mathcal C,
	\]
	and
	\begin{equation}\label{eq:wq-tiling-mass}
		0\leq 1-\sum_{T\in\mathcal T}\nu(T)
		\leq C_2\left(R+\frac{\ell}{R}\right).
	\end{equation}
\end{lemma}
\begin{proof}
	There exists $\varepsilon>0$ such that \(B_Q(E_0,\varepsilon)\subset\{E_X:|X|<1\}\).
	For any $C\in\mathcal C$ and $A\in C^*$, by the definition of $C^*$
	and \eqref{dR}, we have
	\[
	d_Q(E_A,E_0)
	=d_Q(g_CE_A,g_CE_0)
	\leq \frac{3\pi\sqrt{d-1}}{2}R<\varepsilon
	\]
	for all sufficiently large $n$. By injectivity of the graph chart,
	we have $|A|<1$. Hence Lemma~\ref{lem:graph-coordinate-distance} gives
	\[
	|A|\leq 12d_Q(E_A,E_0)
	\leq 18\pi\sqrt{d-1}\,R.
	\]
	Thus the first assertion holds with $C_1=18\pi\sqrt{d-1}$.
	
	Set $\Sigma=\bigcup\limits_{i=1}^k\partial_Q U_{F_i}$. For
	$F\in Q\setminus\Sigma$, let $i$ be the least index such that
	$F\in U_{F_i}$. Then $F\notin U_{F_j}\cup\partial_Q U_{F_j}
	=\overline{U_{F_j}}$ for every $j<i$, so $F\in V_i$. Moreover,
	since $U_{F_j}$ is open, we have
	$\partial_Q\overline{U_{F_j}}\subset\partial_Q U_{F_j}$. Hence
	\begin{align}\label{eq:wq-chart-cover}
		Q\setminus\Sigma
		\subset\bigcup_{i=1}^k V_i,
	\end{align}
	and 
	\begin{align}\label{eq:wq-chart-boundaries}
		\partial_Q V_i
		=\partial_Q\left(
		U_{F_i}\setminus\bigcup_{j<i}\overline{U_{F_j}}\right)
		\subset\partial_Q U_{F_i}
		\cup\bigcup_{j<i}\partial_Q\overline{U_{F_j}}
		\subset\bigcup_{j=1}^i\partial_Q U_{F_j}
		\subset\Sigma.
	\end{align}
	Let $F\in V_i\setminus\bigcup\limits_{C\in\mathcal C}C$ and choose
	$a\in R\mathbb Z^{2(d-1)}$ such that $F\in C_{i,a}$.
	If $C_{i,a}\subset V_i$, then $C_{i,a}\in\mathcal C$,
	contrary to the choice of $F$. Thus $C_{i,a}$ meets both
	$V_i$ and $Q\setminus V_i$. 
	Since $V_i$ is open, we have
	\[
	Q\setminus\partial_Q V_i
	=V_i\cup\bigl(Q\setminus\overline{V_i}\bigr).
	\]
	By connectedness of $C_{i,a}$, it follows that
	$C_{i,a}\cap\partial_Q V_i\neq\varnothing$.
	Choose $G\in C_{i,a}\cap\partial_Q V_i$.
	By
	\eqref{eq:wq-chart-boundaries}, we have $G\in\Sigma$.
	 Since $F\in C_{i,a}\cap U_{F_i}$, the coordinate cube
	 $a+[0,R]^{2(d-1)}$ meets $(-\delta,\delta)^{2(d-1)}$.
	 Then for sufficiently small $R$, we have
	 \[
	 a+[0,R]^{2(d-1)}
	 \subset(-\delta-R,\delta+R)^{2(d-1)}
	 \subset\{A:|A|<1\}.
	 \]
	Thus
	Lemma~\ref{lem:graph-coordinate-distance} gives
	\begin{equation}\label{3piR}
		d_Q(F,\Sigma)
		\leq d_Q(F,G)
		\leq \frac{3\sqrt{2}\pi}{2}\sqrt{2(d-1)}\,R
		=3\pi\sqrt{d-1}\,R.
	\end{equation}
	Since $d_Q(F,\Sigma)=0$ for $F\in\Sigma$, it follows from
	\eqref{eq:wq-chart-cover} and \eqref{3piR} that
	\begin{equation}\label{eq:wq-coarse-boundary-neighborhood}
		Q\setminus\bigcup_{C\in\mathcal C}C
		\subset
		\bigl\{F\in Q\mid d_Q(F,\Sigma)\leq 3\pi\sqrt{d-1}\,R\bigr\}.
	\end{equation}
	For each $i$, the boundary $\partial_Q U_{F_i}$ is a compact
	subset of the enlarged chart. Thus, for sufficiently small $R$,
	\begin{equation}\label{p}
		\bigl\{F\in Q\mid d_Q(F,\partial_Q U_{F_i})
		\leq 3\pi\sqrt{d-1}\,R\bigr\}
		\subset \Phi_i\bigl((-2\delta,2\delta)^{2(d-1)}\bigr).
	\end{equation}
	Let $F\in Q$ satisfy $d_Q(F,\Sigma)\leq 3\pi\sqrt{d-1}\,R$.
	By compactness of $\Sigma$, choose
	$\Phi_i(X)\in\partial_Q U_{F_i}$ such that
	$$d_Q(F,\Phi_i(X))=d_Q(F,\Sigma).$$ 
	Then
	$d_Q(F,\partial_Q U_{F_i})\leq 3\pi\sqrt{d-1}\,R$, so by
	\eqref{p} we write $F=\Phi_i(A)$ with
	$A\in(-2\delta,2\delta)^{2(d-1)}$. Since $|X_r|=\delta$
	for some $r$, Lemma~\ref{lem:graph-coordinate-distance} gives
	\begin{equation}\label{subset}
		\begin{aligned}
			\bigl||A_r|-\delta\bigr|
			\leq |A-X|
			\leq 12d_Q(F,\Phi_i(X))
			=12d_Q(F,\Sigma)
			\leq 36\pi\sqrt{d-1}\,R.
		\end{aligned}
	\end{equation}
	For $1\leq r\leq 2(d-1)$, set
	\[
	S_r=\bigl\{A\in(-2\delta,2\delta)^{2(d-1)}
	\mid \bigl||A_r|-\delta\bigr|\leq 36\pi\sqrt{d-1}\,R\bigr\}.
	\]
	By \eqref{subset}, $F\in\Phi_i(S_r)$ for some $i$ and $r$.
	Therefore,
	\begin{equation}\label{capcap}
		\bigl\{F\in Q\mid d_Q(F,\Sigma)\leq 3\pi\sqrt{d-1}\,R\bigr\}
		\subset \bigcup_{i=1}^k\bigcup_{r=1}^{2(d-1)}\Phi_i(S_r).
	\end{equation}
	Since $\operatorname{vol}(S_r)=O(R)$ and $\nu$ has bounded
	coordinate density, \eqref{eq:wq-coarse-boundary-neighborhood}
	and \eqref{capcap} give
	\begin{equation}\label{eq:wq-coarse-cover}
		\begin{aligned}
			\nu\left(Q\setminus\bigcup_{C\in\mathcal C}C\right)
			\leq \sum_{i=1}^k\sum_{r=1}^{2(d-1)}
			\nu\bigl(\Phi_i(S_r)\bigr)
			\leq C_3 k\sum_{r=1}^{2(d-1)}\operatorname{vol}(S_r)
			=O(R).
		\end{aligned}
	\end{equation}
	Moreover, if $C_{i,a}\in\mathcal C$, then
	$\Phi_i(a)\in C_{i,a}\subset U_{F_i}$, so
	$a\in(-\delta,\delta)^{2(d-1)}$. Consequently,
	\begin{equation}\label{eq:wq-coarse-count}
		\begin{aligned}
			\#\mathcal C
			\leq k\,\#\bigl(R\mathbb Z^{2(d-1)}
			\cap(-\delta,\delta)^{2(d-1)}\bigr)
			\leq k\left(2\left\lfloor\frac{\delta}{R}\right\rfloor+1
			\right)^{2(d-1)}
			=O\bigl(R^{-2(d-1)}\bigr).
		\end{aligned}
	\end{equation}
	Fix $C=C_{i,a}\in\mathcal C$ and let $A\in C^*$ satisfy
	$\operatorname{dist}(A,\partial C^*)>
	(1+2\sqrt{2(d-1)})\ell$. Choose a cube
	$\Omega=A_0+[-\ell,\ell]^{2(d-1)}$ of the $2\ell$-grid
	containing $A$. Since
	\[
	\sup_{X\in\Omega}|X-A|
	\leq\operatorname{diam}(\Omega)
	=2\sqrt{2(d-1)}\,\ell
	<\operatorname{dist}(A,\partial C^*),
	\]
	we have $\Omega\subset C^*$ and
	\[
	\operatorname{dist}(\Omega,\partial C^*)
	\geq\operatorname{dist}(A,\partial C^*)
	-\operatorname{diam}(\Omega)>\ell.
	\]
	Thus $A_0\in\mathcal A_C$, and hence
	\begin{equation}\label{eq:wq-fine-boundary-neighborhood}
		C^*\setminus\bigcup_{A_0\in\mathcal A_C}\Omega_{A_0}
		\subset
		\bigl\{A\in C^*\mid \operatorname{dist}(A,\partial C^*)
		\leq (1+2\sqrt{2(d-1)})\ell\bigr\}.
	\end{equation}
	Set $D=a+[0,R]^{2(d-1)}$ and $\Psi_C=\Phi_i^{-1}\circ\Phi_C$.
	By Lemma~\ref{lem:recentered-graph-chart}, $\Psi_C(C^*)=D$.
	Since both charts are homeomorphisms onto open subsets of $Q$,
	we also have $\Psi_C(\partial C^*)=\partial D$. Applying
	Lemma~\ref{lem:graph-coordinate-distance} in both charts gives
	\[
	|\Psi_C(A)-\Psi_C(A')|
	\leq 12d_Q\bigl(\Phi_C(A),\Phi_C(A')\bigr)
	\leq 18\sqrt{2}\pi\,|A-A'|,
	\ A,A'\in C^*.
	\]
	Taking the infimum over $A'\in\partial C^*$ gives
	\begin{align}\label{dist}
	\operatorname{dist}(\Psi_C(A),\partial D)
	\leq 18\sqrt{2}\pi\,\operatorname{dist}(A,\partial C^*).
	\end{align}
	Set \(W_C=\Psi_C\left(
	C^*\setminus\bigcup\limits_{A_0\in\mathcal A_C}\Omega_{A_0}
	\right)\).
	For any $X=\Psi_C(A)\in W_C$, by \eqref{dist} and
	\eqref{eq:wq-fine-boundary-neighborhood}, we have
	\[
	\operatorname{dist}(X,\partial D)
	\leq 18\sqrt{2}\pi\,\operatorname{dist}(A,\partial C^*)
	\leq c_d\ell,
\]
where \(c_d=18\sqrt{2}\pi\bigl(1+2\sqrt{2(d-1)}\bigr)\).
Therefore,
\begin{align}\label{WCsubset}
	W_C \subset \{X\in D\mid\operatorname{dist}(X,\partial D)\leq c_d\ell\}.
\end{align}
Since $\frac{\ell}{R}=n^{-3/20}\to 0$, we have $2c_d\ell<R$ for
all sufficiently large $n$. For $D=a+[0,R]^{2(d-1)}$, we have
\[
\{X\in D\mid \operatorname{dist}(X,\partial D)>c_d\ell\}
=a+(c_d\ell,R-c_d\ell)^{2(d-1)}.
\]
	Since $X\in D=a+[0,R]^{2(d-1)}$, we have
	\[
	\operatorname{dist}(X,\partial D)
	=
	\min_{1\leq r\leq 2(d-1)}
	\min\{X_r-a_r,\ a_r+R-X_r\}.
	\] 
	Hence, by \eqref{WCsubset},
	\begin{equation}\label{eq:wq-wc-volume}
		\begin{aligned}
			\operatorname{vol}(W_C)
			&\leq \operatorname{vol}\bigl(
			\{X\in D\mid \operatorname{dist}(X,\partial D)\leq c_d\ell\}
			\bigr)\\
			&=R^{2(d-1)}-(R-2c_d\ell)^{2(d-1)}\\
			&\leq 4(d-1)c_d\ell R^{2d-3}.
		\end{aligned}
	\end{equation}
	By the definition of $W_C$ and the identity
	$\Phi_i\circ\Psi_C=\Phi_C$, we have
	\[
	\Phi_i(W_C)
	=C\setminus\bigcup_{A_0\in\mathcal A_C}\Phi_C(\Omega_{A_0}).
	\]
	Let $\rho_i$ be the density of $\nu$ in the chart $\Phi_i$.
	Then \eqref{eq:wq-wc-volume} gives
	\begin{equation}\label{eq:wq-fine-loss}
		\begin{aligned}
			\nu\left(
			C\setminus\bigcup_{A_0\in\mathcal A_C}\Phi_C(\Omega_{A_0})
			\right)
			&=\nu\bigl(\Phi_i(W_C)\bigr)\\
			&=\int_{W_C}\rho_i(X)\,dX\\
			&\leq \|\rho_i\|_\infty\,\operatorname{vol}(W_C)\\
			&\leq 4(d-1)c_d\|\rho_i\|_\infty\,\ell R^{2d-3}\\
			&=O(\ell R^{2d-3}).
		\end{aligned}
	\end{equation}
	By Lemma~\ref{lem:graph-cube-disjointness}, the sets in
	$\mathcal T$ have pairwise disjoint interiors. 
	Let $T=\Phi_C(\Omega_{A_0})\in\mathcal T$ and  $\rho_C$ be the density of $\nu$ in the chart $\Phi_C$.
	Then
	\[
	\nu(\partial_Q T)
	=\nu\bigl(\Phi_C(\partial\Omega_{A_0})\bigr)
	=\int_{\partial\Omega_{A_0}}\rho_C(A)\,dA
	=0.
	\]
	Therefore, 
	\begin{align}\label{nu1}
	\nu\left(\bigcup_{T\in\mathcal T}T\right)
	=\nu\left(\bigcup_{T\in\mathcal T}\operatorname{int}_Q T\right)
	=\sum_{T\in\mathcal T}\nu(T).
	\end{align}
	Consequently, by
	$\nu(Q)=1$, \eqref{nu1},\eqref{eq:wq-coarse-cover},
	\eqref{eq:wq-coarse-count} and \eqref{eq:wq-fine-loss}, we have
	\[
	\begin{aligned}
		0\leq 1-\sum_{T\in\mathcal T}\nu(T)
		&=\nu\left(Q\setminus\bigcup_{T\in\mathcal T}T\right)\\
		&\leq \nu\left(Q\setminus\bigcup_{C\in\mathcal C}C\right)
		+\sum_{C\in\mathcal C}
		\nu\left(
		C\setminus\bigcup_{A_0\in\mathcal A_C}\Phi_C(\Omega_{A_0})
		\right)\\
		&\leq C_2\left(R+R^{-2(d-1)}\ell R^{2d-3}\right)
		=C_2\left(R+\frac{\ell}{R}\right).
	\end{aligned}
	\]
	This proves \eqref{eq:wq-tiling-mass}.
\end{proof}

For $T=\Phi_C(\Omega_{A_0})\in\mathcal T$, set $c_T=\Phi_C(A_0)$.
Let $\lambda=(q,\kappa)\in\Lambda_B(\Omega_{A_0})$.
For $\sigma\in\{-1,1\}^{d-1}$, we have
\[
g_CE_{\lambda,\sigma}
=\Phi_C\bigl(A_0+\sqrt{2\pi}\,h[q,-\sigma\kappa]\bigr).
\]

\begin{proposition}\label{prop:packet-assembly}
	There is a constant $C>0$, independent of $n$, such that
	for all sufficiently large $n$ there exists an orthonormal system
	$\{u_j\}_{j=1}^{J_n}$ in
	$\mathcal H_n(S^d)\otimes_{\mathbb R}\mathbb C$ with the
	following properties.
	
	\begin{enumerate}
		\item[\textup{(i)}]
		\(\left(1-Cn^{-\frac{9}{100}}\right)D_n\leq J_n\leq D_n\).
		\item[\textup{(ii)}]
		For each $j=1,\ldots,J_n$, there are $C_j\in\mathcal C$,
		$A_{0,j}\in\mathcal A_{C_j}$ and
		$\lambda_j\in\Lambda_B(\Omega_{A_{0,j}})$ such that, with
		$T_j=\Phi_{C_j}(\Omega_{A_{0,j}})$,
		\[
		\bigl\{g_{C_j}E_{\lambda_j,\sigma}
		\mid \sigma\in\{-1,1\}^{d-1}\bigr\}
		\subset B_Q(c_{T_j},C\ell).
		\]
		
		\item[\textup{(iii)}]
		For all distinct $T,T'\in\mathcal T$, \(d_Q(c_T,c_{T'})\geq \frac{\ell}{12}\).
		
		\item[\textup{(iv)}]
		For every $F\in Q$,
		\[
		\#\bigl\{
		(j,\sigma)\in\{1,\ldots,J_n\}\times\{-1,1\}^{d-1}
		\mid g_{C_j}E_{\lambda_j,\sigma}\in B_Q(F,h)
		\bigr\}<C.
		\]
	\end{enumerate}
\end{proposition}
\begin{proof}
	Let $R_*=\max\{1,C_1\}R$, where $C_1$ is the constant in
	Lemma~\ref{lem:wq-tiling}. For $C\in\mathcal C$ and
	$A_0\in\mathcal A_C$, we have $|A_0|\leq R_*$ and $\ell\leq R_*$.
	Moreover,
	\[
	n(R_*\ell^2+\ell^3)
	=\max\{1,C_1\}n^{-\frac{1}{20}}+n^{-\frac{1}{5}}
	\longrightarrow 0.
	\]
	By Lemma~\ref{lem:local-packet-orthonormalization} with $R=R_*$, for all sufficiently large $n$, we obtain an orthonormal family
	$\{w_\lambda\}_{\lambda\in\Lambda_B(\Omega_{A_0})}$ for every
	$C\in\mathcal C$ and $A_0\in\mathcal A_C$.
	Set
	\[
	J_n=\sum_{C\in\mathcal C}\sum_{A_0\in\mathcal A_C}
	\#\Lambda_B(\Omega_{A_0}).
	\]
	Enumerate the triples $(C,A_0,\lambda)$ with
	$C\in\mathcal C$, $A_0\in\mathcal A_C$, and
	$\lambda\in\Lambda_B(\Omega_{A_0})$ as
	$(C_j,A_{0,j},\lambda_j)$, $1\leq j\leq J_n$. 
	For each $j$, 
	set \(T_j=\Phi_{C_j}(\Omega_{A_{0,j}})\) and \(w_j=w_{\lambda_j}\circ g_{C_j}^{-1}\).
	Write $w_{\lambda_j}=P|_{S^d}$, where $P$ is a complex-valued
	homogeneous harmonic polynomial of degree $n$. Since $g_{C_j}$
	is orthogonal, we have
	\[
	(P\circ g_{C_j}^{-1})(tx)
	=t^n(P\circ g_{C_j}^{-1})(x),
	\]
	and
    \[
	\Delta(P\circ g_{C_j}^{-1})
	=(\Delta P)\circ g_{C_j}^{-1}=0.
	\]
	Thus $w_j=(P\circ g_{C_j}^{-1})|_{S^d}$ belongs to
	$\mathcal H_n(S^d)\otimes_{\mathbb R}\mathbb C$.
	Suppose that $T_i=T_j$. Then
	\[
	\varnothing\neq T_i\subset
	\operatorname{int}_Q C_i\cap\operatorname{int}_Q C_j,
	\]
	so $C_i=C_j$. Since $\Phi_{C_i}$ is injective,
	$\Omega_{A_{0,i}}=\Omega_{A_{0,j}}$ and hence $A_{0,i}=A_{0,j}$.
	Thus $w_{\lambda_i}$ and $w_{\lambda_j}$ belong to the same
	orthonormal family. By rotation invariance of surface measure,
	\[
	\langle w_i,w_j\rangle
	=\langle w_{\lambda_i},w_{\lambda_j}\rangle
	=\delta_{\lambda_i,\lambda_j}
	=\delta_{ij}.
	\]
	For each $j$, write \(T_j=\Phi_{C_j}(\Omega_{A_{0,j}})\) and \(\lambda_j=(q_j,\kappa_j)\in\Lambda_B(\Omega_{A_{0,j}})\).
	
	\medskip
	\noindent\textup{(ii)}
	Write $\lambda_j=(q_j,\kappa_j)$. For $\sigma\in\{-1,1\}^{d-1}$, set
	\begin{equation*}
		A_{j,\sigma}=A_{0,j}+\sqrt{2\pi}\,h[q_j,-\sigma\kappa_j].
	\end{equation*}
	Then $\Phi_{C_j}(A_{j,\sigma})=g_{C_j}E_{\lambda_j,\sigma}$.
	By the definition of
	$\Lambda_B(\Omega_{A_{0,j}})$,
	\begin{equation}\label{eq:wq-assembly-coordinate-bound}
		\|A_{j,\sigma}-A_{0,j}\|_\infty
		=\sqrt{2\pi}\,h\|\lambda_j\|_\infty
		\leq \ell-\sqrt{2\pi}\,Bh.
	\end{equation}
	Thus $A_{j,\sigma}\in\Omega_{A_{0,j}}$ and
	\begin{equation}\label{eq:wq-assembly-margin}
		\operatorname{dist}(A_{j,\sigma},\partial\Omega_{A_{0,j}})
		=\ell-\|A_{j,\sigma}-A_{0,j}\|_\infty
		\geq \sqrt{2\pi}\,Bh.
	\end{equation}
	By Lemma~\ref{lem:graph-coordinate-distance} and
	\eqref{eq:wq-assembly-coordinate-bound},
	\[
	\begin{aligned}
		d_Q(g_{C_j}E_{\lambda_j,\sigma},c_{T_j})
		&=d_Q\bigl(\Phi_{C_j}(A_{j,\sigma}),\Phi_{C_j}(A_{0,j})\bigr)\\
		&\leq \frac{3\sqrt{2}\pi}{2}|A_{j,\sigma}-A_{0,j}|\\
		&\leq 3\pi\sqrt{d-1}\,\ell.
	\end{aligned}
	\]
	This proves \textup{(ii)} for any $C>3\pi\sqrt{d-1}$.

	\medskip
	\noindent\textup{(iii)}
	There is an $\varepsilon>0$, independent of $C$, such that
	\[
	B_Q(\Phi_C(0),\varepsilon)\subset\Phi_C(\{A:|A|<1\}).
	\]
	Fix $T=\Phi_C(\Omega_{A_0})$ and
	$A\in\operatorname{int}(\Omega_{A_0})$. Let
	$0<t\leq\operatorname{dist}(A,\partial\Omega_{A_0})$.
	If $F\in B_Q(\Phi_C(A),\frac{t}{12})$, then by \eqref{dR}, we have
	\[
	d_Q(F,\Phi_C(0))
	\leq d_Q(F,\Phi_C(A))+d_Q(\Phi_C(A),\Phi_C(0))
	<\frac{\ell}{12}+\frac{3\pi\sqrt{d-1}}{2}R
	<\varepsilon
	\]
	for all sufficiently large $n$. Thus $F=\Phi_C(X)$ for some
	$|X|<1$. 
	By Lemma~\ref{lem:wq-tiling}, $|A|\leq C_1R<1$. Hence
	Lemma~\ref{lem:graph-coordinate-distance} gives
	\[
	|X-A|\leq 12d_Q(F,\Phi_C(A))<t
	\leq\operatorname{dist}(A,\partial\Omega_{A_0}).
	\]
	Consequently $X\in\operatorname{int}(\Omega_{A_0})$, and
	\begin{equation}\label{eq:wq-assembly-interior-ball}
		B_Q\left(\Phi_C(A),\frac{t}{12}\right)
		\subset\operatorname{int}_Q T.
	\end{equation}
	Taking $A=A_0$ and $t=\ell$ in
	\eqref{eq:wq-assembly-interior-ball}, we obtain
	\[
	B_Q\left(c_T,\frac{\ell}{12}\right)\subset\operatorname{int}_Q T.
	\]
	Therefore, for distinct $T,T'\in\mathcal T$, \(d_Q(c_T,c_{T'})\geq\frac{\ell}{12}\).

	\medskip
	\noindent\textup{(iv)}
	Since $Bh=o(\ell)$, it follows from
	\eqref{eq:wq-assembly-margin} and \eqref{eq:wq-assembly-interior-ball}
	that, for all sufficiently large $n$,
	\[
	B_Q\left(g_{C_i}E_{\lambda_i,\sigma},
	\frac{\sqrt{2\pi}}{12}Bh\right)
	\subset\operatorname{int}_Q T_i.
	\]
	If $T_i\neq T_j$, then \eqref{eq:wq-assembly-margin} and
	Lemma~\ref{lem:graph-cube-disjointness} give
	\[
	g_{C_j}E_{\lambda_j,\tau}
	\in\operatorname{int}_Q T_j
	\subset Q\setminus\operatorname{int}_Q T_i.
	\]
	Hence, if \(T_i\neq T_j\),
	\begin{equation}\label{eq:wq-assembly-plane-separation}
		d_Q(g_{C_i}E_{\lambda_i,\sigma},g_{C_j}E_{\lambda_j,\tau})
		\geq\frac{\sqrt{2\pi}}{12}Bh.
	\end{equation}
	Fix $F\in Q$ and suppose that
	$g_{C_i}E_{\lambda_i,\sigma}\in B_Q(F,h)$. Set
	$(C,A_0)=(C_i,A_{0,i})$ and $T=T_i=\Phi_C(\Omega_{A_0})$.
	If $g_{C_j}E_{\lambda_j,\tau}\in B_Q(F,h)$, then, for all
	sufficiently large $n$,
	\[
	d_Q\bigl(g_{C_i}E_{\lambda_i,\sigma},
	g_{C_j}E_{\lambda_j,\tau}\bigr)
	<2h<\frac{\sqrt{2\pi}}{12}Bh.
	\]
	Hence \eqref{eq:wq-assembly-plane-separation} gives $T_j=T$.
	For each such index $(j,\sigma)$,
	\begin{equation}\label{phic}
	\begin{aligned}
	\Phi_C^{-1}(g_{C_j}E_{\lambda_j,\sigma})
	&=A_{j,\sigma}
	=A_0+\sqrt{2\pi}\,h[q_j,-\sigma\kappa_j]\\
	&\in A_0+\sqrt{2\pi}\,h
	\left(\left(\frac{1}{2}\mathbb Z\right)^{d-1}
	\times\mathbb Z^{d-1}\right).
	\end{aligned}
	\end{equation}
	Fix $X=A_{i,\sigma}$ with $\Phi_C(X)\in B_Q(F,h)$.
	For every $(j,\tau)$ such that
	$g_{C_j}E_{\lambda_j,\tau}\in B_Q(F,h)$,
	by Lemma~\ref{lem:graph-coordinate-distance},
	\[
	\|A_{j,\tau}-X\|_\infty
	\leq |A_{j,\tau}-X|
	\leq 12d_Q\bigl(\Phi_C(A_{j,\tau}),\Phi_C(X)\bigr)
	<24h.
	\]
	Thus $A_{j,\tau}\in X+[-24h,24h]^{2(d-1)}$. 
	Consequently, by \eqref{phic} and the
	injectivity of $(j,\tau)\mapsto A_{j,\tau}$, we have
	\[
	\begin{aligned}
		&\#\bigl\{(j,\tau)\in\{1,\ldots,J_n\}\times\{-1,1\}^{d-1}
		\mid g_{C_j}E_{\lambda_j,\tau}\in B_Q(F,h)\bigr\}\\
		\leq&
		\left(1+\left\lfloor
		\frac{48h}{\frac{\sqrt{2\pi}}{2}h}
		\right\rfloor\right)^{2(d-1)}
		=
		\left(1+\left\lfloor\frac{96}{\sqrt{2\pi}}\right\rfloor
		\right)^{2(d-1)}.
	\end{aligned}
	\]
	This proves \textup{(iv)} with $C$ strictly larger than the last bound.

	\medskip
	\noindent\textup{(i)}
	Define $W$ by $We_j=w_j$ and set $G=W^*W$. 
	Then $G_{ij}=\delta_{ij}$ whenever $T_i=T_j$.
	Fix $F\in Q$ and $t>0$. For each $i$ with
	$U_{F_i}\cap B_Q(F,th)\neq\varnothing$, choose $X_i$ such that
	$\Phi_i(X_i)\in U_{F_i}\cap B_Q(F,th)$. By
	Lemma~\ref{lem:graph-coordinate-distance},
	\[
	\|X-X_i\|_\infty
	\leq 12d_Q\bigl(\Phi_i(X),\Phi_i(X_i)\bigr)<24th.
	\]
	Divide each coordinate interval of $X_i+[-24th,24th]^{2(d-1)}$
	 into intervals of length at most $\frac{h}{6\pi\sqrt{d-1}}$
	 and take their Cartesian products.
	If $X$ and $Y$ lie in the
	same box and $\Phi_i(X),\Phi_i(Y)\in U_{F_i}$, then
	\[
	d_Q(\Phi_i(X),\Phi_i(Y))
	\leq \frac{3\sqrt{2}\pi}{2}|X-Y|
	\leq \frac{h}{2}.
	\]
	For each box $D$ for which $\Phi_i(D)\cap U_{F_i}\cap B_Q(F,th)$
	is nonempty, choose a point $P_r$ in this intersection. Then
	\[
	B_Q(F,th)\subset\bigcup_{r=1}^N B_Q(P_r,h),
	\]
	where
	\[
	N\leq k\left(1+\left\lceil\frac{48th}{\frac{h}{6\pi\sqrt{d-1}}}\right\rceil
	\right)^{2(d-1)}
	\leq C(1+t)^{2(d-1)}.
	\]
	Hence, by \textup{(iv)},
	\begin{equation}\label{eq:wq-assembly-count}
		\begin{aligned}
			\#\bigl\{(j,\sigma)\mid
			g_{C_j}E_{\lambda_j,\sigma}\in B_Q(F,th)\bigr\}
			&\leq \sum_{r=1}^N
			\#\bigl\{(j,\sigma)\mid
			g_{C_j}E_{\lambda_j,\sigma}\in B_Q(P_r,h)\bigr\}\\
			&\leq CN
			\leq C'(1+t)^{2(d-1)}.
		\end{aligned}
	\end{equation}
	Fix $s>0$, put $M=2(d-1)+s+2$ and set
	$b=\frac{\sqrt{2\pi}}{12}$. Since $(G-I)_{ji}=0$ whenever
	$T_i=T_j$, \eqref{eq:wq-orthonormal-cross} gives
	\[
	\sum_i |(G-I)_{ji}|
	\leq C_{d,M}\sum_{\sigma,\tau}
	\sum_{i:\,T_i\neq T_j}
	\left(
	1+\frac{
		d_Q(g_{C_j}E_{\lambda_j,\sigma},g_{C_i}E_{\lambda_i,\tau})
	}{h}
	\right)^{-M}.
	\]
	Consequently, by \eqref{eq:wq-assembly-plane-separation} and \eqref{eq:wq-assembly-count},
	\[
	\begin{aligned}
		\max_j\sum_i |(G-I)_{ji}|
		&\leq C_{d,M}\sum_{k\geq\lfloor bB\rfloor}
		(k+2)^{2(d-1)}(1+k)^{-M}\\
		&\leq C_{d,s}\sum_{k\geq\lfloor bB\rfloor}(1+k)^{-s-2}
		\leq C_{d,s}'B^{-s}.
	\end{aligned}
	\]
	Since $G-I$ is Hermitian, its absolute row and column sums have
	the same maximum. Hence, for all sufficiently large $n$,
	\[
	\|G-I\|_{2\to2}\leq\|G-I\|_0
	=\max_j\sum_i |(G-I)_{ji}|
	\leq C_{d,s}B^{-s}<\frac{1}{2}.
	\]
	Thus $G$ is positive definite. 
	Put \(K=G^{-\frac{1}{2}}-I\) and \(u_j=W(I+K)e_j\), \(1\leq j\leq J_n\).
	By the binomial series,
	\begin{equation}\label{eq:wq-assembly-correction}
		\|K\|_0
		\leq (1-\|G-I\|_0)^{-\frac{1}{2}}-1
		\leq C_{d,s}'B^{-s}.
	\end{equation}
	Since $G=W^*W$,
	\[
	\bigl(W(I+K)\bigr)^*W(I+K)
	=G^{-\frac{1}{2}}GG^{-\frac{1}{2}}=I.
	\]
	Hence $\{u_j\}_{j=1}^{J_n}$ is orthonormal and $J_n\leq D_n$.
	Use the same indices $(C_j,A_{0,j},\lambda_j)$ for $u_j$ as for $w_j$.
	By local orthonormality and \eqref{eq:wq-orthonormal-field},
	$\|w_i\|_2=1$ and $\|w_i\|_\infty\leq C_d n^{\frac{d-1}{4}}$.
	Since $u_j-w_j=WKe_j=\sum\limits_i K_{ij}w_i$, it follows from
	\eqref{eq:wq-assembly-correction} that
	\begin{equation}\label{eq:wq-assembly-error}
		\begin{aligned}
			\|u_j-w_j\|_2+n^{-\frac{d-1}{4}}\|u_j-w_j\|_\infty
			&\leq (1+C_d)\sum_i |K_{ij}|\\
			&\leq (1+C_d)\|K\|_0\\
			&\leq C_{d,s}B^{-s}.
		\end{aligned}
	\end{equation}
	For $T=\Phi_C(\Omega_{A_0})\in\mathcal T$, set
	$N_T=\#\Lambda_B(\Omega_{A_0})$. 
	By rotation invariance of $\nu$ and Lemma~\ref{lem:local-packet-count},
	\[
	\left|\frac{N_T}{D_n}-\nu(T)\right|
	\leq C_d\nu(T)
	\left(R+\frac{1}{n}+\frac{(B+1)h}{\ell}\right).
	\]
	Since $J_n=\sum\limits_{T\in\mathcal T}N_T\leq D_n$ and
	$\sum\limits_{T\in\mathcal T}\nu(T)\leq 1$, it follows from
	\eqref{eq:wq-tiling-mass} that
	\[
	\begin{aligned}
		0\leq 1-\frac{J_n}{D_n}
		&=1-\sum_{T\in\mathcal T}\nu(T)
		+\sum_{T\in\mathcal T}\left(\nu(T)-\frac{N_T}{D_n}\right)\\
		&\leq C_d\left(R+\frac{\ell}{R}+h+\frac{(B+1)h}{\ell}\right)\\
		&=O_d\left(n^{-\frac{9}{100}}\right).
	\end{aligned}
	\]
	Taking $C$ large enough for \textup{(i)}, \textup{(ii)}, and
	\textup{(iv)} completes the proof.
\end{proof}
\begin{proposition}\label{prop:assembled-packet-estimates}
	Let $\{u_j\}_{j=1}^{J_n}$ and $(C_j,\lambda_j)_{j=1}^{J_n}$
	be as in Proposition~\ref{prop:packet-assembly}. For every $M>0$ and $A>0$,
	there exist constants $C_{d,M},C_{d,M,A}>0$ such that the following
	estimates hold for all sufficiently large $n$.\\
	\textup{(i)} For every $1\leq j\leq J_n$ and $x\in S^d$,
	\begin{equation}\label{eq:assembled-pointwise}
		\begin{aligned}
			|u_j(x)|
			\leq C_{d,M}n^{\frac{d-1}{4}}
			\sum_{\sigma\in\{-1,1\}^{d-1}}
			\left(1+\sqrt n\,s_{g_{C_j}E_{\lambda_j,\sigma}}(x)\right)^{-M}
			+ C_{d,M,A}n^{-A}.
		\end{aligned}
	\end{equation}
	\textup{(ii)} For every $1\leq i,j\leq J_n$ and all $U,V$ as in \eqref{def_UV},
	\begin{equation}\label{eq:assembled-cross}
		\begin{aligned}
			|\langle Uu_i,Vu_j\rangle|
			&\leq C_{d,M}
			\sum_{\sigma,\tau\in\{-1,1\}^{d-1}}
			\left(
			1+\sqrt n\,d_Q\bigl(
			Ug_{C_i}E_{\lambda_i,\sigma},
			Vg_{C_j}E_{\lambda_j,\tau}
			\bigr)
			\right)^{-M}\\
			&\quad + C_{d,M,A}n^{-A}.
		\end{aligned}
	\end{equation}
\end{proposition}
\begin{proof}
	Recall that $w_j=w_{\lambda_j}\circ g_{C_j}^{-1}$, where
	$w_{\lambda_j}$ is the function constructed on $\Omega_{A_{0,j}}$.
	Fix $M,A>0$ and take $s=100\left(A+\frac{d-1}{4}\right)$.
	By \eqref{eq:wq-assembly-error}, we have
	\begin{equation}\label{eq:assembled-norm-error}
		\|u_j-w_j\|_2+\|u_j-w_j\|_\infty
		\leq C_{d,s}n^{\frac{d-1}{4}}B^{-s}
		=C_{d,s}n^{-A},
		\ 1\leq j\leq J_n.
	\end{equation}
	
	\medskip
	\noindent\textup{(i)}
	Since $w_j(x)=w_{\lambda_j}(g_{C_j}^{-1}x)$ and \(s_{E_{\lambda_j,\sigma}}(g_{C_j}^{-1}x)
	=s_{g_{C_j}E_{\lambda_j,\sigma}}(x)\),
	by \eqref{eq:wq-orthonormal-field} and
	\eqref{eq:assembled-norm-error}, we have
	\[
	\begin{aligned}
		|u_j(x)|
		&\leq |w_j(x)|+\|u_j-w_j\|_\infty\\
		&\leq C_{d,M}n^{\frac{d-1}{4}}
		\sum_{\sigma\in\{-1,1\}^{d-1}}
		\left(1+\sqrt n\,s_{g_{C_j}E_{\lambda_j,\sigma}}(x)\right)^{-M}
		+C_{d,M,A}n^{-A}.
	\end{aligned}
	\]
	This proves \eqref{eq:assembled-pointwise}.
	
	\medskip
	\noindent\textup{(ii)}
By the Cauchy--Schwarz inequality and
	\eqref{eq:assembled-norm-error},
	\[
	\begin{aligned}
		\bigl|\langle Uu_i,Vu_j\rangle-\langle Uw_i,Vw_j\rangle\bigr|
		&\leq
		\bigl|\langle U(u_i-w_i),Vu_j\rangle\bigr|
		+\bigl|\langle Uw_i,V(u_j-w_j)\rangle\bigr|\\
		&\leq \|u_i-w_i\|_2+\|u_j-w_j\|_2
		\leq 2C_{d,s}n^{-A}.
	\end{aligned}
	\]
	Write $R_g f=f\circ g^{-1}$. Since
	$R_{g_1}R_{g_2}=R_{g_1g_2}$ and
	$R_g\overline{f}=\overline{R_gf}$, the operators
	$UR_{g_{C_i}}$ and $VR_{g_{C_j}}$ satisfy \eqref{def_UV}.
	Thus, by Lemma~\ref{lem:orthonormal-packet-cross}, we obtain 
	\[
	\begin{aligned}
		&|\langle Uu_i,Vu_j\rangle|\\
		\leq& |\langle Uw_i,Vw_j\rangle|+2C_{d,s}n^{-A}\\
		\leq& C_{d,M}
		\sum_{\sigma,\tau\in\{-1,1\}^{d-1}}
		\left(
		1+\sqrt n\,d_Q\bigl(
		Ug_{C_i}E_{\lambda_i,\sigma},
		Vg_{C_j}E_{\lambda_j,\tau}
		\bigr)
		\right)^{-M}+2C_{d,s}n^{-A}.
	\end{aligned}
	\]
	This proves \eqref{eq:assembled-cross}.
\end{proof}

\section{A bound for the projection kernel}

Recall that $h=n^{-\frac{1}{2}}$ and $\ell=n^{-\frac{2}{5}}$. 
Set $m=d-1\geq 2$
and $D=D_n$. 
Let $\mu$ be the rotation-invariant probability measure on the
group of rotations of $\mathbb R^{m+2}$.
We write $E_n=O(n^{-\infty})$ if
$E_n=O_{m,M}(n^{-M})$ for every $M>0$.
Set $a=n^{-\frac{1}{16}}$, $H=n^{-\frac{1}{4}}$ and $J=\binom{m+2}{2}$.
Identify each $E\in Q$ with its unit bivector
$(E_{ij})_{i<j}\in\mathbb R^J$. For $\rho>0$, set
\[
U_\rho=
\left\{E\in Q:
\min_{i<j}\operatorname{dist}
\bigl(E_{ij},a(\mathbb Z+1/2)\bigr)\leq\rho
\right\}.
\]
For $\mathbf k=(k_{ij})_{1\leq i<j\leq m+2}\in\mathbb Z^J$, let \(A_{\mathbf k}\) be
\[
\left\{E\in Q:
a\left(k_{ij}-\frac12\right)+H<E_{ij}<
a\left(k_{ij}+\frac12\right)-H,\ 1\leq i<j\leq m+2
\right\}.
\]
Write
\[
\mathcal A
=
\{A_{\mathbf k}:\mathbf k\in\mathbb Z^J,\
A_{\mathbf k}\neq\varnothing\}.
\]
For $A\subseteq Q$, write $-A=\{\overline E:E\in A\}$.
For nonempty subsets $A,A'\subset Q$, set
\[
d_Q(A,A')=\inf_{E\in A,\;F\in A'}d_Q(E,F).
\]

\begin{lemma}\label{lem:wc-geometric-estimates}
	For all sufficiently large $n$, the following statements hold.
	\begin{enumerate}
		\item[(i)] For every $A\in\mathcal A$, we have \(-A\in\mathcal A\), \(A\neq -A\) and \(\operatorname{diam}_{d_Q} A\leq C_m a\).
		For distinct $A,A'\in\mathcal A$, \(d_Q(A,A')\geq c_m H\).
		
		\item[(ii)] For every fixed $C_0>0$, there exists $C_{m,C_0}>0$
		such that
		\begin{equation*}
			\nu\bigl(\{E\in B_Q(E_0,C_0a):s_E(x)\leq t\}\bigr)
			\leq C_{m,C_0}a^m\min\{t^m,a^m\}
		\end{equation*}
		for all $E_0\in Q$, $x\in S^d$, and $t\geq0$.
		
		\item[(iii)] For every $L>m$, there exists $C_{m,L}>0$ such that
		\begin{equation*}
			\sup_{x\in S^d}
			h^{-m}\int_{U_{2H}}
			\left(1+\frac{s_E(x)}h\right)^{-L}\,d\nu(E)
			\leq C_{m,L}
			\left(\frac Ha+\left(\frac hH\right)^{L-m}\right).
		\end{equation*}
	\end{enumerate}
\end{lemma}
\begin{proof}
	\noindent\textup{(i)}
	For $E,F\in A_{\mathbf k}$, we have $|E_{ij}-F_{ij}|\leq a$
	for every $i<j$. Thus, by \eqref{lem:grassmannian-comparisons},
	\[
	d_Q(E,F)\leq C|E-F|\leq C\sqrt J\,a.
	\]
	Taking the supremum over $E,F\in A_{\mathbf k}$ gives
	$\operatorname{diam}_{d_Q}A_{\mathbf k}\leq C_m a$.
	By the definition of $A_{\mathbf k}$, we have
	\begin{equation}\label{eq:wc-cell-reflection}
		-A_{\mathbf k}=A_{-\mathbf k}.
	\end{equation}
	If $E\in A_{\mathbf 0}$, then $|E_{ij}|<\frac a2$ for every $i<j$.
	Hence
	\[
	1=|E|^2=\sum_{i<j}|E_{ij}|^2<\frac{Ja^2}{4}.
	\]
	Since $a\to0$, for all sufficiently large $n$ we have
	\begin{equation}\label{eq:wc-central-cell-empty}
		A_{\mathbf 0}=\varnothing.
	\end{equation}
	For distinct $\mathbf k,\mathbf k'\in\mathbb Z^J$, choose $i<j$
	with $k_{ij}\neq k'_{ij}$. By interchanging $\mathbf k$ and
	$\mathbf k'$ if necessary, assume that $k'_{ij}\geq k_{ij}+1$.
	For $E\in A_{\mathbf k}$ and $F\in A_{\mathbf k'}$, we have
	\begin{equation}\label{eq:wc-coordinate-separation}
		F_{ij}-E_{ij}
		>
		a(k'_{ij}-k_{ij}-1)+2H
		\geq 2H.
	\end{equation}
	In particular, $A_{\mathbf k}\cap A_{\mathbf k'}=\varnothing$.
	By \eqref{eq:wc-cell-reflection}, $-A\in\mathcal A$ for every
	$A\in\mathcal A$. If $A_{\mathbf k}=-A_{\mathbf k}$ is nonempty,
	then \eqref{eq:wc-cell-reflection} and the disjointness give
	$\mathbf k=-\mathbf k$, which contradicts
	\eqref{eq:wc-central-cell-empty}. Hence $A\neq -A$ for every
	$A\in\mathcal A$.
	Finally, by \eqref{lem:grassmannian-comparisons} and \eqref{eq:wc-coordinate-separation},
	\[
	d_Q(E,F)\geq c|E-F|\geq c|F_{ij}-E_{ij}|\geq c_m H.
	\]
	Taking the infimum over $E\in A_{\mathbf k}$ and
	$F\in A_{\mathbf k'}$ gives
	$d_Q(A_{\mathbf k},A_{\mathbf k'})\geq c_m H$.
	
	\medskip
	\noindent\textup{(ii)}
	Fix $x\in S^d$. Let $E_*$ denote the plane $E_0$ defined in
	Section~\ref{sec:global-system}.
	By left and right invariance of $\mu$, we have
	\begin{equation}\label{eq:rotation-measure-invariance}
		\mu(R\mathcal B)=\mu(\mathcal B R)=\mu(\mathcal B)
	\end{equation}
	for every rotation $R$ and every Borel set $\mathcal B$ of rotations.
	Let $g$ have distribution $\mu$ and set
	$E=gE_*$, $v=s_E(x)$, and $\xi=g^{-1}x$.
	For every rotation $R$ and all Borel sets $B\subset Q$ and
	$C\subset S^d$, by \eqref{eq:rotation-measure-invariance},
	\[
	\mu\{g:RE\in B\}
	=\mu\{g:(Rg)E_*\in B\}
	=\mu\{g:E\in B\}
	\]
	and
	\[
	\mu\{g:R\xi\in C\}
	=\mu\{g:(gR^{-1})^{-1}x\in C\}
	=\mu\{g:\xi\in C\}.
	\]
	Thus the distributions of $E$ and $\xi$ are rotation invariant
	and hence equal to $\nu$ and $\sigma$, respectively.
	Since $E_*=\mathbb R^2\oplus\{0\}$, we have
	\[
	1-v^2=|P_E x|^2=|P_{E_*}\xi|^2=\xi_1^2+\xi_2^2.
	\]
	The density of $\xi_1^2+\xi_2^2$ is
	$\frac{m}{2}(1-u)^{\frac{m}{2}-1}$ for $0<u<1$.
	Hence, for $0\leq t\leq1$,
	\begin{equation}\label{eq:wc-v-distribution}
		\nu\{E:s_E(x)\leq t\}
		=\mu\{g:\xi_1^2+\xi_2^2\geq 1-t^2\}
		=\frac{m}{2}\int_{1-t^2}^{1}(1-u)^{\frac{m}{2}-1}\,du
		=t^m.
	\end{equation}
	By \eqref{eq:wc-v-distribution}, $0<v<1$ for $\nu$-almost every $E$.
	For such $E$, set $e=\frac{P_E x}{|P_E x|}$ and write
	$\omega_E=e\wedge b$, where $(e,b)$ is a positively oriented
	orthonormal basis of $E$. Then
	$\langle b,x\rangle=\langle b,P_E x\rangle=0$ and
	$\langle e,x\rangle=|P_E x|=\sqrt{1-v^2}$.
	Since $e-x\perp b$ and $|b|=1$, we have
	\begin{equation}\label{eq:wc-bivector-error}
		|\omega_E-x\wedge b|^2
		=|(e-x)\wedge b|^2
		=|e-x|^2
		=2\bigl(1-\sqrt{1-v^2}\bigr)
		\leq 2v^2.
	\end{equation}
	Set $S_x=S^d\cap x^\perp$ and let $\sigma_x$ be the normalized
	surface measure on $S_x$. For every rotation $R$ with $Rx=x$,
	we have $v(RE)=v(E)$ and $b(RE)=Rb(E)$.
	For a Borel set $I\subset(0,1)$, define
	\[
	\nu_I(F)=\nu\{E:v(E)\in I,\ b(E)\in F\},
	\ F\subset S_x \text{ Borel}.
	\]
	By rotation invariance of $\nu$, $\nu_I(RF)=\nu_I(F)$ whenever
	$Rx=x$. If $\nu_I(S_x)>0$, then
	$\frac{\nu_I}{\nu_I(S_x)}$ is a rotation-invariant probability
	measure on $S_x$ and hence equals $\sigma_x$. Thus
	\[
	\nu_I(F)=\nu_I(S_x)\,\sigma_x(F)
	\]
	for every Borel set $F\subset S_x$. The same identity holds
	when $\nu_I(S_x)=0$. By \eqref{eq:wc-v-distribution},
	\[
	\nu_I(S_x)=\nu\{E:v(E)\in I\}=\int_I mv^{m-1}\,dv.
	\]
	Therefore,
	\begin{equation}\label{eq:wc-v-b-law}
		\nu\{E:v(E)\in I,\ b(E)\in F\}
		=\sigma_x(F)\int_I mv^{m-1}\,dv.
	\end{equation}
	Let $E_0\in Q$. Suppose that $0\leq t\leq a$.
	If $d_Q(E,E_0)\leq C_0a$ and $v\leq t$, then, by (2.2) and
	\eqref{eq:wc-bivector-error},
	\[
	|x\wedge b-\omega_{E_0}|
	\leq |\omega_E-x\wedge b|+|\omega_E-\omega_{E_0}|
	\leq \sqrt{2}\,t+2C_0a
	\leq C_1a,
	\]
	where $C_1=\sqrt{2}+2C_0$. Set
	\[
	F=\{b\in S_x:|x\wedge b-\omega_{E_0}|\leq C_1a\}.
	\]
	If $F\neq\varnothing$, choose $b_0\in F$. For every $b\in F$,
	\[
	|b-b_0|
	=|x\wedge(b-b_0)|
	\leq |x\wedge b-\omega_{E_0}|
	+|x\wedge b_0-\omega_{E_0}|
	\leq 2C_1a.
	\]
	Thus $\sigma_x(F)\leq C_{m,C_0}a^m$. By
	\eqref{eq:wc-v-b-law},
	\begin{align*}
		\nu\{E\in B_Q(E_0,C_0a):s_E(x)\leq t\}
		&\leq \nu\{E:v(E)\leq t,\ b(E)\in F\}\\
		&=\sigma_x(F)\int_0^t mv^{m-1}\,dv\\
		&\leq C_{m,C_0}a^mt^m.
	\end{align*}
	For $t>a$, write $d\nu=\rho(A)\,dA$ in the graph coordinates
	near $E_*$. By rotation invariance and
	Lemma~\ref{lem:graph-coordinate-distance},
	\begin{align*}
		\nu\{E\in B_Q(E_0,C_0a):s_E(x)\leq t\}
		&\leq \nu(B_Q(E_0,C_0a))\\
		&=\nu(B_Q(E_*,C_0a))\\
		&\leq \int_{|A|\leq C_m C_0a}\rho(A)\,dA\\
		&\leq C_{m,C_0}a^{2m}.
	\end{align*}
	Combining the two cases gives \textup{(ii)}.
	
	\medskip
	\noindent\textup{(iii)}
	Fix $x\in S^d$ and use the notation in \textup{(ii)}.
	Set $\delta=(2+\sqrt{2})H$. 
	For each $i<j$, write \(r_{ij}=(x_i^2+x_j^2)^{\frac{1}{2}}\) and
	\[
	D_{ij}
	=
	\left\{u\in[-r_{ij},r_{ij}]:
	\operatorname{dist}\left(u,a\left(\mathbb Z+\frac{1}{2}\right)\right)
	\leq\delta\right\}.
	\]
	Since $|b|=1$, we have
	$|(x\wedge b)_{ij}|=|x_i b_j-x_j b_i|\leq r_{ij}$.
	If $E\in U_{2H}$ and $0<v(E)\leq H$, then for some $i<j$,
	by \eqref{eq:wc-bivector-error},
	\[
	\operatorname{dist}
	\left((x\wedge b)_{ij},a\left(\mathbb Z+\frac{1}{2}\right)\right)
	\leq 2H+|\omega_E-x\wedge b|
	\leq 2H+\sqrt{2}\,v(E)
	\leq\delta.
	\]
	Thus
	\begin{equation}\label{eq:wc-grid-inclusion}
		\{E\in U_{2H}:0<v(E)\leq H\}
		\subseteq
		\bigcup_{i<j}
		\{E\in Q:0<v(E)\leq H,\ (x\wedge b(E))_{ij}\in D_{ij}\}.
	\end{equation}
	We claim that for all sufficiently large $n$ and every $i<j$,
	\begin{equation}\label{eq:wc-coordinate-grid-bound}
		\sigma_x\{b\in S_x:(x\wedge b)_{ij}\in D_{ij}\}
		\leq C_m\frac{H}{a}.
	\end{equation}
	In fact, fix $i<j$. We consider two cases.
	
	\noindent\textup{Case 1.}
	$r_{ij}<\frac{a}{4}$. For every $u\in[-r_{ij},r_{ij}]$,
	\[
	\operatorname{dist}\left(u,a\left(\mathbb Z+\frac{1}{2}\right)\right)
	\geq \frac{a}{2}-|u|
	\geq \frac{a}{2}-r_{ij}
	> \frac{a}{4}.
	\]
	Since $H=o(a)$, we have $\delta<\frac{a}{4}$ for all
	sufficiently large $n$. Hence $D_{ij}=\varnothing$, which implies that \eqref{eq:wc-coordinate-grid-bound} holds.
	
	\noindent\textup{Case 2.}
	$r_{ij}\geq\frac{a}{4}$. Write $(x\wedge b)_{ij}=r_{ij}T(b)$,
	where
	\[
	T(b)=\frac{x_i b_j-x_j b_i}{r_{ij}}.
	\]
	The coefficient vector of $T$ belongs to $x^\perp$ and has norm one.
	By rotation invariance of $\sigma_x$, the density of $T$ under
	$\sigma_x$ is
	\begin{equation}\label{eq:wc-coordinate-density}
		f_T(u)=c_m(1-u^2)^{\frac{m-2}{2}}\leq c_m,
		\ -1<u<1,
	\end{equation}
	where the inequality follows from $m\geq2$.
	Since $\delta=o(a)$, the interval
	$[-r_{ij}-\delta,r_{ij}+\delta]$ contains at most
	$C\bigl(1+\frac{r_{ij}}{a}\bigr)$ points of
	$a\bigl(\mathbb Z+\frac{1}{2}\bigr)$. Each such point contributes
	an interval of length at most $2\delta$ to $D_{ij}$. Hence, by
	\eqref{eq:wc-coordinate-density},
	\begin{align*}
		\sigma_x\{b\in S_x:(x\wedge b)_{ij}\in D_{ij}\}
		&=\int_{D_{ij}}\frac{1}{r_{ij}}
		f_T\left(\frac{u}{r_{ij}}\right)\,du\\
		&\leq \frac{c_m}{r_{ij}}\,|D_{ij}|\\
		&\leq C_m\frac{H}{r_{ij}}
		\left(1+\frac{r_{ij}}{a}\right)\\
		&\leq C_m\frac{H}{a}.
	\end{align*}
	This proves the claim.
	By \eqref{eq:wc-grid-inclusion}, \eqref{eq:wc-v-b-law} and \eqref{eq:wc-coordinate-grid-bound},
	for every Borel set $I\subset(0,H]$,
	\begin{equation}\label{eq:wc-grid-measure}
		\begin{aligned}
			\nu\{E\in U_{2H}:v(E)\in I\}
			&\leq
			\sum_{i<j}\sigma_x\{b:(x\wedge b)_{ij}\in D_{ij}\}
			\int_I mv^{m-1}\,dv\\
			&\leq C_m\frac{H}{a}\int_I mv^{m-1}\,dv.
		\end{aligned}
	\end{equation}
	Set $L>m$. By \eqref{eq:wc-grid-measure} on $(0,H]$,
	\eqref{eq:wc-v-distribution} on $(H,1)$ and setting $v=hw$,
	we obtain
	\begin{align*}
		&h^{-m}\int_{U_{2H}}
		\left(1+\frac{s_E(x)}{h}\right)^{-L}\,d\nu(E)\\
		&\leq C_m\frac{H}{a}
		\int_0^{H/h}\frac{w^{m-1}}{(1+w)^L}\,dw
		+m\int_{H/h}^{1/h}\frac{w^{m-1}}{(1+w)^L}\,dw\\
		&\leq C_m\frac{H}{a}
		\int_0^\infty\frac{w^{m-1}}{(1+w)^L}\,dw
		+m\int_{H/h}^\infty w^{m-1-L}\,dw\\
		&\leq C_{m,L}
		\left(\frac{H}{a}+\left(\frac{h}{H}\right)^{L-m}\right).
	\end{align*}
	The bound is uniform in $x\in S^d$.
\end{proof}

Let $\{f_j\}_{j=1}^{J_n}$ be the orthonormal system constructed
in Proposition~\ref{prop:packet-assembly}.
This system satisfies the estimates in
Proposition~\ref{prop:assembled-packet-estimates}.
For $T\in\mathcal T$, write
\[
I_T=\{j\in\{1,\ldots,J_n\}:T_j=T\}.
\]
Set $\mathcal P=\{\{A,-A\}:A\in\mathcal A\}$.
For each $\mathfrak a\in\mathcal P$, choose a rotation
$g_{\mathfrak a}$ and a member $A_{\mathfrak a}\in\mathfrak a$.
Write \(I_{\mathfrak a}
=\bigcup\limits_{T\in\mathcal T_{\mathfrak a}} I_T\),
where
\[
\mathcal T_{\mathfrak a}
=\{T\in\mathcal T:I_T\neq\varnothing,\
g_{\mathfrak a}c_T\in A_{\mathfrak a}\}.
\]
For a rotation $g$ and a function $f$ on $S^d$, write
$(gf)(x)=f(g^{-1}x)$.
For $j\in I_{\mathfrak a}$, set
\(w_{\mathfrak a,j,+1}=g_{\mathfrak a}f_j\) and
\(w_{\mathfrak a,j,-1}=\overline{g_{\mathfrak a}f_j}\).
Write $\{w_\iota\}_{\iota\in I}$ for this family, where
\[
I=\{(\mathfrak a,j,\varepsilon):
\mathfrak a\in\mathcal P,\ j\in I_{\mathfrak a},\
\varepsilon\in\{-1,1\}\}.
\]
For $T\in\mathcal T_{\mathfrak a}$, let
\[
\widetilde c_{\mathfrak a,T,\varepsilon}
=
\begin{cases}
	g_{\mathfrak a}c_T, & \varepsilon=+1,\\
	\overline{g_{\mathfrak a}c_T}, & \varepsilon=-1.
\end{cases}
\]
For $\iota=(\mathfrak a,j,\varepsilon)\in I$ and
$\sigma\in\{-1,1\}^m$, let
\[
\widetilde E_{\iota,\sigma}
=
\begin{cases}
	g_{\mathfrak a}g_{C_j}E_{\lambda_j,\sigma},
	& \varepsilon=+1,\\
	\overline{g_{\mathfrak a}g_{C_j}E_{\lambda_j,\sigma}},
	& \varepsilon=-1.
\end{cases}
\]

\begin{lemma}\label{lem:wc-selected-gram}
	There exists a constant $C_d>0$, independent of $n$ and of the
	choice of $(g_{\mathfrak a},A_{\mathfrak a})_{\mathfrak a\in\mathcal P}$,
	such that the following statements hold for all sufficiently large
	$n$ and every such choice.
	\begin{enumerate}
		\item[(i)] \(\begin{aligned}
			\sup\limits_{P\in Q}
			\#\bigl\{(\mathfrak a,T,\varepsilon):
			\mathfrak a\in\mathcal P,
			T\in\mathcal T_{\mathfrak a},
			\varepsilon\in\{-1,1\},
			d_Q(\widetilde c_{\mathfrak a,T,\varepsilon},P)
			\leq\ell\bigr\}
			\leq C_d.
		\end{aligned}\)
		
		\item[(ii)] \(\sup\limits_{P\in Q}
		\#\bigl\{(\iota,\sigma)\in I\times\{-1,1\}^m:
		d_Q(\widetilde E_{\iota,\sigma},P)\leq h\bigr\}
		\leq C_d\).
		
		\item[(iii)] Let $G=(\langle w_\iota,w_{\iota'}\rangle)_{\iota,\iota'\in I}$. 
		Then \(\#I\leq C_d h^{-2m}\) and 
		\[
		\|G-\mathrm{Id}\|=O(n^{-\infty}).
		\]
	\end{enumerate}
\end{lemma}
\begin{proof}
	Fix $(g_{\mathfrak a},A_{\mathfrak a})_{\mathfrak a\in\mathcal P}$.
	For $\mathfrak a\in\mathcal P$, set \(A_{\mathfrak a,+1}=A_{\mathfrak a}\) and \(A_{\mathfrak a,-1}=-A_{\mathfrak a}\).
	Then $\widetilde c_{\mathfrak a,T,\varepsilon}
	\in A_{\mathfrak a,\varepsilon}$ for every
	$T\in\mathcal T_{\mathfrak a}$. Moreover, for
	$\iota=(\mathfrak a,j,\varepsilon)\in I$, we have
	$T_j\in\mathcal T_{\mathfrak a}$. By
	Proposition~\ref{prop:packet-assembly}(ii),
	\begin{equation}\label{eq:wc-selected-plane-centre}
		d_Q\bigl(\widetilde E_{\iota,\sigma},
		\widetilde c_{\mathfrak a,T_j,\varepsilon}\bigr)
		=d_Q(g_{C_j}E_{\lambda_j,\sigma},c_{T_j})
		\leq C_d\ell.
	\end{equation}
	For $(\mathfrak a,\varepsilon)\neq
	(\mathfrak a',\varepsilon')$, the sets
	$A_{\mathfrak a,\varepsilon}$ and
	$A_{\mathfrak a',\varepsilon'}$ are distinct. Hence, for $T\in\mathcal T_{\mathfrak a}$ and
	$T'\in\mathcal T_{\mathfrak a'}$, by
	Lemma~\ref{lem:wc-geometric-estimates}(i),
	\begin{equation}\label{eq:wc-selected-centre-separation}
		d_Q\bigl(\widetilde c_{\mathfrak a,T,\varepsilon},
		\widetilde c_{\mathfrak a',T',\varepsilon'}\bigr)
		\geq c_m H.
	\end{equation}
	Thus, for
	$\iota=(\mathfrak a,j,\varepsilon)$ and
	$\iota'=(\mathfrak a',j',\varepsilon')$ with
	$(\mathfrak a,\varepsilon)\neq(\mathfrak a',\varepsilon')$,
	\eqref{eq:wc-selected-plane-centre} and
	\eqref{eq:wc-selected-centre-separation} give
	\begin{equation}\label{eq:wc-selected-plane-separation}
		d_Q(\widetilde E_{\iota,\sigma},
		\widetilde E_{\iota',\tau})
		\geq c_m H-2C_d\ell
		\geq \frac{c_mH}{2}
	\end{equation}
	for all sufficiently large $n$.
	
	\medskip
	\noindent\textup{(i)}
	Fix $P\in Q$ and set
	\[
	\mathcal I_P
	=
	\bigl\{(\mathfrak a,T,\varepsilon):
	\mathfrak a\in\mathcal P,\ T\in\mathcal T_{\mathfrak a},\
	\varepsilon\in\{-1,1\},\
	d_Q(\widetilde c_{\mathfrak a,T,\varepsilon},P)\leq\ell
	\bigr\}.
	\]
	For $(\mathfrak a,T,\varepsilon),
	(\mathfrak a',T',\varepsilon')\in\mathcal I_P$, we have
	\[
	d_Q\bigl(\widetilde c_{\mathfrak a,T,\varepsilon},
	\widetilde c_{\mathfrak a',T',\varepsilon'}\bigr)
	\leq 2\ell<c_mH
	\]
	for all sufficiently large $n$. Thus, by
	\eqref{eq:wc-selected-centre-separation},
	$(\mathfrak a,\varepsilon)=(\mathfrak a',\varepsilon')$.
	If the two triples are distinct, then $T\neq T'$ and
	Proposition~\ref{prop:packet-assembly}(iii) gives
	\begin{equation}\label{eq:wc-same-cell-centre-separation}
		d_Q\bigl(\widetilde c_{\mathfrak a,T,\varepsilon},
		\widetilde c_{\mathfrak a,T',\varepsilon}\bigr)
		=d_Q(c_T,c_{T'})
		\geq \frac{\ell}{12}.
	\end{equation}
	By the graph coordinates and rotation invariance of $\nu$,
	there exist constants $c_m,C_m,r_m>0$ such that
	\begin{equation}\label{eq:wc-small-ball-volume}
		c_m r^{2m}\leq \nu(B_Q(P,r))\leq C_m r^{2m}
	\end{equation}
	for every $P\in Q$ and $0<r\leq r_m$.
	By \eqref{eq:wc-same-cell-centre-separation}, the open balls
	$B_Q(\widetilde c_{\mathfrak a,T,\varepsilon},\frac{\ell}{24})$,
	with $(\mathfrak a,T,\varepsilon)\in\mathcal I_P$, are pairwise
	disjoint. By the definition of $\mathcal I_P$, they are contained
	in $B_Q(P,2\ell)$. Hence, by \eqref{eq:wc-small-ball-volume},
	\begin{align*}
	\#\mathcal I_P\left(\frac{\ell}{24}\right)^{2m}
	&\leq C_m
	\sum_{(\mathfrak a,T,\varepsilon)\in\mathcal I_P}
	\nu\left(B_Q\left(\widetilde c_{\mathfrak a,T,\varepsilon},
	\frac{\ell}{24}\right)\right)\\
	&\leq C_m\nu(B_Q(P,2\ell))\\
	&\leq C'_m(2\ell)^{2m}.
	\end{align*}
	Taking the supremum over $P\in Q$ proves \textup{(i)}.
	
	\medskip
	\noindent\textup{(ii)}
	Fix $P\in Q$ and set
	\[
	\mathcal J_P
	=
	\bigl\{(\iota,\sigma)\in I\times\{-1,1\}^m:
	d_Q(\widetilde E_{\iota,\sigma},P)\leq h\bigr\}.
	\]
	For $(\iota,\sigma),(\iota',\tau)\in\mathcal J_P$, we have
	\[
	d_Q(\widetilde E_{\iota,\sigma},\widetilde E_{\iota',\tau})
	\leq 2h<\frac{c_mH}{2}
	\]
	for all sufficiently large $n$. Thus, by
	\eqref{eq:wc-selected-plane-separation}, all elements of
	$\mathcal J_P$ have the same $(\mathfrak a,\varepsilon)$.
	Suppose that $\mathcal J_P\neq\varnothing$. For this
	$(\mathfrak a,\varepsilon)$, set
	\[
	P_0=
	\begin{cases}
		g_{\mathfrak a}^{-1}P, & \varepsilon=+1,\\
		g_{\mathfrak a}^{-1}\overline P, & \varepsilon=-1.
	\end{cases}
	\]
	Then, for every $((\mathfrak a,j,\varepsilon),\sigma)\in\mathcal J_P$,
	\[
	d_Q(g_{C_j}E_{\lambda_j,\sigma},P_0)
	=d_Q(\widetilde E_{(\mathfrak a,j,\varepsilon),\sigma},P)
	\leq h.
	\]
	Choose $P_1,\ldots,P_N\in\{F\in Q:d_Q(F,P_0)\leq h\}$ such that
	\begin{equation}\label{eq:wc-cover-separation}
		d_Q(P_k,P_j)\geq\frac{h}{2},\ k\neq j
	\end{equation}
	and whenever \(d_Q(F,P_0)\leq h\),
	\begin{equation}\label{eq:wc-cover-distance}
		\min_{1\leq k\leq N}d_Q(F,P_k)<\frac{h}{2}.
	\end{equation}
	By \eqref{eq:wc-cover-separation}, the balls
	$B_Q(P_k,\frac{h}{4})$ are pairwise disjoint. Since
	$d_Q(P_k,P_0)\leq h$, they are contained in $B_Q(P_0,2h)$.
	Thus, by \eqref{eq:wc-small-ball-volume},
	\[
	N\left(\frac{h}{4}\right)^{2m}\leq C_m(2h)^{2m},
	\  N\leq C'_m.
	\]
	By \eqref{eq:wc-cover-distance},
	\[
	\{F\in Q:d_Q(F,P_0)\leq h\}
	\subseteq \bigcup_{k=1}^{N}B_Q(P_k,h).
	\]
	For every $((\mathfrak a,j,\varepsilon),\sigma)\in\mathcal J_P$,
	we have $d_Q(g_{C_j}E_{\lambda_j,\sigma},P_0)\leq h$.
	Hence Proposition~\ref{prop:packet-assembly}(iv) gives
	\[
	\#\mathcal J_P
	\leq
	\sum_{k=1}^{N}
	\#\bigl\{(j,\sigma)\in I_{\mathfrak a}\times\{-1,1\}^m:
	g_{C_j}E_{\lambda_j,\sigma}\in B_Q(P_k,h)\bigr\}
	\leq C_d.
	\]
	Taking the supremum over $P\in Q$ proves \textup{(ii)}.

	\medskip
	\noindent\textup{(iii)}
	Choose $P_1,\ldots,P_{N_h}\in Q$ such that 
	\begin{align}\label{h}
	d_Q(P_\alpha,P_\beta)\geq h,\ \alpha\neq\beta,
	\end{align} 
	and
	\begin{align}\label{cover}
		Q=\bigcup_{\alpha=1}^{N_h}B_Q(P_\alpha,h).
	\end{align}
	By \eqref{h}, the balls $B_Q(P_\alpha,\frac{h}{2})$ are pairwise
	disjoint. Hence \eqref{eq:wc-small-ball-volume} gives
	$N_h\leq C_m h^{-2m}$.
	Moreover, by \eqref{cover} and \textup{(ii)},
	\begin{equation}\label{eq:wc-selected-cardinality}
		\#I
		\leq \sum_{\alpha=1}^{N_h}
		\#\bigl\{(\iota,\sigma)\in I\times\{-1,1\}^m:
		d_Q(\widetilde E_{\iota,\sigma},P_\alpha)\leq h\bigr\}
		\leq C_d h^{-2m}.
	\end{equation}
	For $\iota=(\mathfrak a,j,\varepsilon)$ and
	$\iota'=(\mathfrak a,j',\varepsilon)$, orthonormality of
	$\{f_j\}_{j=1}^{J_n}$ gives \(\langle w_\iota,w_{\iota'}\rangle=\delta_{jj'}\).
	If $(\mathfrak a,\varepsilon)\neq
	(\mathfrak a',\varepsilon')$, then, by \eqref{eq:assembled-cross} and
	\eqref{eq:wc-selected-plane-separation}, for every $K,M>0$,
	\begin{align*}
		|\langle w_\iota,w_{\iota'}\rangle|
		&\leq C_{d,K}
		\sum_{\sigma,\tau\in\{-1,1\}^m}
		\left(1+\frac{d_Q(\widetilde E_{\iota,\sigma},
			\widetilde E_{\iota',\tau})}{h}\right)^{-K}
		+C_{d,K,M}n^{-M}\\
		&\leq C_{d,K,M}
		\left(\left(\frac{h}{H}\right)^K+n^{-M}\right).
	\end{align*}
	Since $G-\mathrm{Id}$ is Hermitian, by
	\eqref{eq:wc-selected-cardinality} we have
	\begin{align*}
		\|G-\mathrm{Id}\|
		\leq
		\max_{\iota\in I}
		\sum_{\iota'\in I}
		|(G-\mathrm{Id})_{\iota,\iota'}|
		\leq C_{d,K,M}n^m
		\left(n^{-\frac{K}{4}}+n^{-M}\right).
	\end{align*}
	For any $A>0$, take $K=4(m+A)$ and $M=m+A$. Then
	\(\|G-\mathrm{Id}\|\leq C_{d,A}n^{-A}\).
\end{proof}

\begin{lemma}\label{lem:wc-diagonal-lower}
	For all sufficiently large $n$, there is a choice of
	$(g_{\mathfrak a},A_{\mathfrak a})_{\mathfrak a\in\mathcal P}$
	such that
	\begin{equation}\label{eq:wc-selected-diagonal}
		\inf_{x\in S^d}\sum_{\iota\in I}|w_\iota(x)|^2
		\geq \left(1-C_d n^{-\frac{1}{20}}\right)D.
	\end{equation}
\end{lemma}
\begin{proof}
	Put $\beta=\frac{1}{20}$. Since $m=d-1\geq2$, we have \(\beta<\min\left\{\frac{9}{100},\frac{3}{16},\frac{m}{32}\right\}\).
	For each $\mathfrak a\in\mathcal P$, choose $g_{\mathfrak a}$
	according to $\mu$ and $A_{\mathfrak a}$ uniformly from $\mathfrak a$,
	with all choices independent. Set
	\[
	X_{\mathfrak a}(x)
	=\frac{2}{D}\sum_{j\in I_{\mathfrak a}}
	|(g_{\mathfrak a}f_j)(x)|^2,
	\]
	and \[
	S(x)=\sum_{\iota\in I}|w_\iota(x)|^2
	=D\sum_{\mathfrak a\in\mathcal P}X_{\mathfrak a}(x).
	\]
	By \eqref{eq:assembled-pointwise} with exponent $\frac{L}{2}$ and the
	Cauchy--Schwarz inequality, for every $L>0$,
	\begin{equation}\label{eq:wc-rotated-square}
		|(gf_j)(x)|^2
		\leq C_{d,L}h^{-m}
		\sum_{\sigma\in\{-1,1\}^m}
		\left(1+\frac{s_{gg_{C_j}E_{\lambda_j,\sigma}}(x)}{h}
		\right)^{-L}
		+O(n^{-\infty}),
	\end{equation}
	uniformly in $g,j,x$.
	
	For $\mathfrak a=\{A,-A\}$ and fixed $g$, we have
	\[
	\mathbb E_{A_{\mathfrak a}}
	\mathbf{1}_{A_{\mathfrak a}}(gc_{T_j})
	=\frac{1}{2}\mathbf{1}_{A\cup(-A)}(gc_{T_j}).
	\]
	Since the sets $A\cup(-A)$, with $\{A,-A\}\in\mathcal P$,
	partition $Q\setminus U_H$, we have
	\begin{align*}
		D\sum_{\mathfrak a\in\mathcal P}\mathbb E X_{\mathfrak a}(x)
		&=\sum_{j=1}^{J_n}
		\int_{\{g:gc_{T_j}\notin U_H\}}
		|(gf_j)(x)|^2\,d\mu(g)\\
		&=J_n-\sum_{j=1}^{J_n}
		\int_{\{g:gc_{T_j}\in U_H\}}
		|(gf_j)(x)|^2\,d\mu(g),
	\end{align*}
	where the last equality follows from
	$\int |(gf_j)(x)|^2\,d\mu(g)=\|f_j\|_2^2=1$.
	By Proposition~\ref{prop:packet-assembly}(ii) and \eqref{we-wf},
	\begin{align}\label{distance}
	\bigl|\omega_{gg_{C_j}E_{\lambda_j,\sigma}}
	-\omega_{gc_{T_j}}\bigr|\leq C_d\ell.
	\end{align}
	Since $C_d\ell\leq H$ for all sufficiently large $n$, \eqref{distance} and the definition of $U_H$ give
	\begin{equation}\label{eq:wc-rotation-region-inclusion}
		\{g:gc_{T_j}\in U_H\}
		\subseteq
		\{g:gg_{C_j}E_{\lambda_j,\sigma}\in U_{2H}\},
		\ \sigma\in\{-1,1\}^m.
	\end{equation}
	For each $j$ and $\sigma$, the plane
	$gg_{C_j}E_{\lambda_j,\sigma}$ has distribution $\nu$.
	Thus, by \eqref{eq:wc-rotated-square} and
	Lemma~\ref{lem:wc-geometric-estimates}(iii) with $L=m+2$,
	\begin{align*}
		&\int_{\{g:gc_{T_j}\in U_H\}}
		|(gf_j)(x)|^2\,d\mu(g)\\
		\leq& C_d h^{-m}\int_{U_{2H}}
		\left(1+\frac{s_E(x)}{h}\right)^{-m-2}\,d\nu(E)
		+O(n^{-\infty})\\
		\leq& C_d\left(\frac{H}{a}
		+\left(\frac{h}{H}\right)^2\right)
		+O(n^{-\infty}).
	\end{align*}
	By Proposition~\ref{prop:packet-assembly}(i), we obtain
	\begin{equation}\label{eq:wc-expected-diagonal}
		0\leq
		1-\sum_{\mathfrak a\in\mathcal P}\mathbb E X_{\mathfrak a}(x)
		\leq C_d\left(
		n^{-\frac{9}{100}}+n^{-\frac{3}{16}}+n^{-\frac{1}{2}}
		\right).
	\end{equation}
	
	Fix $\mathfrak a\in\mathcal P$ and choose
	$P_{\mathfrak a}\in A_{\mathfrak a}$. By
	Lemma~\ref{lem:wc-geometric-estimates}(i) and
	\eqref{eq:wc-selected-plane-centre},
	\[
	\widetilde E_{(\mathfrak a,j,+1),\sigma}
	\in B_Q(P_{\mathfrak a},C_da),
	\ j\in I_{\mathfrak a},\ \sigma\in\{-1,1\}^m.
	\]
	For $t\geq0$, set
	\[
	N_{\mathfrak a,x}(t)=
	\#\bigl\{(j,\sigma)\in I_{\mathfrak a}\times\{-1,1\}^m:
	s_{\widetilde E_{(\mathfrak a,j,+1),\sigma}}(x)\leq th
	\bigr\}.
	\]
	If $(j,\sigma)$ is counted by $N_{\mathfrak a,x}(t)$ and
	$F\in B_Q(\widetilde E_{(\mathfrak a,j,+1),\sigma},h)$,
	then, by Lemma~\ref{lem:grassmannian-comparisons}, we have \(F\in B_Q(P_{\mathfrak a},C'_da)\) and \(s_F(x)\leq(t+C_d)h\).
	 By Lemma~\ref{lem:wc-selected-gram}(ii), for every $F\in Q$,
	 \[
	 \sum_{(j,\sigma)\in I_{\mathfrak a}\times\{-1,1\}^m}
	 \mathbf{1}_{B_Q(\widetilde E_{(\mathfrak a,j,+1),\sigma},h)}(F)
	 \leq C_d.
	 \]
	 Hence, by
	\eqref{eq:wc-small-ball-volume} and
	Lemma~\ref{lem:wc-geometric-estimates}(ii),
	\begin{equation}\label{eq:wc-local-plane-count}
		\begin{aligned}
			N_{\mathfrak a,x}(t)h^{2m}
			&\leq C_d\nu\bigl\{F\in B_Q(P_{\mathfrak a},C'_da):
			s_F(x)\leq(t+C_d)h\bigr\}\\
			&\leq C_d a^m h^m(1+t)^m.
		\end{aligned}
	\end{equation}
	By \eqref{eq:wc-rotated-square} with $L=m+2$, summing over
	the ranges
	$kh\leq s_{\widetilde E_{(\mathfrak a,j,+1),\sigma}}(x)<(k+1)h$
	gives
	\begin{align*}
		0\leq X_{\mathfrak a}(x)
		&\leq \frac{C_dh^{-m}}{D}
		\sum_{k=0}^{\infty}
		\frac{N_{\mathfrak a,x}(k+1)}{(1+k)^{m+2}}
		+O(n^{-\infty})\\
		&\leq \frac{C_da^mh^{-2m}}{D}
		\sum_{k=0}^{\infty}
		\frac{(2+k)^m}{(1+k)^{m+2}}
		+O(n^{-\infty}),
	\end{align*}
	where the last inequality follows from
	\eqref{eq:wc-local-plane-count}. 
	As $h^{-2m}\leq C_dD$, we obtain
	\begin{equation}\label{eq:wc-cell-sum-bound}
		0\leq X_{\mathfrak a}(x)\leq C_d a^m.
	\end{equation}
	
	Fix $x\in S^d$ and set \(Y_{\mathfrak a}(x)=X_{\mathfrak a}(x)
	-\mathbb E X_{\mathfrak a}(x)\).
	By \eqref{eq:wc-cell-sum-bound} and
	\eqref{eq:wc-expected-diagonal}, we have \(|Y_{\mathfrak a}(x)|\leq C_da^m\) and
	\[
	\sum_{\mathfrak a\in\mathcal P}
	\mathbb E|Y_{\mathfrak a}(x)|^2
	\leq C_da^m\sum_{\mathfrak a\in\mathcal P}
	\mathbb E X_{\mathfrak a}(x)
	\leq C_da^m.
	\]
	For $t=n^{-\beta}$, Bernstein's inequality gives
	\begin{equation}\label{eq:wc-diagonal-tail}
		\mathbb P\left\{
		\sum_{\mathfrak a\in\mathcal P}Y_{\mathfrak a}(x)<-t
		\right\}
		\leq \exp\left(-c_d\frac{t^2}{a^m}\right)
		=\exp\left(-c_dn^{\frac{m}{16}-2\beta}\right).
	\end{equation}
	Choose $\mathcal N\subset S^d$ such that
	$\#\mathcal N\leq C_dn^{3d}$ and
	\[
	\min_{y\in\mathcal N}d_{S^d}(x,y)\leq n^{-3},
	\ x\in S^d,
	\]
	where $d_{S^d}$ denotes the geodesic distance. Since
	$\frac{m}{16}-2\beta>0$, by \eqref{eq:wc-diagonal-tail},
	\[
	\mathbb P\left\{
	\sum_{\mathfrak a\in\mathcal P}Y_{\mathfrak a}(y)<-t
	\text{ for some }y\in\mathcal N
	\right\}
	\leq C_dn^{3d}\exp\left(-c_dn^{\frac{m}{16}-2\beta}\right)
	<1
	\]
	for all sufficiently large $n$. Fix a choice for which
	$\sum\limits_{\mathfrak a\in\mathcal P}Y_{\mathfrak a}(y)\geq-t$
	for every $y\in\mathcal N$. By \eqref{eq:wc-expected-diagonal},
	\begin{equation}\label{eq:wc-diagonal-on-net}
		S(y)\geq (1-C_dn^{-\beta})D,
		\  y\in\mathcal N.
	\end{equation}
	
	Let $k_x$ be the evaluation vector at $x$ in
	$\mathcal H_n(S^d)\otimes_{\mathbb R}\mathbb C$. By rotation
	invariance, $\|k_x\|_2^2=D$. Hence, by
	Lemma~\ref{lem:wc-selected-gram}(iii),
	\[
	S(x)=\sum_{\iota\in I}|\langle w_\iota,k_x\rangle|^2
	\leq \|G\|\|k_x\|_2^2\leq 2D.
	\]
	Since the restriction of $S$ to each great circle is a
	trigonometric polynomial of degree at most $2n$, Bernstein's
	inequality gives
	\[
	\|\nabla S\|_\infty\leq 2n\|S\|_\infty\leq 4nD.
	\]
	For every $x\in S^d$, choose $y\in\mathcal N$ such that
	$d_{S^d}(x,y)\leq n^{-3}$. By \eqref{eq:wc-diagonal-on-net},
	\[
	S(x)\geq S(y)-4n^{-2}D
	\geq (1-C'_dn^{-\beta})D.
	\]
	Taking the infimum over $x\in S^d$ gives
	\eqref{eq:wc-selected-diagonal}.
\end{proof}
\begin{lemma}\label{lem:wilson-homogenization}
	For all sufficiently large $n$, there exist a choice of
	$(g_{\mathfrak a},A_{\mathfrak a})_{\mathfrak a\in\mathcal P}$
	and an orthonormal family $\{u_\iota\}_{\iota\in I}$ in
	$\mathcal H_n(S^d)\otimes_{\mathbb R}\mathbb C$ satisfying
	the following properties.
	\begin{enumerate}
		\item[(i)] For every $\mathfrak a\in\mathcal P$ and
		$j\in I_{\mathfrak a}$, we have \(u_{\mathfrak a,j,-1}
		=\overline{u_{\mathfrak a,j,+1}}\)
		and
		\[
		\max_{\iota\in I}\|u_\iota-w_\iota\|_\infty
		=O(n^{-\infty}).
		\]
		
		\item[(ii)] For every $x\in S^d$,
		\begin{equation}\label{yiban}
			0\leq D-\sum_{\iota\in I}|u_\iota(x)|^2
			\leq C_d n^{-\frac{1}{20}}D.
		\end{equation}
	\end{enumerate}
\end{lemma}
\begin{proof}
	Fix a choice of $(g_{\mathfrak a},A_{\mathfrak a})_
	{\mathfrak a\in\mathcal P}$ satisfying
	\eqref{eq:wc-selected-diagonal}. By
	Lemma~\ref{lem:wc-selected-gram}(iii),
	$\|G-\mathrm{Id}\|=O(n^{-\infty})$. Hence, for all sufficiently
	large $n$, the eigenvalues of $G$ belong to
	$\left[\frac{1}{2},\frac{3}{2}\right]$. Set $M=G^{-\frac{1}{2}}$.
	Then
	\begin{equation}\label{eq:wc-gram-inverse-bounds}
		\|M-\mathrm{Id}\|+\|G^{-1}-\mathrm{Id}\|
		=O(n^{-\infty}).
	\end{equation}
	For $\iota\in I$, define
	\[
	u_\iota=\sum_{\kappa\in I}M_{\iota,\kappa}w_\kappa.
	\]
	Then we have \((\langle u_\iota,u_\kappa\rangle)_{\iota,\kappa\in I}
	=MGM^*=\mathrm{Id}\).
	
	\medskip
	\noindent\textup{(i)}
	For $\iota=(\mathfrak a,j,\varepsilon)$, write
	$\iota^\ast=(\mathfrak a,j,-\varepsilon)$. Since
	$w_{\iota^\ast}=\overline{w_\iota}$, \(G_{\iota^\ast,\kappa^\ast}
	=\overline{G_{\iota,\kappa}}\).
	Since $\|G-\mathrm{Id}\|<1$, we have
	\[
	M=G^{-\frac{1}{2}}
	=\sum_{k=0}^{\infty}
	\binom{-\frac{1}{2}}{k}(G-\mathrm{Id})^k.
	\]
	Since $\binom{-\frac{1}{2}}{k}\in\mathbb R$ for every $k\geq0$,
	it follows that
	$M_{\iota^\ast,\kappa^\ast}=\overline{M_{\iota,\kappa}}$. 
	Thus
	\[
	u_{\iota^\ast}
	=\sum_{\kappa\in I}M_{\iota^\ast,\kappa^\ast}
	w_{\kappa^\ast}
	=\overline{\sum_{\kappa\in I}M_{\iota,\kappa}w_\kappa}
	=\overline{u_\iota}.
	\]
	For $x\in S^d$, set
	$\mathbf w(x)=(w_\iota(x))_{\iota\in I}$ and
	$\mathbf u(x)=(u_\iota(x))_{\iota\in I}$, viewed as column vectors.
	Let $k_x$ be the evaluation vector in
	$\mathcal H_n(S^d)\otimes_{\mathbb R}\mathbb C$.
	By the addition formula, $\|k_x\|_2^2=D$. Therefore,
	\begin{equation}\label{eq:wc-before-orthogonalization-bound}
		\|\mathbf w(x)\|_{\ell^2(I)}^2
		=\sum_{\iota\in I}|\langle w_\iota,k_x\rangle|^2
		\leq \|G\|\|k_x\|_2^2
		\leq 2D.
	\end{equation}
	Since $\mathbf u(x)=M\mathbf w(x)$, by
	\eqref{eq:wc-gram-inverse-bounds} and
	\eqref{eq:wc-before-orthogonalization-bound},
	\[
	\max_{\iota\in I}|u_\iota(x)-w_\iota(x)|
	\leq \|(M-\mathrm{Id})\mathbf w(x)\|_{\ell^2(I)}
	\leq (2D)^{\frac{1}{2}}\|M-\mathrm{Id}\|
	=O(n^{-\infty}),
	\]
	uniformly in $x\in S^d$. This proves
	\textup{(i)}.
	
	\medskip
	\noindent\textup{(ii)}
	By \eqref{eq:wc-gram-inverse-bounds} and
	\eqref{eq:wc-before-orthogonalization-bound},
	\begin{equation}\label{eq:wc-diagonal-perturbation}
		\begin{aligned}
			\left|\sum_{\iota\in I}|u_\iota(x)|^2
			-\sum_{\iota\in I}|w_\iota(x)|^2\right|
			&=\left|\mathbf w(x)^*(G^{-1}-\mathrm{Id})
			\mathbf w(x)\right|\\
			&\leq 2D\|G^{-1}-\mathrm{Id}\|
			=O(n^{-\infty}).
		\end{aligned}
	\end{equation}
	By \eqref{eq:wc-selected-diagonal},
	\eqref{eq:wc-diagonal-perturbation} and Bessel's inequality,
	\[
	0\leq D-\sum_{\iota\in I}|u_\iota(x)|^2
	\leq C_d n^{-\frac{1}{20}}D,
	\ x\in S^d.
	\]
	This proves \eqref{yiban}.
\end{proof}

\section{Construction of a real orthonormal system}

Throughout this section, $d\geq 3$, $m=d-1$, and
$D_n=\dim\mathcal H_n(S^d)$. All $L^2$ norms are taken with respect
to normalized surface measure on $S^d$.
We construct a real orthonormal system $\mathcal U_n$ with the same
complex span as the packets of Lemma~\ref{lem:wilson-homogenization}
and such that
\[
\sum_{u\in\mathcal U_n}\|u\|_\infty^2 \leq C_d D_n.
\]

For $A=[U,V]\in M_{m,2}(\mathbb R)$, write
\[
\operatorname{graph}A=\{(y,Ay):y\in\mathbb R^2\}.
\]
For $S=\operatorname{diag}(s_1,\ldots,s_m)$ with $s_i=\pm1$, put
$A^S=[U,SV]$.
We write
\[
s_A(x)=\operatorname{dist}(x,\operatorname{graph}A).
\]
We write $|A|$ for the Frobenius norm of a matrix $A$.
In the following lemmas,
finite lists and their subsets are counted with multiplicity.

\begin{lemma}\label{lem:reflected-lattices}
	There exists a full lattice $\Lambda\subset M_2(\mathbb R)$
	such that, for every diagonal sign matrix $S$ of order two,
	\[
	|\det H^S|\ge1,\quad0\ne H\in\Lambda.
	\]
\end{lemma}

\begin{proof}
	We first construct $\Lambda$. Set
	\[
	R=\mathbb Z[\sqrt2]
	=\{r+s\sqrt2:r,s\in\mathbb Z\},
	\qquad
	\tau(r+s\sqrt2)=r-s\sqrt2.
	\]
	For $a=a_1+a_2\sqrt2$ and $b=b_1+b_2\sqrt2$, where
	$a_1,a_2,b_1,b_2\in\mathbb Z$, define
	\[
	H(a,b)
	=
	\begin{pmatrix}
		a&3\tau(b)\\
		b&\tau(a)
	\end{pmatrix}
	=
	\begin{pmatrix}
		a_1+a_2\sqrt2&3b_1-3b_2\sqrt2\\
		b_1+b_2\sqrt2&a_1-a_2\sqrt2
	\end{pmatrix}.
	\]
	Thus
	\[
	\begin{aligned}
		\Lambda
		&:=\{H(a,b):a,b\in R\}\\
		&=\mathbb Z
		\begin{pmatrix}1&0\\0&1\end{pmatrix}
		+\mathbb Z
		\begin{pmatrix}\sqrt2&0\\0&-\sqrt2\end{pmatrix}+\mathbb Z
		\begin{pmatrix}0&3\\1&0\end{pmatrix}
		+\mathbb Z
		\begin{pmatrix}0&-3\sqrt2\\\sqrt2&0\end{pmatrix}.
	\end{aligned}
	\]
	The four generating matrices are linearly independent over
	$\mathbb R$. Hence $\Lambda$ is a full lattice in $M_2(\mathbb R)$.
	
	For $a=r+s\sqrt2\in R$, write
	\[
	N(a):=a\tau(a)
	=(r+s\sqrt2)(r-s\sqrt2)
	=r^2-2s^2\in\mathbb Z.
	\]
	For $S=\operatorname{diag}(s_1,s_2)$, where
	$s_1,s_2\in\{-1,1\}$, we have
	\[
	H(a,b)^S
	=
	\begin{pmatrix}
		a&3s_1\tau(b)\\
		b&s_2\tau(a)
	\end{pmatrix}.
	\]
	Hence
	\begin{equation}\label{eq:many-sign-det}
		\begin{aligned}
			\det H(a,b)^S
			=s_2a\tau(a)-3s_1b\tau(b)=s_2N(a)-3s_1N(b)\in\mathbb Z.
		\end{aligned}
	\end{equation}
	
	For $c=r+s\sqrt2\in R$, observe that
	\[
	N(c)=r^2-2s^2\equiv r^2+s^2\mod3.
	\]
	Since the only squares modulo $3$ are $0$ and $1$,
	{\color{black}
		\begin{equation}\label{eq:quadratic-norm-divisibility}
			3\mid N(c)
			\quad\Longleftrightarrow\quad 3\mid r\ \text{and}\ 3\mid s
			\quad\Longleftrightarrow\quad c\in3R.
		\end{equation}
	}

	{\color{black}Suppose first that $ab\ne0$. With
		$\mathbb N_0=\{0,1,2,\ldots\}$, write}
	\[
	a=3^ja_0,\qquad b=3^kb_0,\qquad
	j,k\in\mathbb N_0,\quad a_0,b_0\notin3R.
	\]
	{\color{black}By \eqref{eq:quadratic-norm-divisibility},}
	\[
	N(a)=3^{2j}N(a_0),\qquad
	N(b)=3^{2k}N(b_0),\qquad
	3\nmid N(a_0)N(b_0).
	\]
	For a nonzero integer $q$, let $v_3(q)$ denote the exponent
	of $3$ in $q$. {\color{black}Since $3\nmid N(a_0)N(b_0)$, we have}
	{\color{black}
		\begin{equation}\label{eq:quadratic-sign-valuations}
			v_3\bigl(s_2N(a)\bigr)=2j,\qquad
			v_3\bigl(3s_1N(b)\bigr)=2k+1.
		\end{equation}
	}
	{\color{black}By \eqref{eq:many-sign-det}, if $\det H(a,b)^S=0$, then}
	\[
	s_2N(a)=3s_1N(b),
	\]
	{\color{black}contradicting \eqref{eq:quadratic-sign-valuations}, since it would
		give $2j=2k+1$.}
	
	If $b=0$ and $a\ne0$, then
	$\det H(a,b)^S=s_2N(a)\ne0$.
	Similarly, if $a=0$ and $b\ne0$, then
	$\det H(a,b)^S=-3s_1N(b)\ne0$.
	Therefore, by \eqref{eq:many-sign-det},
	\[
	\det H(a,b)^S\in\mathbb Z\setminus\{0\},
	\qquad (a,b)\ne(0,0),
	\]
	which implies
	$
	|\det H(a,b)^S|\ge1. 
	$
\end{proof}
\begin{lemma}\label{lem:reflected-cubic-lattice}
	There exists a full lattice $\Lambda_3\subset M_3(\mathbb R)$
	with the following property. For every diagonal sign matrix $S$ of order
	three, there is a diagonal sign matrix $S'$, depending only on $S$, such that
	\begin{equation}\label{eq:s4-combined-signed-volume}
		|\det[U,SV,S'W]|\ge1,
		\quad0\ne[U,V,W]\in\Lambda_3.
	\end{equation}
\end{lemma}
\begin{proof}
	We construct $\Lambda_3$; cf.\
	\cite[Proposition~6.7.4]{WM}. Set \(\zeta=e^{\frac{2\pi i}{7}}\)and \(\theta=\zeta+\zeta^{-1}\).
	Since $\zeta^7=1$ and $\zeta\ne1$,
	\[
	1+\zeta+\cdots+\zeta^6
	=\frac{\zeta^7-1}{\zeta-1}=0.
	\]
	Multiplying by $\zeta^{-3}$, we obtain
	\[
	1+(\zeta+\zeta^{-1})
	+(\zeta^2+\zeta^{-2})
	+(\zeta^3+\zeta^{-3})=0.
	\]
	On the other hand,
	\begin{gather*}
		\zeta^2+\zeta^{-2}
		=(\zeta+\zeta^{-1})^2-2
		=\theta^2-2,\\
		\zeta^3+\zeta^{-3}
		=(\zeta+\zeta^{-1})^3-3(\zeta+\zeta^{-1})
		=\theta^3-3\theta.
	\end{gather*}
	Consequently,
	\[
	0=1+\theta+(\theta^2-2)+(\theta^3-3\theta)
	=\theta^3+\theta^2-2\theta-1.
	\]
	{\color{black}Set
		\begin{equation}\label{eq:cubic-defining-polynomial}
			f(X)=X^3+X^2-2X-1,\qquad f(\theta)=0.
		\end{equation}
	}
	
	Let $\mathbb F_2=\mathbb Z/2\mathbb Z$, and let $\overline f$
	denote the polynomial obtained by reducing the coefficients of $f$
	modulo $2$. {\color{black}Then}
	\[
	\overline f(X)=X^3+X^2+1\in\mathbb F_2[X].
	\]
	Since
	$
	\overline f(0)=1$ and $\overline f(1)=1,
	$
	{\color{black}the cubic polynomial $\overline f$ is irreducible over $\mathbb F_2$.
		Since $f$ is monic, Gauss's lemma implies that $f$ is irreducible
		over $\mathbb Q$. By \eqref{eq:cubic-defining-polynomial}, it is
		the minimal polynomial of $\theta$. Set $R_3=\mathbb Z[\theta]$. Then}
	{\color{black}
		\begin{equation}\label{eq:cubic-quotient-field}
			\begin{aligned}
				R_3/2R_3
				&=\mathbb Z[\theta]/2\mathbb Z[\theta]\cong\mathbb F_2[X]/(X^3+X^2+1)
				\cong\mathbb F_8.
			\end{aligned}
		\end{equation}
	}
	Here $\mathbb F_8$ denotes the field with eight elements.
	
	Let $\tau$ now denote the automorphism of $\mathbb Q(\zeta)$
	given by $\tau(\zeta)=\zeta^2$. {\color{black}Using $\zeta^7=1$ and
		\eqref{eq:cubic-defining-polynomial}, we obtain}
	{\color{black}
		\begin{equation}\label{eq:cubic-automorphism-action}
			\tau(\theta)=\theta^2-2,\qquad
			\tau^2(\theta)=\theta^3-3\theta=1-\theta-\theta^2,\qquad
			\tau^3(\theta)=\theta.
		\end{equation}
	}
	{\color{black}By \eqref{eq:cubic-automorphism-action}, $\tau$ restricts to an
		automorphism of $R_3$. Its restriction}
	has order three, since $\theta$ has degree three and hence
	$\tau(\theta)-\theta=\theta^2-\theta-2\ne0$.
	
	For $a,b,c\in R_3$, define
	\[
	H(a,b,c)=
	\begin{pmatrix}
		a&2\tau(c)&2\tau^2(b)\\
		b&\tau(a)&2\tau^2(c)\\
		c&\tau(b)&\tau^2(a)
	\end{pmatrix},
	\]
	and set
	\[
	\Lambda_3=\{H(a,b,c):a,b,c\in R_3\}.
	\]
	
	Put $u_i=\tau^i(\theta)$, $0\le i\le2$.
	{\color{black}By \eqref{eq:cubic-automorphism-action} and the fact that $\tau$
		has order three on $R_3$, the $u_i$ are real and distinct. Hence}
	{\color{black}
		\begin{equation}\label{eq:cubic-vandermonde}
			V=(u_i^j)_{i,j=0}^2,\qquad
			\det V=\prod_{0\le i<j\le2}(u_j-u_i)\ne0.
		\end{equation}
	}
	Since $f$ is monic and is the minimal polynomial of $\theta$,
	each $a\in R_3$ has a unique expression
	$a=r_0+r_1\theta+r_2\theta^2$, with $r_0,r_1,r_2\in\mathbb Z$.
	For this expression,
	\[
	\begin{pmatrix}
		a\\ \tau(a)\\ \tau^2(a)
	\end{pmatrix}
	=V\begin{pmatrix}r_0\\r_1\\r_2\end{pmatrix}.
	\]
	Thus the image of $R_3$ under its three real embeddings is
	$V\mathbb Z^3$, a full lattice in $\mathbb R^3$.
	
	Moreover, by the definition of $H$,
	\[
	\Lambda_3
	=\sum_{j=0}^2
	\bigl(
	\mathbb ZH(\theta^j,0,0)
	+\mathbb ZH(0,\theta^j,0)
	+\mathbb ZH(0,0,\theta^j)
	\bigr).
	\]
	The three groups of matrices have disjoint supports.
	{\color{black}Within each group, linear independence follows from
		\eqref{eq:cubic-vandermonde}.}
	Hence these nine matrices form a real basis of $M_3(\mathbb R)$,
	and $\Lambda_3$ is a full lattice.
	
	For $u\in R_3$, we now use $N(u)$ and $\operatorname{Tr}(u)$ for
	its norm and trace from $\mathbb Q(\theta)$ to $\mathbb Q$, respectively:
	\[
	N(u)=u\tau(u)\tau^2(u),\qquad
	\operatorname{Tr}(u)=u+\tau(u)+\tau^2(u).
	\]
	Since the elements of $R_3$ are algebraic integers,
	$N(u),\operatorname{Tr}(u)\in\mathbb Z$.
	Expanding the determinant, we obtain
	\[
	\begin{aligned}
		\det H(a,b,c)
		&=a\tau(a)\tau^2(a)
		+2b\tau(b)\tau^2(b)
		+4c\tau(c)\tau^2(c)\\
		&\quad-2\bigl(
		a\tau(b)\tau^2(c)
		+\tau(a)\tau^2(b)c
		+\tau^2(a)b\tau(c)
		\bigr).
	\end{aligned}
	\]
	Set
	\[
	T=\operatorname{Tr}\bigl(a\tau(b)\tau^2(c)\bigr)\in\mathbb Z.
	\]
	Then
	\begin{equation}\label{eq:sec3-cubic-determinant}
		\det H(a,b,c)=N(a)+2N(b)+4N(c)-2T\in\mathbb Z.
	\end{equation}
	
	We show that $\det H(a,b,c)=0$ implies $a=b=c=0$.
	{\color{black}Let $\overline u=u+2R_3$ for $u\in R_3$. By
		\eqref{eq:cubic-quotient-field}, $R_3/2R_3$ is a field.
		Let $\overline\tau$ denote the}
	automorphism of this field induced by $\tau$. Thus
	$\overline\tau(\overline u)=\overline{\tau(u)}$, and
	\[
	\overline{N(u)}
	=\overline u\,
	\overline\tau(\overline u)\,
	\overline\tau^{\,2}(\overline u).
	\]
	{\color{black}Since $\overline\tau$ is an automorphism of this field, the
		product is zero if and only if $\overline u=0$. Therefore,}
	{\color{black}
		\begin{equation}\label{eq:cubic-norm-divisibility}
			2\mid N(u)
			\quad\Longleftrightarrow\quad
			\overline u=0
			\quad\Longleftrightarrow\quad
			u\in2R_3.
		\end{equation}
	}
	
	Suppose that $\det H(a,b,c)=0$.
	{\color{black}Reducing \eqref{eq:sec3-cubic-determinant} modulo $2$ and using
		\eqref{eq:cubic-norm-divisibility}, we get}
	\[
	N(a)\equiv0\mod2,\qquad a=2a_1,\quad a_1\in R_3.
	\]
	Consequently,
	\[
	N(a)=8N(a_1),\qquad
	T=2\operatorname{Tr}\bigl(a_1\tau(b)\tau^2(c)\bigr)
	\in2\mathbb Z.
	\]
	{\color{black}Since $N(a)\in8\mathbb Z$ and $T\in2\mathbb Z$,
		reducing \eqref{eq:sec3-cubic-determinant} modulo $4$ yields}
	\[
	0\equiv2N(b)\mod4.
	\]
	{\color{black}By \eqref{eq:cubic-norm-divisibility}, we may write $b=2b_1$
		with $b_1\in R_3$.}
	It follows that
	\[
	N(b)=8N(b_1),\qquad
	T=4\operatorname{Tr}\bigl(a_1\tau(b_1)\tau^2(c)\bigr)
	\in4\mathbb Z.
	\]
	{\color{black}Since $N(a),N(b)\in8\mathbb Z$ and $T\in4\mathbb Z$,
		reducing \eqref{eq:sec3-cubic-determinant} modulo $8$ gives}
	$
	0\equiv4N(c)\mod8.
	$
	{\color{black}By \eqref{eq:cubic-norm-divisibility}, $c=2c_1$ for some
		$c_1\in R_3$.}
	
	By the definition of $H$,
	\[
	H(a,b,c)=2H(a_1,b_1,c_1),
	\]
	and therefore
	{\color{black}
		\begin{equation}\label{eq:cubic-descent-step}
			0=\det H(a,b,c)=8\det H(a_1,b_1,c_1).
		\end{equation}
	}
	{\color{black}Iterating \eqref{eq:cubic-descent-step} gives}
	$
	a,b,c\in2^jR_3(j\ge1).
	$
	Since $R_3=\mathbb Z\oplus\mathbb Z\theta\oplus\mathbb Z\theta^2$,
	\[
	\bigcap_{j\ge1}2^jR_3
	=\left(\bigcap_{j\ge1}2^j\mathbb Z\right)
	\oplus\left(\bigcap_{j\ge1}2^j\mathbb Z\right)\theta
	\oplus\left(\bigcap_{j\ge1}2^j\mathbb Z\right)\theta^2
	=\{0\}.
	\]
	Hence $a=b=c=0$. Consequently,
	\[
	\det H(a,b,c)\in\mathbb Z\setminus\{0\}
	\qquad\text{whenever }(a,b,c)\ne(0,0,0),
	\]
	which proves
	{\color{black}
		\begin{equation}\label{eq:cubic-unsigned-determinant}
			|\det H|\ge1\qquad(0\ne H\in\Lambda_3).
		\end{equation}
	}
	
	It remains to include the sign matrices. In this final step, write
	$H=H(a,b,c)=[U,V,W]$, where $U,V,W$ are its three column vectors.
	For $p,q,r\in\{-1,1\}$, put
	$D(p,q,r)=\operatorname{diag}(p,q,r)$.
	Since $p^2=q^2=r^2=1$, we have
	\begin{equation}\label{eq:sec3-cubic-signs}
		\begin{aligned}
			D(p,q,r)H(pa,qb,rc)
			&=
			\begin{pmatrix}
				a&2pr\tau(c)&2pq\tau^2(b)\\
				b&pq\tau(a)&2qr\tau^2(c)\\
				c&qr\tau(b)&rp\tau^2(a)
			\end{pmatrix}\\
			&=[U,D(pr,qp,rq)V,D(pq,qr,rp)W].
		\end{aligned}
	\end{equation}
	
	{\color{black}Given $S=\operatorname{diag}(s_1,s_2,s_3)$ with
		$s_1,s_2,s_3\in\{-1,1\}$, put}
	\[
	\varepsilon=s_1s_2s_3,\qquad
	p=1,\qquad q=\varepsilon s_2,\qquad r=\varepsilon s_1.
	\]
	Then
	\[
	pr=\varepsilon s_1,\qquad
	qp=\varepsilon s_2,\qquad
	rq=s_1s_2=\varepsilon s_3,
	\]
	and hence
	$
	D(pr,qp,rq)=\varepsilon S.
	$
	Define
	\[
	S'=D(pq,qr,rp)
	=\varepsilon\operatorname{diag}(s_2,s_3,s_1).
	\]
	In particular, $S'$ depends only on $S$.
	By \eqref{eq:sec3-cubic-signs},
	\[
	[U,\varepsilon SV,S'W]=D(p,q,r)H(pa,qb,rc).
	\]
	Since $|\varepsilon|=1$ and $|\det D(p,q,r)|=1$,
	\[
	\begin{aligned}
		|\det[U,SV,S'W]|
		&=|\det[U,\varepsilon SV,S'W]|\\
		&=|\det D(p,q,r)|\,|\det H(pa,qb,rc)|\\
		&=|\det H(pa,qb,rc)|.
	\end{aligned}
	\]
	If $H(a,b,c)\ne0$, then $(pa,qb,rc)\ne(0,0,0)$.
	{\color{black}By \eqref{eq:cubic-unsigned-determinant}, we obtain}
	\[
	|\det[U,SV,S'W]|=|\det H(pa,qb,rc)|\ge1.
	\]
	This proves the assertion.
\end{proof}
We fix lattices $\Lambda$ and $\Lambda_3$ as in
Lemmas~\ref{lem:reflected-lattices} and~\ref{lem:reflected-cubic-lattice},
respectively.

\begin{lemma}\label{lem:reflected-partition}
	Let $m\geq 2$, $K,J\geq 1$ and $0<\eta\leq 1$. Let $\mathcal A$ be a finite list
	in $M_{m,2}(\mathbb R)$ such that $|A|\le K$ for every entry $A$.
	Suppose that every ball of radius $\eta$ contains at most $J$ entries.
	Then $\mathcal A$ can be partitioned into at most $C\eta^{-m}$
	classes, and matrices $B_A$ can be chosen so that
	\begin{equation}\label{eq:sec3-approximating-matrices}
		|A-B_A|\le C\eta
	\end{equation}
	and, for distinct entries $A,A'$ in the same class,
	\begin{equation}\label{eq:sec3-reflected-separation}
		|(B_A^S-B_{A'}^S)y|\ge c\eta|y|,
		\ y\in\mathbb R^2
	\end{equation}
	for every diagonal sign matrix $S$. Here $c,C>0$ depend only on $m,K,J$.
\end{lemma}

\begin{proof}
	We consider two cases.
	
	\smallskip
	\noindent\textit{Case 1: $m=2k$.}
	
	Write $A$ in two-row blocks $A_1,\ldots,A_k\in M_2(\mathbb R)$,
	stacked vertically. Let $C_0$ be a covering radius for $\Lambda$.
	For each $i$, choose $H_i\in\Lambda$ such that
	\[
	|A_i/\eta-H_i|\le C_0,
	\]
	and put
	\[
	B_A=
	\begin{pmatrix}
		\eta H_1\\ \vdots\\ \eta H_k
	\end{pmatrix}.
	\]
	By the definition of the Frobenius norm,
	\[
	\begin{aligned}
		|A-B_A|^2
		&=\sum_{i=1}^k|A_i-\eta H_i|^2=\eta^2\sum_{i=1}^k|A_i/\eta-H_i|^2\le kC_0^2\eta^2.
	\end{aligned}
	\]
	Hence
	{\color{black}
		\begin{equation}\label{eq:sec3-even-approximation}
			|A-B_A|\le C_0\sqrt{k}\,\eta.
		\end{equation}
	}
	
	Put $q=\lceil\eta^{-1/2}\rceil$. Since $0<\eta\le1$,
	$\eta^{-1/2}\le q\le2\eta^{-1/2}$. The lattice $\Lambda$ has rank four,
	so $\Lambda/q\Lambda\cong(\mathbb Z/q\mathbb Z)^4$. Thus
	{\color{black}
		\begin{equation}\label{eq:sec3-residue-count}
			\#\bigl((\Lambda/q\Lambda)^k\bigr)
			=q^{4k}\le 2^{4k}\eta^{-2k}=2^{4k}\eta^{-m}.
		\end{equation}
	}

	For a fixed tuple $\mathbf H=(H_1,\ldots,H_k)$, let
	$\mathcal A_{\mathbf H}$ denote the list of entries assigned
	to this tuple.
	{\color{black}
		By \eqref{eq:sec3-even-approximation},
	}
	\[
	\mathcal A_{\mathbf H}
	\subset B(\eta\mathbf H,C_0\sqrt{k}\,\eta),
	\]
	where the blocks of $\eta\mathbf H$ are stacked vertically.
	
	In dimension $4k$, this ball can be covered by at most
	$C_k=(1+2C_0\sqrt{k})^{4k}$ balls of radius $\eta$.
	{\color{black}
		Each ball contains at most $J$ entries. Hence
		\begin{equation}\label{eq:sec3-even-local-count}
			\#\mathcal A_{\mathbf H}\le C_kJ.
		\end{equation}
	}
	
	For each fixed $\mathbf H$, assign distinct local numbers
	\[
	j(A)\in\{1,\ldots,\#\mathcal A_{\mathbf H}\}
	\qquad(A\in\mathcal A_{\mathbf H}).
	\]
	Classify each entry $A$ by the label
	\[
	\bigl(j(A),H_1\bmod q\Lambda,\ldots,H_k\bmod q\Lambda\bigr).
	\]
	{\color{black}
		Combining \eqref{eq:sec3-residue-count} and
		\eqref{eq:sec3-even-local-count}, the number of classes is at most
	}
	\[
	C_kJ\,q^{4k}
	\le 2^{4k}C_kJ\,\eta^{-m}.
	\]
	
	Let $A,A'$ be distinct entries in the same class, with associated
	lattice tuples $\mathbf H=(H_r)_{r=1}^k$ and
	$\mathbf H'=(H_r')_{r=1}^k$.
	By the definition of the class labels,
	\[
	j(A)=j(A'),\qquad
	H_r-H_r'\in q\Lambda,\ 1\le r\le k.
	\]
	If $\mathbf H=\mathbf H'$, the equality of the local numbers
	would imply that $A$ and $A'$ are the same entry.
	Thus $H_i\ne H_i'$ for some $i$, and we may write
	\[
	H_i-H_i'=qH,\qquad 0\ne H\in\Lambda.
	\]
	
	{\color{black}
		Fix an arbitrary diagonal sign matrix $S$ of order $m$.
		Let $S_i$ be its restriction to this two-row block, and put
	}
	\[
	M=(\eta qH)^{S_i}
	=\eta(H_i-H_i')^{S_i}.
	\]
	Let $s_1\ge s_2$ be the singular values of $M$.
	Since the sign operation preserves the Frobenius norm,
	{\color{black}
		\begin{equation}\label{eq:sec3-even-first-singular-value}
			\begin{aligned}
				s_1
				&\le|M|
				=\eta|H_i-H_i'|\le|B_A-B_{A'}|\\
				&\le|B_A-A|+|A|+|A'|+|A'-B_{A'}|\\
				&\le2C_0\sqrt{k}\,\eta+2K\le2(K+C_0\sqrt{k}).
			\end{aligned}
		\end{equation}
	}
	On the other hand, Lemma~\ref{lem:reflected-lattices} gives
	$
	|\det H^{S_i}|\ge1.
	$
	Using $q\ge\eta^{-\frac 12}$, we therefore obtain
	{\color{black}
		\begin{equation}\label{eq:sec3-even-singular-product}
			\begin{aligned}
				s_1s_2
				=|\det M|=(\eta q)^2|\det H^{S_i}|
				\ge(\eta q)^2
				\ge\eta^2\eta^{-1}
				=\eta.
			\end{aligned}
		\end{equation}
	}
	
	{\color{black}
		Combining \eqref{eq:sec3-even-first-singular-value}
		and \eqref{eq:sec3-even-singular-product}, we obtain
	}
	\[
	s_2=\frac{s_1s_2}{s_1}
	\ge\frac{\eta}{2(K+C_0\sqrt{k})}
	=:c\eta.
	\]
	For each $r$, let $S_r$ be the restriction of $S$ to the
	$r$th two-row block. Then, for every $y\in\mathbb R^2$,
	\[
	\begin{aligned}
		|(B_A^S-B_{A'}^S)y|^2
		&=\sum_{r=1}^k
		|\eta(H_r-H_r')^{S_r}y|^2\ge|\eta(H_i-H_i')^{S_i}y|^2\\
		&=|My|^2\ge s_2^2|y|^2
		\ge c^2\eta^2|y|^2.
	\end{aligned}
	\]
	Taking square roots proves
	\eqref{eq:sec3-reflected-separation}.
	
	\smallskip
	\noindent\textit{Case 2: $m=3+2b$, where $b\ge0$ is an integer.}
	
	Write $A$ in vertically stacked row blocks $A_0,A_1,\ldots,A_b$,
	where $A_0\in M_{3,2}(\mathbb R)$ and
	$A_i\in M_2(\mathbb R)$ for $1\le i\le b$.
	
	Set $t=\sqrt{\eta}$, and let $[A_0,0]$ denote the
	$3\times3$ matrix obtained by adjoining a zero column to $A_0$.
	{\color{black}
		Let $C>0$ be a common covering radius for $\Lambda_3$ and
		$\Lambda$. By their covering properties, choose
	}
	\[
	[P_0,w]\in t\Lambda_3,\qquad
	P_i\in t\Lambda, \ 1\le i\le b,
	\]
	where $P_0$ consists of the first two columns of $[P_0,w]$
	and $w$ is its third column, so that
	{\color{black}
		\begin{equation}\label{eq:sec3-odd-lattice-approximation}
			\begin{aligned}
				|[A_0,0]-[P_0,w]|
				&=t\bigl|t^{-1}[A_0,0]-t^{-1}[P_0,w]\bigr|
				\le Ct,\\
				|A_i-P_i|
				&=t\bigl|t^{-1}A_i-t^{-1}P_i\bigr|
				\le Ct,
				\  1\le i\le b.
			\end{aligned}
		\end{equation}
	}

	Put $P=(P_0,\ldots,P_b)$, with the row blocks stacked vertically.
	{\color{black}
		By \eqref{eq:sec3-odd-lattice-approximation},
		with a constant $C_m>0$ depending only on $m$,
	}
	\[
	\begin{aligned}
		|A-P|^2
		&=|A_0-P_0|^2+\sum_{i=1}^b|A_i-P_i|^2\le(b+1)C^2t^2
		\le C_m^2t^2,
	\end{aligned}
	\]
	and
	\[
	|w|^2
	\le|[A_0,0]-[P_0,w]|^2
	\le C^2t^2.
	\]
	Hence
	\[
	|A-P|\le C_mt,\qquad |w|\le Ct.
	\]
	
	Identify $M_{m,2}(\mathbb R)$ with $\mathbb R^{2m}$, and let
	$Q$ be the unique half-open cube of side $\eta$ in the
	standard grid containing $A-P$.
	Then
	\[
	A-P\in Q,\qquad A\in P+Q.
	\]
	Fix a lattice tuple
	\[
	\mathcal L=([P_0,w],P_1,\ldots,P_b)
	\]
	and a cube $Q$, and denote the corresponding list of entries by
	$\mathcal A_{\mathcal L,Q}$.
	Since a cube of side $\eta$ in $\mathbb R^{2m}$ can be covered
	by $C_m$ balls of radius $\eta$, the multiplicity assumption gives
	{\color{black}
		\begin{equation}\label{eq:sec3-odd-local-count}
			\#\mathcal A_{\mathcal L,Q}\le C_mJ.
		\end{equation}
	}
	Assign distinct local numbers
	\[
	j(A)\in\{1,\ldots,\#\mathcal A_{\mathcal L,Q}\}
	\qquad(A\in\mathcal A_{\mathcal L,Q}),
	\]
	and use $(Q,j(A))$ as the class label.
	
		It remains to count the possible cubes $Q$.
		The local numbering bound \eqref{eq:sec3-odd-local-count} will
		then give the number of classes.
	Since $|A-P|\le C_mt$, each coordinate of $A-P$ belongs
	to $[-C_mt,C_mt]$.
	There are at most $2+2C_mt/\eta$ grid intervals of length $\eta$
	meeting this interval. Therefore,
	\[
	\#\{Q\}
	\le(2+2C_mt/\eta)^{2m}
	\le C_m'(1+t/\eta)^{2m}.
	\]
	Finally, using $t=\sqrt{\eta}$ and $0<\eta\le1$, we obtain
	{\color{black}
		\begin{equation}\label{eq:sec3-odd-cube-count}
			\begin{aligned}
				\#\{Q\}
				&\le C_m'(1+\eta^{-1/2})^{2m}\le 2^{2m}C_m'\eta^{-m}
				\le C_m''\eta^{-m}.
			\end{aligned}
		\end{equation}
	}
	
	{\color{black}
		Combining \eqref{eq:sec3-odd-local-count} and
		\eqref{eq:sec3-odd-cube-count}, the number of class labels is
		at most $C_mJ\,C_m''\eta^{-m}$. We now define $B_A$ by
	}
	\[
	B_A=P+\operatorname{center}(Q).
	\]
	Since $A-P\in Q$ and $Q$ has side length $\eta$,
	\[
	|A-B_A|
	=|A-P-\operatorname{center}(Q)|
	\le\frac{\sqrt{2m}}2\,\eta.
	\]
	This proves \eqref{eq:sec3-approximating-matrices}.
	
	Let $A,A'$ be distinct entries in the same class, and use
	primes for the quantities associated with $A'$.
	The two entries have the same cube $Q$ and the same local number.
	Their lattice tuples must therefore be different.
	Write
	\[
	\Delta B=B_A-B_{A'},\qquad
	\Delta P=P-P',\qquad
	\Delta w=w-w',
	\]
	and let $\Delta P_i=P_i-P_i'$ denote the corresponding row blocks.
	Since the cubes agree,
	{\color{black}
		\begin{equation}\label{eq:sec3-common-cube-difference}
			\begin{aligned}
				\Delta B
				&=(P+\operatorname{center}(Q))
				-(P'+\operatorname{center}(Q))=P-P'=\Delta P.
			\end{aligned}
		\end{equation}
	}
	
	With $C$ as in the covering estimates for both lattices, put
	\[
	C_K=2(K+C),\qquad c_K=\min\{(2CC_K)^{-1},C_K^{-1}\}.
	\]
	
	{\color{black}
		Fix an arbitrary diagonal sign matrix $S$ of order $m$.
		Suppose first that the three-row lattice block changes. Then
	}
	\[
	[P_0,w]\ne[P_0',w'],\qquad
	[\Delta P_0,\Delta w]\in t\Lambda_3\setminus\{0\}.
	\]
	Let $S_0$ be the restriction of $S$ to the first three rows,
	and set
	$
	M=(\Delta P_0)^{S_0}.
	$
	{\color{black}
		Applying Lemma~\ref{lem:reflected-cubic-lattice} to
	}
	$t^{-1}[\Delta P_0,\Delta w]\in\Lambda_3\setminus\{0\}$,
	we obtain a diagonal sign matrix $S'$ such that
	\[
	\begin{aligned}
		|\det[M,S'\Delta w]|
		&=t^3\left|
		\det\bigl[(t^{-1}\Delta P_0)^{S_0},
		S'(t^{-1}\Delta w)\bigr]
		\right|\ge t^3.
	\end{aligned}
	\]
	
	Let $s_1\ge s_2$ be the singular values of $M$.
	By the volume inequality and the orthogonality of $S'$,
		\begin{equation}\label{eq:sec3-cubic-singular-volume}
			\begin{aligned}
				t^3
				\le|\det[M,S'\Delta w]|\le\sqrt{\det(M^\top M)}\,|S'\Delta w|=s_1s_2|\Delta w|.
			\end{aligned}
		\end{equation}
	Moreover, since $t=\sqrt{\eta}\le1$,
		\begin{equation}\label{eq:sec3-cubic-first-singular-value}
			\begin{aligned}
				s_1
				\le|M|=|\Delta P_0|\le|P_0-A_0|+|A_0|
				+|A_0'|+|A_0'-P_0'|\le2K+2Ct
				\le C_K,
			\end{aligned}
		\end{equation}
	while
	\[
	|\Delta w|
	\le|w|+|w'|
	\le2Ct.
	\]
		By \eqref{eq:sec3-cubic-singular-volume},
		\eqref{eq:sec3-cubic-first-singular-value}, and the bound
		$|\Delta w|\le2Ct$ from \eqref{eq:sec3-odd-lattice-approximation},
	\[
	t^3\le s_1s_2|\Delta w|
	\le 2CC_Kt\,s_2.
	\]
	It follows that
	\[
	s_2\ge\frac{t^2}{2CC_K}\ge c_Kt^2.
	\]
	
	If a two-row lattice block changes, say
	$P_i\ne P_i'$ for some $1\le i\le b$, then
	$
	t^{-1}\Delta P_i\in\Lambda\setminus\{0\}.
	$
	Put $M_i=(\Delta P_i)^{S_i}$, where $S_i$ is the restriction
	of $S$ to this row block, and let $s_1\ge s_2$ now denote
	the singular values of $M_i$.
	By Lemma~\ref{lem:reflected-lattices},
	\[
	\begin{aligned}
		s_1s_2
		=|\det M_i|=t^2\left|\det(t^{-1}\Delta P_i)^{S_i}\right|\ge t^2.
	\end{aligned}
	\]
	Moreover,
	\[
	\begin{aligned}
		s_1
		&\le|M_i|=|\Delta P_i|\\
		&\le|P_i-A_i|+|A_i|
		+|A_i'|+|A_i'-P_i'|\\
		&\le2K+2Ct
		\le C_K.
	\end{aligned}
	\]
	Consequently,
	\[
	s_2
	\ge\frac{t^2}{C_K}
	\ge c_Kt^2.
	\]
	
	At least one lattice block changes.
	Thus, for some $i\in\{0,\ldots,b\}$,
	{\color{black}
		\begin{equation}\label{eq:sec3-odd-common-separation}
			|(\Delta P_i)^{S_i}y|
			\ge c_Kt^2|y|,
			\ y\in\mathbb R^2.
		\end{equation}
	}
	{\color{black}
		Writing $S_r$ for the restriction of $S$ to the $r$th row
		block, we combine \eqref{eq:sec3-common-cube-difference}
		and \eqref{eq:sec3-odd-common-separation} to obtain
	}
	\[
	\begin{aligned}
		|(B_A^S-B_{A'}^S)y|^2
		&=\sum_{r=0}^b|(\Delta P_r)^{S_r}y|^2\ge|(\Delta P_i)^{S_i}y|^2\ge c_K^2t^4|y|^2=c_K^2\eta^2|y|^2.
	\end{aligned}
	\]
	Taking square roots proves
	\eqref{eq:sec3-reflected-separation}
	and completes the proof.
\end{proof}

\begin{lemma}\label{lem:translated-incidences}
	Let $\mathcal C$ be one of the classes in
	Lemma~\ref{lem:reflected-partition}. If $|T_0|\le K$ and
	$0<\lambda\le1$, then for every
	$x\in S^{m+1}$ and $t\ge0$,
	\begin{equation}\label{eq:odd-affine-incidences}
		\#\{A\in\mathcal C:\min_S s_{T_0+\lambda A^S}(x)
		\le t\lambda\eta\}
		\le C(1+t)^m,
	\end{equation}
	where $C$ depends only on $m,K,J$.
\end{lemma}

\begin{proof}
	Fix $S$ and write $x=(y,z)\in\mathbb R^2\oplus\mathbb R^m$.
	Put $T_A=T_0+\lambda A^S$ and
	\[
	N_t=\#\{A\in\mathcal C:s_{T_A}(x)\le t\lambda\eta\}.
	\]
	The map $[-T_A,I]$ vanishes on $\operatorname{graph}T_A$. Since
	$|T_A|\le2K$, we have
	\begin{equation}\label{eq:sec3-graph-residual}
		|z-T_Ay|\le C_1s_{T_A}(x),\qquad C_1=\sqrt{1+4K^2}.
	\end{equation}
	Let $b=1/(4(1+2K))$. We again consider two cases.
	
	\smallskip
	\noindent\textbf{Case 1: }$|y|\ge b$.
	By \eqref{eq:sec3-reflected-separation}, the points $B_A^Sy$,
	$A\in\mathcal C$, are $cb\eta$-separated in $\mathbb R^m$.
	For every entry counted by $N_t$, equations
	\eqref{eq:sec3-approximating-matrices} and
	\eqref{eq:sec3-graph-residual} give
	\[
	\begin{aligned}
		\left|B_A^Sy-\frac{z-T_0y}{\lambda}\right|
		&\le|B_A-A|\,|y|+\lambda^{-1}|z-T_Ay|\le C(1+t)\eta.
	\end{aligned}
	\]
	Put $r=\frac{cb\eta}{3}$. For distinct counted entries $A,A'$, we have
	\[
	|B_A^Sy-B_{A'}^Sy|\ge cb\eta=3r>2r.
	\]
	Thus the balls $B(B_A^Sy,r)$ are pairwise disjoint.
	{\color{black}
		The bounds \eqref{eq:sec3-approximating-matrices}
		and \eqref{eq:sec3-graph-residual} therefore place their union in
	}
	\[
	B\left(\frac{z-T_0y}{\lambda},\,C(1+t)\eta+r\right).
	\]
	Let $\omega_m$ denote the volume of the unit ball in
	$\mathbb R^m$. Comparing volumes gives
	\[
	N_t\omega_m r^m
	\le\omega_m\bigl(C(1+t)\eta+r\bigr)^m.
	\]
	Cancelling $\omega_m$ and dividing by $r^m$, we obtain
	\[
	\begin{aligned}
		N_t
		&\le\left(\frac{C(1+t)\eta+r}{r}\right)^m=\left(1+\frac{3C}{cb}(1+t)\right)^m\\
		&\le\left(1+\frac{3C}{cb}\right)^m(1+t)^m=C'(1+t)^m,
	\end{aligned}
	\]
	where $C'=(1+3C/(cb))^m$ is independent of $t$ and $\eta$.
	
	\smallskip
	\noindent\textbf{Case 2:} $|y|<b$.
	There is nothing to prove if $N_t=0$. Otherwise, for a counted entry,
		\begin{equation}\label{eq:sec3-small-y-obstruction}
			\begin{aligned}
				C_1t\lambda\eta
				&\ge |z-T_Ay|\ge\sqrt{1-|y|^2}-2K|y|
				>\sqrt{1-b^2}-2Kb>\tfrac12.
			\end{aligned}
		\end{equation}
	Fix a unit vector $e\in\mathbb R^2$. By
	\eqref{eq:sec3-approximating-matrices} and
	\eqref{eq:sec3-reflected-separation} with $S=I$, the points
	$B_Ae$, $A\in\mathcal C$, lie in $B(0,K+C)$ and are
	$c\eta$-separated. Put $r=\frac{c\eta}{3}$. The balls $B(B_Ae,r)$ are
	disjoint and lie in $B(0,K+C+r)$. Comparing volumes gives
	{\color{black}
		\begin{equation}\label{eq:sec3-class-packing}
			\#\mathcal C\,r^m\le(K+C+r)^m,\qquad
			\#\mathcal C\le C'\eta^{-m}.
		\end{equation}
	}
	{\color{black}
		Combining \eqref{eq:sec3-class-packing} with
		\eqref{eq:sec3-small-y-obstruction} and $0<\lambda\le1$,
	}
	\[
	N_t\le\#\mathcal C\le C'\eta^{-m}
	\le C'(2C_1\lambda t)^m\le C''(1+t)^m.
	\]
	The constant is independent of $S$. Summing these bounds over the
	$2^m$ choices of $S$ proves \eqref{eq:odd-affine-incidences}.
\end{proof}

\begin{proposition}\label{prop:packet-sup-sum}
	For all sufficiently large $n$, the paired system obtained by applying
	Lemma~\ref{lem:wilson-homogenization} to the construction in
	Proposition~\ref{prop:packet-assembly} spans the complexification of a
	real subspace with a real orthonormal basis $\mathcal U_n$ satisfying
	\begin{equation}\label{eq:sec3-packet-sup-sum}
		\sum_{u\in\mathcal U_n}\|u\|_\infty^2\le C_dD_n.
	\end{equation}
	The change of basis preserves the orthogonal projection onto this
	subspace and its diagonal.
\end{proposition}

\begin{proof}
	Put $\ell=n^{-\frac25}$ and $h=n^{-\frac12}$.
	Fix the choice and the paired family supplied by
	Lemma~\ref{lem:wilson-homogenization}.
	Retain the columns
	$
	u_{\mathfrak a,k,+1},
	\mathfrak a\in\mathcal P, k\in I_{\mathfrak a},
	$
	and reindex them as $(u_j)_{j\in J}$.
	Their conjugates are the columns $u_{\mathfrak a,k,-1}$.
	Since the full paired family is orthonormal,
		\begin{equation}\label{eq:sec3-retained-orthogonality}
			\langle u_i,u_j\rangle=\delta_{ij},
			\quad
			\langle u_i,\overline{u_j}\rangle=0,
			\quad i,j\in J.
		\end{equation}
	In particular, $2\#J\le D_n$.
	
	For $\mathfrak a\in\mathcal P$ and
	$T\in\mathcal T_{\mathfrak a}$, write $B=(\mathfrak a,T)$.
	We identify $B$ with the retained indices corresponding to
	$(\mathfrak a,k)$ with $k\in I_T$, and denote its center by
	\[
	P_B=\widetilde c_{\mathfrak a,T,+1}
	=g_{\mathfrak a}c_T.
	\]
	Fix an enumeration $(\sigma_r)_{r=1}^{2^m}$ of $\{-1,1\}^m$.
	For each $j$ corresponding to $(\mathfrak a,k)$, put
	\[
	E_{j,r}=\widetilde E_{(\mathfrak a,k,+1),\sigma_r},
	\  1\leq r\leq 2^m.
	\]
	By Lemma~\ref{lem:wc-selected-gram}\textup{(i)},
		\begin{equation}\label{eq:sec3-block-multiplicity}
			\sup_{P\in Q}
			\#\{B:d_Q(P_B,P)\le\ell\}\le C.
		\end{equation}
	For each block center $P_B$, choose a chart $U_{F_i}$ from
	Section~4 containing $P_B$. In this chart, write
	\[
	g_i^{-1}P_B=E_{A_{P_B}}, \quad |A_{P_B}|\leq K,
	\]
	where $K$ is independent of $n$.
	The uniform bi-Lipschitz bounds imply that a slope ball of radius
	$\ell$ maps into a ball of radius $C\ell$ in $Q$.
		Such a ball is covered by a bounded number of balls of radius $\ell$.
		Combining this covering with \eqref{eq:sec3-block-multiplicity} gives
	\[
	\sup_A
	\#\{B:P_B\text{ is assigned to chart }i,\ 
	|A_{P_B}-A|\le\ell\}\le C.
	\]
	
	Apply Lemma~\ref{lem:reflected-partition} to these slope matrices
	with $\eta=\ell$, and include the chart index in the label
	$\alpha$. Write $\alpha(B)$ for the label of $B$.
	For a center assigned to chart $i$,
	\[
	s_{P_B}(x)=s_{A_{P_B}}(g_i^{-1} x).
	\]
	Lemma~\ref{lem:translated-incidences}, with
	$(T_0,\lambda)=(0,1)$, therefore gives
	\begin{equation}\label{eq:sec3-block-incidences}
		\begin{aligned}
			\#\{\alpha\}&\le C\ell^{-m},\\
			\#\{B:\alpha(B)=\alpha,\ s_{P_B}(x)\le t\ell\}
			&\le C(1+t)^m,
			\ x\in S^d,\ t\ge0.
		\end{aligned}
	\end{equation}
	
	Fix a block $B=(\mathfrak a,T)$, and write
	$T=\Phi_C(\Omega_{A_0})$.
	Put
	\[
	O_B=(g_{\mathfrak a}g_C)^{-1}.
	\]
	If $j\in B$ corresponds to the original packet parameter
	$\lambda_k=(q_k,\kappa_k)$, set
	\[
	\xi_j=\sqrt{2\pi}\,h q_k,
	\qquad
	\upsilon_j=-\sqrt{2\pi}\,h\kappa_k.
	\]
	By the satellite formula in Section~\ref{sec:global-system},
	\begin{equation}\label{eq:sec3-transported-satellites}
		\begin{gathered}
			O_BE_{j,r}=E_{A_0+[\xi_j,S_r\upsilon_j]},
			\quad S_r=\operatorname{diag}(\sigma_r),\quad
			|\xi_j|+|\upsilon_j|\le C\ell.
		\end{gathered}
	\end{equation}
	Thus the reflected form holds in the coordinates transported
	with this block.
	
	Put
	\[
	M_j=\ell^{-1}[\xi_j,\upsilon_j],
	\qquad
	\eta=\frac h\ell=n^{-\frac{1}{10}}.
	\]
	By \eqref{eq:sec3-transported-satellites}, $|M_j|\le C$, and
	$0<\eta\le1$ for all sufficiently large $n$.
	The distinct packet parameters in this block belong to
	$(\tfrac12\mathbb Z)^m\times\mathbb N^m$.
	Consequently, the matrices $[\xi_j,\upsilon_j]$ have bounded
	multiplicity in balls of radius $h$. Since
	\[
	|M_j-M|\le\eta
	\quad\Longleftrightarrow\quad
	|[\xi_j,\upsilon_j]-\ell M|\le h,
	\]
	we obtain
	\begin{equation}\label{eq:sec3-normalized-satellite-multiplicity}
		|M_j|\le C,
		\qquad
		\sup_M\#\{j\in B:|M_j-M|\le\eta\}\le C.
	\end{equation}
	Thus we can apply
	Lemma~\ref{lem:reflected-partition} to $(M_j)_{j\in B}$.
	Write $\beta_B(j)$ for the resulting label.
	The bound for the number of labels is uniform in $B$, so the same
	label set may be used in every block.
	
	For $\widetilde x=O_Bx$, \eqref{eq:sec3-transported-satellites} and
	orthogonality give
	\[
	s_{E_{j,r}}(x)
	=s_{A_0+[\xi_j,S_r\upsilon_j]}(\widetilde x)
	=s_{A_0+\ell M_j^{S_r}}(\widetilde x).
	\]
		If $\min_r s_{E_{j,r}}(x)\le th$, then
		\[
		\begin{aligned}
			\min_S s_{A_0+\ell M_j^S}(\widetilde x)
			&\le \min_r s_{E_{j,r}}(x)\le th=t\ell\eta.
		\end{aligned}
		\]
		Apply Lemma~\ref{lem:translated-incidences} with
		$(T_0,\lambda)=(A_0,\ell)$. Together with the class count in
		Lemma~\ref{lem:reflected-partition}, this gives
		\begin{equation}\label{eq:sec3-satellite-incidences}
			\begin{gathered}
				\#\{\beta\}\le C\eta^{-m}=C(\ell/h)^m,\\
				\begin{aligned}
					&\#\{j\in B:\beta_B(j)=\beta,\ \min_r s_{E_{j,r}}(x)\le th\}\le C(1+t)^m.
				\end{aligned}
			\end{gathered}
		\end{equation}
	These bounds are uniform in $B$, $x\in S^d$, and $t\ge0$.
	
	Let $B(j)$ be the block containing $j$.
	Classify the retained indices by
	\[
	\mathfrak B_{\alpha,\beta}
	=
	\{j\in J:
	\alpha(B(j))=\alpha,\ 
	\beta_{B(j)}(j)=\beta\}.
	\]
	Let $\mathcal C$ be the collection of nonempty classes
	$\mathfrak B_{\alpha,\beta}$. By
	\eqref{eq:sec3-block-incidences} and
	\eqref{eq:sec3-satellite-incidences},
	\begin{equation}\label{eq:sec3-number-of-classes}
		\#\mathcal C
		\le \#\{\alpha\}\,\#\{\beta\}
		\le C\ell^{-m}\left(\frac{\ell}{h}\right)^m
		= Ch^{-m}.
	\end{equation}
	Fix a class $\mathfrak B=\mathfrak B_{\alpha,\beta}\in\mathcal C$.
	If $j\in\mathfrak B$ and $s_{E_{j,r}}(x)\le th$ for some $r$,
	then \eqref{eq:wc-selected-plane-centre} and
	Lemma~\ref{lem:grassmannian-comparisons} give
	\[
	\begin{aligned}
		s_{P_{B(j)}}(x)
		&\le s_{E_{j,r}}(x)
		+C\,d_Q(E_{j,r},P_{B(j)})\le th+C\ell.
	\end{aligned}
	\]
	By \eqref{eq:sec3-block-incidences}, the number of possible blocks
	with label $\alpha$ is at most
	\[
	\begin{aligned}
		&\#\{B:\alpha(B)=\alpha,\ s_{P_B}(x)\le th+C\ell\}\\
		\le&
		C\left(1+\frac{th+C\ell}{\ell}\right)^m
		\le C'\left(1+\frac{th}{\ell}\right)^m.
	\end{aligned}
	\]
	Within each such block, \eqref{eq:sec3-satellite-incidences}
	gives at most $C(1+t)^m$ indices with label $\beta$
	whose satellites satisfy the required distance bound.
	Summing over the blocks, we obtain
	\begin{equation}\label{eq:odd-two-scale-count}
		\begin{aligned}
			\#\{j\in\mathfrak B:\min_r s_{E_{j,r}}(x)\le th\}
			&\le C\left(1+\frac{th}{\ell}\right)^m(1+t)^m\le C(1+t)^{2m},
		\end{aligned}
	\end{equation}
	where the last inequality uses $h\le\ell$.
	
	For $x\in S^d$, put
	\[
	\rho_j(x)=h^{-1}\min_r s_{E_{j,r}}(x).
	\]
	Apply \eqref{eq:assembled-pointwise} after the block rotation,
	and use Lemma~\ref{lem:wilson-homogenization}\textup{(i)}.
	For every $M,N>0$,
	\[
	\begin{aligned}
		|u_j(x)|
		&\le C_Mh^{-m/2}
		\sum_{r=1}^{2^m}
		(1+h^{-1}s_{E_{j,r}}(x))^{-M}
		+C_Nn^{-N}\\
		&\le C_Mh^{-m/2}(1+\rho_j(x))^{-M}
		+C_Nn^{-N}.
	\end{aligned}
	\]
	Here the fixed factor $2^m$ is absorbed into $C_M$.
	
	Fix $M>2m+1$ and $N>m$.
	For every integer $p\ge0$, \eqref{eq:odd-two-scale-count} gives
	\[
	\begin{aligned}
		\#\{j\in\mathfrak B:p\le\rho_j(x)<p+1\}
		&\le\#\{j\in\mathfrak B:\rho_j(x)\le p+1\}\le C(p+2)^{2m}.
	\end{aligned}
	\]
	Therefore,
	\begin{equation}\label{eq:sec3-class-sum}
		\begin{aligned}
			\sum_{j\in\mathfrak B}|u_j(x)|
			&=\sum_{p=0}^{\infty}
			\sum_{\substack{j\in\mathfrak B\\p\le\rho_j(x)<p+1}}
			|u_j(x)|\\
			&\le C_Mh^{-m/2}
			\sum_{p=0}^{\infty}
			\frac{(p+2)^{2m}}{(p+1)^M}
			+C_ND_nn^{-N}\\
			&\le Ch^{-m/2}+C_Nn^{m-N}\le Ch^{-m/2}.
		\end{aligned}
	\end{equation}
	The series converges because $M>2m+1$.
	We also used $D_n\le C_dn^m$ and
	$n^{m-N}\le1\le h^{-m/2}$.
	Thus the bound is uniform in $x$, $n$, and $\mathfrak B$.
	
	For a class $\mathfrak B=\{j_0,\ldots,j_{L-1}\}$, put
	\begin{equation}\label{eq:sec3-class-fourier}
		v_r=L^{-1/2}
		\sum_{q=0}^{L-1}e^{2\pi i rq/L}u_{j_q},
		\qquad 0\le r<L.
	\end{equation}
		By \eqref{eq:sec3-retained-orthogonality},
		\begin{equation}\label{eq:sec3-class-fourier-orthogonality}
			\begin{aligned}
				\langle v_r,v_s\rangle
				&=\frac1L\sum_{q=0}^{L-1}e^{2\pi i(r-s)q/L}
				=\delta_{rs},\qquad
				\langle v_r,\overline{v_s}\rangle&=0.
			\end{aligned}
		\end{equation}
	Moreover, \eqref{eq:sec3-class-sum} implies
		\begin{equation}\label{eq:sec3-class-fourier-bound}
			\begin{aligned}
				\|v_r\|_\infty
				&\le L^{-1/2}
				\sup_{x\in S^d}\sum_{j\in\mathfrak B}|u_j(x)|\le CL^{-1/2}h^{-m/2}.
			\end{aligned}
		\end{equation}
	
	Define the real functions
	\begin{equation}\label{eq:sec3-realification}
		g_r=\frac{v_r+\overline{v_r}}{\sqrt2},
		\qquad
		k_r=\frac{v_r-\overline{v_r}}{i\sqrt2},
		\qquad 0\le r<L.
	\end{equation}
		By \eqref{eq:sec3-class-fourier-orthogonality} and
		\eqref{eq:sec3-realification},
		\begin{equation}\label{eq:sec3-real-orthogonality}
			\langle g_r,g_s\rangle
			=\langle k_r,k_s\rangle=\delta_{rs},
			\qquad
			\langle g_r,k_s\rangle=0.
		\end{equation}
	Thus these functions form a real orthonormal basis of the real
	subspace corresponding to the class and its conjugate.
	The Fourier matrix is unitary, so
	\begin{equation}\label{eq:sec3-class-diagonal}
		\begin{aligned}
			\sum_{r=0}^{L-1}\bigl(g_r(x)^2+k_r(x)^2\bigr)
			&=2\sum_{r=0}^{L-1}|v_r(x)|^2=2\sum_{j\in\mathfrak B}|u_j(x)|^2.
		\end{aligned}
	\end{equation}
		By \eqref{eq:sec3-realification} and
		\eqref{eq:sec3-class-fourier-bound},
		\begin{equation}\label{eq:sec3-class-sup-sum}
			\begin{aligned}
				\sum_{r=0}^{L-1}
				& \bigl(\|g_r\|_\infty^2+\|k_r\|_\infty^2\bigr)\le4\sum_{r=0}^{L-1}\|v_r\|_\infty^2\\&\le4L\bigl(CL^{-1/2}h^{-m/2}\bigr)^2\le Ch^{-m}.
			\end{aligned}
		\end{equation}
	
	Take the union of these real bases over all classes.
		By \eqref{eq:sec3-retained-orthogonality}, distinct classes and their
		conjugates are mutually orthogonal. Together with
		\eqref{eq:sec3-real-orthogonality}, this shows that the union is a real
		orthonormal basis $\mathcal U_n$ of the required real subspace.
		Combining \eqref{eq:sec3-number-of-classes} and
		\eqref{eq:sec3-class-sup-sum}, we obtain
	\[
	\begin{aligned}
		\sum_{u\in\mathcal U_n}\|u\|_\infty^2
		&\le Ch^{-m}\#\mathfrak B\le Ch^{-2m}=Cn^m\le C_dD_n.
	\end{aligned}
	\]
		This proves \eqref{eq:sec3-packet-sup-sum}.
		Finally, summing \eqref{eq:sec3-class-diagonal} over the classes gives
	\[
	\sum_{u\in\mathcal U_n}|u(x)|^2
	=2\sum_{j\in J}|u_j(x)|^2
	=\sum_{\iota\in I}|u_\iota(x)|^2.
	\]
	The new basis spans the same complexified subspace as the paired
	family. Hence its orthogonal projection and its diagonal are preserved.
\end{proof}

\section{Proof of Theorem~\ref{thm:main}}

We next estimate the sum of the squared supremum norms of a basis for
the remaining subspace. An orthogonal change of basis will then give
a bound for each basis element.

\begin{lemma}\label{lem:gaussian-completion}
	Let $n\ge1$, and let $R\subset\mathcal H_n(S^d)$ be a real subspace
	with reproducing kernel $K_R$. Then $R$ has a real orthonormal basis
	$(v_j)_{j=1}^{\dim R}$ such that
	\begin{equation}\label{eq:sec3-gaussian-basis}
		\sum_j\|v_j\|_\infty^2
		\le C_d\log(n+2)\sup_{x\in S^d}K_R(x,x).
	\end{equation}
\end{lemma}

\begin{proof}
	The conclusion is immediate if $R=\{0\}$. Suppose that $r=\dim R>0$,
	fix a real orthonormal basis $(e_j)_{j=1}^r$ of $R$, and put
	\[
	F=\sum_{j=1}^r\gamma_je_j,\quad
	L=\sup_{x\in S^d}K_R(x,x),
	\]
	where $\gamma_1,\ldots,\gamma_r$ are independent standard real
	Gaussian random variables. For each $x\in S^d$, $F(x)$ is a centered
	real Gaussian and
	\begin{equation}\label{eq:sec3-gaussian-variance}
		\mathbb E|F(x)|^2=\sum_{j=1}^re_j(x)^2=K_R(x,x)\le L.
	\end{equation}
	
	Choose a $(4n)^{-1}$-net $\mathcal N$ for the geodesic distance on
	$S^d$, with
	\[
	M=|\mathcal N|\le C_d(n+1)^d,
	\qquad
	Z=\max_{y\in\mathcal N}|F(y)|.
	\]
	Fix a realization of $F$.
	Since $F\in R\subset\mathcal H_n(S^d)$, for every $x\in S^d$
	and every unit tangent vector $v$ at $x$, the function
	\[
	q(t)=F(\cos t\,x+\sin t\,v)
	\]
	is a trigonometric polynomial of degree at most $n$.
	Writing $\nabla F$ for the spherical gradient, Bernstein's
	inequality gives
	\[
	|\langle\nabla F(x),v\rangle|
	=|q'(0)|
	\le n\|q\|_\infty
	\le n\|F\|_\infty.
	\]
	Taking the supremum over $x$ and $v$, we obtain
	\begin{equation}\label{eq:sec3-gaussian-gradient}
			\|\nabla F\|_\infty\le n\|F\|_\infty.
		\end{equation}
	
	For each $x\in S^d$, choose $y\in\mathcal N$ with
	$\operatorname{dist}_{S^d}(x,y)\le(4n)^{-1}$. Then
	\[
	\begin{aligned}
		|F(x)|
		&\le |F(y)|
		+\operatorname{dist}_{S^d}(x,y)\|\nabla F\|_\infty\le Z+\frac{1}{4n}\|\nabla F\|_\infty.
	\end{aligned}
	\]
		Taking the supremum over $x$ and using
		\eqref{eq:sec3-gaussian-gradient}, we have
	\[
	\|F\|_\infty
	\le Z+\frac{1}{4n}\|\nabla F\|_\infty
	\le Z+\frac14\|F\|_\infty,
	\]
	and hence
		\begin{equation}\label{eq:sec3-gaussian-net-bound}
			\|F\|_\infty\le\frac43 Z\le2Z.
		\end{equation}

	By \eqref{eq:sec3-gaussian-variance} and the two-sided Gaussian tail
	bound, $\mathbb P(|F(x)|^2>t)\le2e^{-t/(2L)}$ for $t\ge0$.
	A union bound over the $M$ net points gives
	\begin{equation}\label{eq:sec3-gaussian-tail}
		\mathbb P(Z^2>t)
		\le\sum_{x\in\mathcal N}\mathbb P(|F(x)|^2>t)
		\le2Me^{-t/(2L)}.
	\end{equation}
		Set $t_0=2L\log(2M)$. By \eqref{eq:sec3-gaussian-net-bound}
		and \eqref{eq:sec3-gaussian-tail}, integration of the tail gives
	\begin{equation}\label{eq:sec3-gaussian-expectation}
		\begin{aligned}
			\mathbb E\|F\|_\infty^2
			&\le4\mathbb EZ^2
			=4\int_0^{\infty}\mathbb P(Z^2>t)\,dt\\
			&\le4t_0+8M\int_{t_0}^{\infty}e^{-t/(2L)}\,dt\\
			&=8L\log(2M)+16ML e^{-t_0/(2L)}\\
			&=8L\bigl(1+\log(2M)\bigr)
			\le C_dL\log(n+2).
		\end{aligned}
	\end{equation}
	
	By orthonormality,
	\[
	\|F\|_2^2=\sum_{j=1}^r\gamma_j^2,
	\qquad
	\mathbb E\|F\|_2^2
	=\sum_{j=1}^r\mathbb E\gamma_j^2=r.
	\]
	Put
	$
	U=\frac{F}{\|F\|_2}.
	$
	By rotational invariance of the Gaussian coefficient vector,
	$U$ is uniformly distributed on the unit sphere of $R$ and is
	independent of $\|F\|_2$. Since $F=\|F\|_2U$, it follows that
		\begin{equation}\label{eq:sec3-gaussian-radial-factorization}
			\begin{aligned}
				\mathbb E\|F\|_\infty^2
				&=\mathbb E\bigl(\|F\|_2^2\|U\|_\infty^2\bigr)=\bigl(\mathbb E\|F\|_2^2\bigr)
				\bigl(\mathbb E\|U\|_\infty^2\bigr)=r\,\mathbb E\|U\|_\infty^2.
			\end{aligned}
		\end{equation}
	
	Let $Q=(q_{ij})_{i,j=1}^r$ have normalized Haar distribution on
	$O(r)$, and put
	$
	v_j=\sum_{i=1}^r q_{ij}e_i
	(1\le j\le r).
	$
	Since $Q^TQ=I$,
	\[
	\langle v_i,v_j\rangle
	=\sum_{k=1}^r q_{ki}q_{kj}
	=\delta_{ij}.
	\]
	Thus $(v_j)_{j=1}^r$ is an orthonormal basis of $R$.
	Each column of $Q$ is uniformly distributed on the unit sphere of
	$\mathbb R^r$, so each $v_j$ has the same distribution as $U$.
		By \eqref{eq:sec3-gaussian-radial-factorization} and
		\eqref{eq:sec3-gaussian-expectation},
		$$
		\mathbb E\sum_{j=1}^r\|v_j\|_\infty^2
		=\sum_{j=1}^r\mathbb E\|v_j\|_\infty^2=r\,\mathbb E\|U\|_\infty^2=\mathbb E\|F\|_\infty^2
		\le C_dL\log(n+2).
		$$
		Thus some choice of $Q$ gives
		a basis satisfying \eqref{eq:sec3-gaussian-basis}.
\end{proof}

\begin{lemma}\label{lem:haar-completion}
	Suppose that a real orthonormal basis $(f_j)_{j=1}^D$ satisfies
	\begin{equation}\label{eq:sec3-average-sup}
		\frac1D\sum_{j=1}^D\|f_j\|_\infty^2\le A^2
	\end{equation}
	for some $A>0$. Its span has a real orthonormal basis
	$(Y_j)_{j=1}^D$ satisfying
	\begin{equation}\label{eq:sec3-uniform-sup}
		\max_{1\le j\le D}\|Y_j\|_\infty\le(1+2\sqrt2)A.
	\end{equation}
\end{lemma}

\begin{proof}
	Write $b_j=\|f_j\|_\infty$ and $J=\{j:b_j>A\}$. If $J$ is empty,
	the given basis already satisfies \eqref{eq:sec3-uniform-sup}.
	For each $j\in J$, choose $m_j=2^{k_j}$, where $k_j\in\mathbb N_0$,
	such that
	\begin{equation}\label{eq:sec3-dyadic-size}
		m_j\le b_j^2/A^2<2m_j.
	\end{equation}
		By \eqref{eq:sec3-average-sup} and \eqref{eq:sec3-dyadic-size},
	\[
	\sum_{j\in J}(m_j-1)
	\le A^{-2}\sum_{j\in J}b_j^2-|J|
	\le D-|J|=|J^c|.
	\]
	So we may choose sets $E_j\subset J^c$ such that
	\[
	|E_j|=m_j-1,
	\qquad
	E_i\cap E_j=\varnothing
	\quad(i,j\in J,\ i\ne j).
	\]
	For each $j\in J$, the group indexed by $\{j\}\cup E_j$
	contains exactly $m_j$ functions, with
	\[
	\|f_j\|_\infty=b_j>A,
	\qquad
	\|f_i\|_\infty=b_i\le A
	\quad(i\in E_j).
	\]
	
	Let $m=2^k$. Take
	\[
	h_0=m^{-1/2}(1,\ldots,1)^T
	\]
	as the constant column. For $1\le p\le k$ and
	$0\le q<m/2^p$, define $h_{p,q}\in\mathbb R^m$ by
		\begin{equation}\label{eq:sec3-haar-columns}
			h_{p,q}(i)=
			\begin{cases}
				2^{-p/2},
				&q2^p<i\le q2^p+2^{p-1},\\
				-2^{-p/2},
				&q2^p+2^{p-1}<i\le(q+1)2^p,\\
				0,&\text{otherwise},
			\end{cases}
			\qquad 1\le i\le m.
		\end{equation}
		By \eqref{eq:sec3-haar-columns},
	\[
	\sum_{i=1}^m h_{p,q}(i)=0,
	\qquad
	\|h_{p,q}\|_2^2=2^p\,2^{-p}=1.
	\]
	The supports of two distinct vectors $h_{p,q}$ are either disjoint
	or nested. In the latter case, the vector with larger support is
	constant on the smaller support. Hence
		\begin{equation}\label{eq:sec3-haar-orthogonality}
			\langle h_0,h_{p,q}\rangle=0,
			\qquad
			\langle h_{p,q},h_{p',q'}\rangle
			=\delta_{pp'}\delta_{qq'}.
		\end{equation}
	Since $\|h_0\|_2=1$ and
	\[
	1+\sum_{p=1}^k\frac{m}{2^p}=m,
	\]
		\eqref{eq:sec3-haar-orthogonality} shows that these vectors form
		the columns of a real orthogonal matrix.
	
	Fix $j\in J$, and write $m_j=2^{k_j}$.
	Assign $f_j$ to the constant column of the Haar matrix of order $m_j$,
	and label the functions indexed by $E_j$ as $g_{p,q}$, in correspondence
	with its nonconstant columns. Since $E_j\subset J^c$,
	\[
	\|g_{p,q}\|_\infty\le A
	\qquad
	\left(1\le p\le k_j,\quad 0\le q<m_j/2^p\right).
	\]
	Applying this matrix to the functions indexed by $\{j\}\cup E_j$
	gives
		\begin{equation}\label{eq:sec3-haar-transform}
			Y_r
			=m_j^{-1/2}f_j
			+\sum_{p=1}^{k_j}
			\sum_{q=0}^{m_j/2^p-1}h_{p,q}(r)g_{p,q},
			\qquad 1\le r\le m_j.
		\end{equation}
	
		By \eqref{eq:sec3-haar-columns}, for each fixed $p$, the supports of
		$h_{p,q}$ partition $\{1,\ldots,m_j\}$. Thus each row has exactly
		one nonzero coefficient at scale $p$, and
		\begin{equation}\label{eq:sec3-haar-row-sum}
			\sum_{q=0}^{m_j/2^p-1}|h_{p,q}(r)|=2^{-p/2}.
		\end{equation}
		Using \eqref{eq:sec3-haar-transform}, \eqref{eq:sec3-haar-row-sum},
		$\|g_{p,q}\|_\infty\le A$, and \eqref{eq:sec3-dyadic-size},
		the triangle inequality gives
		\begin{equation}\label{eq:sec3-haar-transformed-bound}
			\begin{aligned}
				\|Y_r\|_\infty
				&\le b_jm_j^{-1/2}+A\sum_{p=1}^{k_j}2^{-p/2}\le\sqrt2A+A\sum_{p=1}^{\infty}2^{-p/2}\\
				&=A\left(\sqrt2+\frac1{\sqrt2-1}\right)
				=(1+2\sqrt2)A.
			\end{aligned}
		\end{equation}
	
	If $m_j=1$, the Haar matrix is $(1)$, so
	\[
	Y_1=f_j,\qquad
	\|Y_1\|_\infty=b_j<\sqrt{2m_j}\,A=\sqrt2\,A
	\]
	by \eqref{eq:sec3-dyadic-size}.
	For every unassigned index $i$,
	\[
	i\notin\bigcup_{j\in J}(\{j\}\cup E_j)
	\quad\Longrightarrow\quad i\in J^c
	\quad\Longrightarrow\quad \|f_i\|_\infty\le A.
	\]
	
	Since the groups are pairwise disjoint, the Haar transformations,
	extended by the identity on the unassigned functions, define a real
	orthogonal change of basis. By
	\eqref{eq:sec3-haar-transformed-bound} and the definition of $J$,
	the resulting real orthonormal basis $(Y_i)_{i=1}^{D}$ satisfies
	\[
	\max_{1\le i\le D}\|Y_i\|_\infty
	\le \max\bigl\{(1+2\sqrt{2})A,A\bigr\}
	= (1+2\sqrt{2})A.
	\qedhere
	\]
\end{proof}

\begin{proof}[Proof of Theorem~\ref{thm:main}]
	Fix $d\ge3$.
	For all sufficiently large $n$, let $\mathcal U_n$
	be the real orthonormal system of
	Proposition~\ref{prop:packet-sup-sum}. Set
	\[
	W_n=\operatorname{span}_{\mathbb R}\mathcal U_n,
	\qquad R_n=\mathcal H_n(S^d)\ominus W_n.
	\]
		Then
		\begin{equation}\label{eq:sec3-main-decomposition}
			\mathcal H_n(S^d)=W_n\oplus R_n.
		\end{equation}
	
	Let $K_{R_n}$ be the reproducing kernel of $R_n$.
		By Lemma~\ref{lemma2.3}\textup{(iii)},
		\eqref{eq:sec3-main-decomposition}, and the orthonormality of
		$\mathcal U_n$,
		\begin{gather*}
			D_n=K_{\mathcal H_n(S^d)}(x,x)
			=K_{W_n}(x,x)+K_{R_n}(x,x),\\
			K_{W_n}(x,x)
			=\sum_{u\in\mathcal U_n}|u(x)|^2.
		\end{gather*}
	The change of basis in Proposition~\ref{prop:packet-sup-sum} preserves
	the projection diagonal. Thus \eqref{yiban} gives
	\begin{equation}\label{eq:sec3-complement-diagonal}
		0\le K_{R_n}(x,x)
		=D_n-\sum_{u\in\mathcal U_n}|u(x)|^2
		\le C_dn^{-1/20}D_n,
		x\in S^d.
	\end{equation}
		By \eqref{eq:sec3-gaussian-basis}, $R_n$ has a real orthonormal
		basis $\mathcal V_n$ satisfying
		\begin{equation}\label{eq:sec3-complement-sup-sum}
			\begin{aligned}
				\sum_{v\in\mathcal V_n}\|v\|_\infty^2
				\le C_d\log(n+2)\sup_{x\in S^d}K_{R_n}(x,x)
				\le C_d'n^{-1/20}\log(n+2)D_n,
			\end{aligned}
		\end{equation}
		where the second inequality uses \eqref{eq:sec3-complement-diagonal}.
		By \eqref{eq:sec3-main-decomposition},
		$\mathcal U_n\cup\mathcal V_n$ is a real orthonormal basis of
		$\mathcal H_n(S^d)$. Combining \eqref{eq:sec3-packet-sup-sum}
		and \eqref{eq:sec3-complement-sup-sum}, we obtain
	\begin{equation}\label{eq:sec3-full-average-bound}
		\begin{aligned}
			\sum_{f\in\mathcal U_n\cup\mathcal V_n}\|f\|_\infty^2
			&=\sum_{u\in\mathcal U_n}\|u\|_\infty^2
			+\sum_{v\in\mathcal V_n}\|v\|_\infty^2\\
			&\le C_dD_n+C_d'n^{-1/20}\log(n+2)D_n\\
			&\le C_d''D_n.
		\end{aligned}
	\end{equation}
	The constant $C_d''$ is independent of $n$, since
	$n^{-1/20}\log(n+2)\to0$ as $n\to\infty$.
	
	Set $A_d=\sqrt{C_d''}$ and enumerate
	$\mathcal U_n\cup\mathcal V_n$ as $(f_j)_{j=1}^{D_n}$.
	Dividing \eqref{eq:sec3-full-average-bound} by $D_n$ gives
	\[
	\frac1{D_n}\sum_{j=1}^{D_n}\|f_j\|_\infty^2\le A_d^2.
	\]
		Applying Lemma~\ref{lem:haar-completion} with $D=D_n$ and $A=A_d$,
		we obtain a real orthonormal basis $(Y_{n,j})_{j=1}^{D_n}$ such that
		\begin{equation}\label{eq:sec3-large-degree-bound}
			\max_{1\le j\le D_n}\|Y_{n,j}\|_\infty
			\le(1+2\sqrt2)A_d.
		\end{equation}
		Choose $n_0\ge1$ so that \eqref{eq:sec3-large-degree-bound} holds
		for every $n\ge n_0$.
	For each $0\le n<n_0$, choose any real orthonormal basis
	$(Y_{n,j})_{j=1}^{D_n}$ of $\mathcal H_n(S^d)$.
	Since these finitely many functions are continuous on $S^d$,
	\[
	M_d=\max_{0\le n<n_0}\max_{1\le j\le D_n}
	\|Y_{n,j}\|_\infty<\infty.
	\]
		Together with \eqref{eq:sec3-large-degree-bound},
		this gives, for $C_d=\max\{M_d,(1+2\sqrt2)A_d\}$,
		\[
		\sup_{n\ge0}\max_{1\le j\le D_n}\|Y_{n,j}\|_\infty\le C_d.
		\qedhere
		\]
\end{proof}

\medskip
\noindent\textbf{Declaration of AI Assistance.}
The authors used ChatGPT (OpenAI) to discuss mathematical arguments,
check calculations and improve the exposition. 
The authors carefully reviewed every suggestion given by AI and independently checked the mathematical content. 
The authors take full responsibility for the final manuscript.

\bigskip
\begingroup\small\parindent=0pt\relax
\textbf{Miao He}\\
Research Center for Mathematics and Interdisciplinary Sciences, Shandong University, Qingdao 266237, P. R. China\\
E-mail: \href{mailto:hm293428@163.com}{\texttt{hm293428@163.com}}\par
\medskip
\textbf{Penghui Wang}\\
School of Mathematics, Shandong University, Jinan 250100, P. R. China\\
E-mail: \href{mailto:phwang@sdu.edu.cn}{\texttt{phwang@sdu.edu.cn}}\par
\medskip
\textbf{Fan Zhang}\\
School of Mathematics, Shandong University, Jinan 250100, P. R. China\\
E-mail: \href{mailto:202211823@mail.sdu.edu.cn}{\texttt{202211823@mail.sdu.edu.cn}}\par
\medskip
\textbf{Zeyou Zhu}\\
School of Mathematical Sciences, Fudan University, Shanghai 200433, P. R. China\\
E-mail: \href{mailto:zhuzeyou@fudan.edu.cn}{\texttt{zhuzeyou@fudan.edu.cn}}\par
\endgroup

\begin{thebibliography}{99}
\bibitem{Axler2024}
S.~Axler,
\newblock \emph{Linear Algebra Done Right},
\newblock 4th ed., Undergraduate Texts in Mathematics, Springer, 2024.

\bibitem{Bourgain85}
J.~Bourgain,
\newblock Applications of the spaces of homogeneous polynomials to some
problems on the ball algebra,
\newblock \emph{Proc. Amer. Math. Soc.} \textbf{93} (1985), 277--283.

\bibitem{Bourgain}
J.~Bourgain,
\newblock On uniformly bounded bases in spaces of holomorphic functions,
\newblock \emph{Amer. J. Math.} \textbf{138} (2016), 571--584.
\href{https://arxiv.org/abs/1506.05694}{arXiv:1506.05694}.

\bibitem{BL}
N.~Burq and G.~Lebeau,
\newblock Injections de Sobolev probabilistes et applications,
\newblock \emph{Ann. Sci. \'Ec. Norm. Sup\'er.} (4) \textbf{46}
(2013), 917--962.
\href{https://doi.org/10.24033/asens.2206}{doi:10.24033/asens.2206}.

\bibitem{DaiXu2013}
F.~Dai and Y.~Xu,
\newblock \emph{Approximation Theory and Harmonic Analysis on Spheres and Balls},
\newblock Springer Monographs in Mathematics, Springer, New York, 2013.

\bibitem{Wilson}
I.~Daubechies, S.~Jaffard and J.-L.~Journ\'e,
\newblock A simple Wilson orthonormal basis with exponential decay,
\newblock \emph{SIAM J. Math. Anal.} \textbf{22} (1991), 554--572.

\bibitem{GFR}
F.~Gon\c{c}alves, J.~P.~Ferreira and M.~Radziwi\l\l,
\newblock Uniformly bounded orthonormal bases of homogeneous polynomials
on the complex sphere,
\newblock Preprint (2026),
\href{https://arxiv.org/abs/2610.14160}{arXiv:2610.14160}.

\bibitem{GrochenigLeinert}
K.~Gr\"ochenig and M.~Leinert,
\newblock Symmetry and inverse-closedness of matrix algebras and
functional calculus for infinite matrices,
\newblock \emph{Trans. Amer. Math. Soc.} \textbf{358} (2006),
no.~6, 2695--2711.

\bibitem{GrochenigWalnut}
K.~Gr\"ochenig and D.~Walnut,
\newblock A Riesz basis for Bargmann--Fock space related to sampling
and interpolation,
\newblock \emph{Ark. Mat.} \textbf{30} (1992), 283--295.

\bibitem{GrochenigZimmermann}
K.~Gr\"ochenig and G.~Zimmermann,
\newblock Hardy's theorem and the short-time Fourier transform
of Schwartz functions,
\newblock \emph{J. Lond. Math. Soc.} (2) \textbf{63} (2001), 205--214.

\bibitem{Hall2013}
B.~C.~Hall,
\newblock \emph{Quantum Theory for Mathematicians},
\newblock Graduate Texts in Mathematics \textbf{267},
Springer, New York, 2013.

\bibitem{Han2}
X.~Han,
\newblock Spherical harmonics and point configurations on the sphere,
\newblock \emph{Geom. Funct. Anal.} \textbf{36} (2026), 1559--1586.

\bibitem{HelzelTzavaras2017}
C.~Helzel and A.~E.~Tzavaras,
\newblock A kinetic model for the sedimentation of rod-like particles,
\newblock \emph{Multiscale Model. Simul.} \textbf{15} (2017), 500--536.

\bibitem{JNT}
D.~Jakobson, N.~Nadirashvili and J.~Toth,
\newblock Geometric properties of eigenfunctions,
\newblock \emph{Russian Math. Surveys} \textbf{56} (2001),
no.~6, 1085--1105.

\bibitem{Kirwin2010}
W.~D.~Kirwin,
\newblock Higher asymptotics of Laplace's approximation,
\newblock \emph{Asymptot. Anal.} \textbf{70} (2010),
no.~3--4, 231--248.

\bibitem{Rudin}
W.~Rudin,
\newblock \emph{Function Theory in the Unit Ball of $\mathbb C^n$},
\newblock Grundlehren Math. Wiss. \textbf{241},
Springer-Verlag, New York, 1980.

\bibitem{RyllWojtaszczyk}
J.~Ryll and P.~Wojtaszczyk,
\newblock On homogeneous polynomials on a complex ball,
\newblock \emph{Trans. Amer. Math. Soc.} \textbf{276} (1983), 107--116.

\bibitem{Shiffman}
B.~Shiffman,
\newblock Uniformly bounded orthonormal sections of positive line
bundles on complex manifolds,
\newblock \emph{Contemp. Math.} \textbf{644} (2015), 227--240.

\bibitem{TZ}
J.~A.~Toth and S.~Zelditch,
\newblock Riemannian manifolds with uniformly bounded eigenfunctions,
\newblock \emph{Duke Math. J.} \textbf{111} (2002), no.~1, 97--132.

\bibitem{VK}
J.~M.~VanderKam,
\newblock $L^\infty$ norms and quantum ergodicity on the sphere,
\newblock \emph{Int. Math. Res. Not.} (1997), no.~7, 329--347.
\newblock Correction, \emph{ibid.} (1998), no.~1, 65.

\bibitem{WZZ}
P.~Wang, F.~Zhang and Z.~Zhu,
\newblock Rudin's problem on uniformly bounded homogeneous orthonormal bases,
\newblock Preprint (2026),
\href{https://arxiv.org/abs/2610.29946}{arXiv:2610.29946}.

\bibitem{WM}
D.~Witte Morris,
\newblock \emph{Introduction to Arithmetic Groups},
\newblock Deductive Press, 2015.
\end{thebibliography}
\end{document}